\documentclass[11pt]{article}

\usepackage[margin=1in]{geometry}
\usepackage{amsmath,amssymb,amsthm,mathtools,mathrsfs}
\usepackage{booktabs,enumitem,microtype,etoolbox}
\setlist[enumerate]{label=(\roman*), leftmargin=*}
\usepackage[font=small,labelfont=bf,width=.85\textwidth]{caption}
\usepackage{subcaption}
\usepackage[section]{placeins}
\usepackage{needspace}
\usepackage{xcolor,tikz}
\usetikzlibrary{arrows.meta,calc,positioning,fit,backgrounds,decorations.pathreplacing}
\usepackage{aliascnt}
\definecolor{linkcyan}{RGB}{0,124,174}
\usepackage[colorlinks=true,linkcolor=linkcyan,citecolor=linkcyan,
urlcolor=linkcyan,bookmarksopen=true]{hyperref}
\usepackage[nameinlink,noabbrev]{cleveref}

\newtheorem{theorem}{Theorem}[section]
\newtheorem*{theoremA}{Theorem A}
\newcommand{\sharedtheorem}[2]{%
\newaliascnt{#1}{theorem}\newtheorem{#1}[#1]{#2}\aliascntresetthe{#1}}
\sharedtheorem{proposition}{Proposition}
\sharedtheorem{lemma}{Lemma}
\sharedtheorem{corollary}{Corollary}
\theoremstyle{definition}
\sharedtheorem{definition}{Definition}
\sharedtheorem{example}{Example}
\sharedtheorem{problem}{Problem}
\sharedtheorem{remark}{Remark}
\crefname{theorem}{Theorem}{Theorems}
\crefname{proposition}{Proposition}{Propositions}
\crefname{lemma}{Lemma}{Lemmas}
\crefname{corollary}{Corollary}{Corollaries}
\crefname{definition}{Definition}{Definitions}
\crefname{example}{Example}{Examples}
\crefname{problem}{Problem}{Problems}
\crefname{remark}{Remark}{Remarks}
\crefname{section}{Section}{Sections}
\crefname{appendix}{Appendix}{Appendices}
\crefname{subsection}{Section}{Sections}
\crefname{figure}{Figure}{Figures}
\crefname{equation}{equation}{equations}
\Crefname{equation}{Equation}{Equations}

\newcommand{\R}{\mathbb R}
\newcommand{\C}{\mathbb C}
\newcommand{\Q}{\mathbb Q}
\newcommand{\Z}{\mathbb Z}
\newcommand{\N}{\mathbb N}
\newcommand{\Qbar}{\overline{\mathbb Q}}
\newcommand{\PP}{\mathbb P}
\newcommand{\A}{\mathbb A}
\newcommand{\Reg}{\operatorname{Reg}}

\newcommand{\height}{\operatorname{H}}

\DeclarePairedDelimiter{\abs}{\lvert}{\rvert}
\DeclarePairedDelimiter{\card}{\lvert}{\rvert}
\DeclarePairedDelimiter{\norm}{\lVert}{\rVert}

\newcommand{\addendmarker}[1]{%
  \AtBeginEnvironment{#1}{\pushQED{\qed}\renewcommand{\qedsymbol}{\(\lozenge\)}}%
  \AtEndEnvironment{#1}{\popQED}}
\addendmarker{definition}
\addendmarker{example}
\addendmarker{remark}
\addendmarker{problem}

\title{Arithmetic Visibility in Polygons:\\
Curved Sight in the Art Gallery%
\thanks{The author thanks Satyan Devadoss and Joseph O'Rourke for their
MathFest 2025 workshop and for their book~\cite{DevadossORourke2025}.}}
\author{Malik Barrett}
\date{}

\begin{document}
\maketitle

\begin{abstract}
    In the classical art gallery problem, a guard sees a target when the segment joining them lies entirely within the polygon. We explore what changes when sight may bend along regular arcs on absolutely irreducible curves, but only those defined by primitive integer polynomials of total degree at most \(d\) and coefficient height at most \(H\). We call such curves admissible sight carriers.

    Arcs can dramatically improve sight. In a one-sided comb, no straight-line guard sees selected targets in two different teeth, yet a single central guard sees a chosen target in every tooth along parabolic carriers.

    Sight can also fail in two different ways. Geometrically, the polygon may block every sight arc of the allowed degree between two points. For instance, in an alternating comb, the number of guards needed for a row of targets spanning the comb grows linearly with the number of sides at each fixed degree, even for carriers with arbitrary real coefficients. Arithmetically, no admissible carrier may pass through the points at all. Even in the unit square, at fixed \(d\) and \(H\), the proportion of rational targets that must guard themselves tends to one as the height bound grows. An algebraic target may lie on no rational curve of degree at most \(d\), however large \(H\) is.

    Both failures can be overcome for finite target sets by allowing the degree and coefficient height to grow. We show that every finite set of algebraic targets in a rational polygon becomes visible from any rational interior guard, with explicit bounds. We also construct, for every \(d\), arbitrarily large Galois orbits in the unit square, all of projective height \(4^{d+1}\), whose points must each guard themselves at degree \(d\) for every \(H\). Yet at degree \(d+1\) and coefficient height one, a single guard sees the whole orbit along \(Z(Y-X^{d+1})\).
\end{abstract}


\section{Introduction}

In the classical art gallery problem a guard in a polygonal ``gallery'' sees a target when the segment joining them stays inside the polygon. One asks for the minimum number of guards needed to see the whole gallery, its \emph{guard number}. Chv\'atal's theorem, in Fisk's proof, then gives \(\lfloor n/3\rfloor\) vertex guards for a simple \(n\)-gon \cite{Chvatal1975,Fisk1978}; see also \cite[Theorem~1.1]{ORourke1987}. We ask what becomes of the guard number if sight is allowed to bend, but only along curves that can be written down with a bounded amount of arithmetic information. The question was prompted by the comb galleries of \cite[Section~1.4, Figure~1.17]{DevadossORourke2025} and by the edge-guarding question posed there as Unsolved Problem~5.

We let sight travel along a regular arc of an absolutely irreducible plane curve defined by a primitive integer polynomial of total degree at most \(d\) and largest absolute coefficient at most \(H\). We call such a curve an \emph{admissible sight carrier} and the connecting arc a \emph{sight arc}; only the arc, not the whole carrier, has to stay inside the gallery. When no confusion can arise, we simply say \emph{carrier} and \emph{arc}. Different targets may use different carriers, just as a classical guard uses different lines for different targets. Arithmetic thus enters the sight rule itself, and not only the choice of guard locations, which is where it has usually entered the art gallery literature (see \nameref{par:related}). Morally, \(d\) measures how much a sight arc may bend, while \(H\) measures how much arithmetic information it may carry. Throughout, the question is the classical one: how few guards suffice to
see a given set of targets (\cref{def:guard-number}) within a given gallery?

Our choice of rational coefficients makes the art gallery problem arithmetic. A rational polynomial vanishing at an algebraic target also vanishes at every conjugate of that target, so choosing a carrier through one point imposes conditions on its whole Galois orbit. However, a polynomial over the target field need not vanish at those conjugates; at degree one such polynomials would simply recover segment visibility (\cref{rem:coefficient-field}).

\paragraph{Sight requirements.}
Sight, in our sense, has to satisfy two requirements: an admissible polynomial must vanish at the guard and the target, and a regular arc on its real curve must join them inside the polygon. The first requirement is arithmetic and the second is geometric. 

Regarding the second requirement, comb-shaped galleries demonstrate why straight-line sight can require many guards: separate teeth obscure different parts of the gallery. In a one-sided comb we will exhibit targets such that no two can be seen by a single straight-line guard, yet a single central guard sees them all along parabolic carriers. Curved sight does not always overcome geometric obstructions, however. If the notches of the comb are deep and alternating, a sight arc between suitably placed targets must repeatedly cross a fixed line, and a carrier of degree \(d\) has at most \(d\) such intersections (see \cref{fig:alternating-comb}). At fixed degree, the guard number therefore grows linearly with the number of sides. This limitation holds even for real coefficients, with no coefficient bound; the obstruction is geometric.

Regarding the first requirement, consider the unit square. Since every segment joining two of its points stays inside it, a single straight-line guard sees the whole square. Any increase in the guard number there is thus arithmetic, and arises because the degree and coefficient height are fixed. In fact, almost all rational points, ordered by height, must then guard themselves. For algebraic targets, increasing the coefficient-height budget alone may never help. The sharpest instance of this is proved at the end of the paper (\cref{thm:main-jump}):

\begin{theoremA}
    For every \(d \ge 1\) and every \(M \ge 1\) there is a Galois-stable set \(T\) of more than \(M\) algebraic points in the interior of the unit square, all of projective height \(4^{d+1}\), such that every point of \(T\) must guard itself at degree \(d\), at every coefficient height, whereas at degree \(d+1\) and coefficient height one a single point of \(T\) sees all of \(T\).
\end{theoremA}

\paragraph{Main results and plan of the paper.}
After setting up the model (\cref{sec:model}) and the combinatorics of finite target sets (\cref{sec:finite}), our first comparisons investigate how much of classical guarding the new sight rule retains (\cref{sec:linear}). We recover the classical theory restricted to rational points when we consider degree-one carriers with unrestricted coefficient height. Similarly, for any rational polygon, degree-one edge guarding of rational targets recovers the classical minimum once the target-height bound is sufficiently large, with coefficient height linear in the target height. Edge guarding of quadratic targets turns out to be purely combinatorial. As each nonrational quadratic point lies on exactly one rational line, the degree-one edge guard number is the size of a minimum vertex cover of a graph on the polygon's edges determined by its chords, whatever the quadratic field. As with point guards, the edge guard number for a given target set can collapse with an increase in degree: in a convex \(n\)-gon, guarding all quadratic targets requires \(n-1\) edges at degree one but only one at degree two, with coefficient height unrestricted in both cases. We then use the combs described above to demonstrate both the advantages and the limitations of curved sight in arithmetic guarding.

Restricting our attention to an L-shaped gallery, we compare the effects of degree, guard placement, and coefficient field (\cref{sec:l-gallery}). Rational degree-two carriers reduce the edge guard number for quadratic targets from four to one, even in this nonconvex gallery. Furthermore, we construct rational conics joining a fixed guard outside the classical kernel to an explicit family of quadratic targets. If carrier coefficients are instead allowed in a fixed real quadratic field \(K\), the same guard sees every \(K\)-point target at degree two.

We then examine the arithmetic restrictions on their own. Increasing height can overcome the scarcity of carriers through rational targets (\cref{sec:scarcity}), but not the absence of a rational polynomial relation of an allowed degree through an algebraic target (\cref{sec:relation}). Allowing \emph{both} bounds to grow permits every algebraic target in a rational polygon to be joined to any chosen rational interior guard by the arithmetic arc lifting theorem. If the target set is finite, then this can be done for some finite budget \(d, H\). An effective refinement of arithmetic arc lifting gives coarse explicit degree and height bounds from interpolation data and estimates controlling how much the sight arc may be perturbed while remaining regular and inside the polygon (\cref{thm:lifting,thm:effective-lifting}). Mere computability would already follow by enumerating carriers and testing sight; the formulas instead quantify the required degree and height. Finally, \cref{sec:jump} proves Theorem A, which shows that a degree jump from \(d\) to \(d+1\) can reduce an arbitrarily large guard number to one, with target height independent of the orbit size. We end with questions suggested by the examples (\cref{sec:questions}).

Longer proofs are collected in the appendices; the main text is meant to be readable without them.

\paragraph[Related sight rules]{Related sight rules.}\label{par:related}
Arithmetic restrictions were previously examined in the classical art gallery problem. Requiring rational guard coordinates can increase the minimum number of guards \cite[Theorem~1]{AbrahamsenAdamaszekMiltzow2017} and \cite[Theorem~1]{MeijerMiltzow2026}. Other work restricts guard locations to a grid \cite{VaeziGhodsi2021}, or studies the computational difficulty of finding a minimum set of point guards \cite{AbrahamsenAdamaszekMiltzow2022}. In these problems, sight still follows straight segments.

Changing the sight rule has also been studied. Circular and conic visibility allow arcs from particular families of curves \cite{ChouWoo1995,GarciaLopezRamos2000}, while conic link paths allow several conic pieces and minimize the number needed to connect two points \cite{ChouChou2016}. By contrast, in our model each sight arc must lie on a single absolutely irreducible curve, though different targets may use different carriers. Rather than fixing a family of circles or conics, we bound the degree and coefficient height of its primitive integer polynomial.

There is a related arithmetic literature on visibility along polynomial graphs \cite{GoinsHarrisKubikMbirika2018,ChaubeyPandey2021,LobsenzPhillips2026}. There, sight follows specified polynomial families, and an intervening lattice point can block visibility. Our targets do not block sight. Instead, arithmetic restricts which carriers are available, and the polygon restricts which connecting arcs can be used. We ask how these restrictions change the minimum number of guards.

\section{Arithmetic visibility}\label{sec:model}

We adopt the art gallery problem convention of a \emph{polygon} \(P\) being any closed region in the real plane bounded by simple closed piecewise-linear curves. Throughout this paper, a polygon means a simple polygon with nonempty interior. A \emph{rational polygon} has vertices in \(\Q^2\). For any subfield \(K \subseteq \C\) write \(P(K) \coloneqq P \cap K^2\) to refer to the \(\mathit{K}\)\emph{-points of \(P\)}, including its boundary points.\footnotemark{} Similarly, \((\operatorname{int}P)(K)\) means the \(K\)-points in the interior of \(P\). We use \({\N=\{1,2,\ldots\}}\), and write \(H \in [1,\infty]\) for a coefficient budget, allowing \(H=\infty\) to mean that no bound is imposed.

\footnotetext{Here \(P(K)\) denotes coordinate membership in a polygonal region, rather than \(K\)-points of an algebraic variety as in \(\mathrm{C}_F(K)\). In this paper, \(K\) is either a number field or \(\Qbar\).}

We write the upright symbol \(\height\) for the absolute multiplicative Weil height and abbreviate \({\height([1:x:y])}\) to \(\height(q)\) for \(q = (x,y)\). For a reduced rational number \(a/b\), this gives \(\height(a/b) = \max\{\abs{a},\abs{b}\}\). The italic letters \(d,H\) bound the degree and coefficient height of a carrier. As such, these parameters may be referred to as the \emph{budget} for admissible carriers. We use \(B\) to bound the Weil height of a target. These are different restrictions. All degrees and coefficient heights are measured in the chosen affine plane coordinates.

We use \(O(g)\) for an upper bound by a constant multiple of \(g\), and \(\Theta(g)\) for upper and lower bounds by positive constant multiples. Subscripts specify the parameters on which these constants may depend. An absolute constant is independent of all parameters in the statement. Named polygon constructions use calligraphic letters; \(P\) denotes a general polygon and \(\mathcal{S}=[0,1]^2\) the unit square.

\begin{definition}[Admissible polynomials and their carriers]\label{def:admissible}
    For \(d \in \N\) and \(H \in [1,\infty]\), let \(\mathcal F_{d,H}\) be the set, modulo multiplication by \(-1\), of primitive absolutely irreducible polynomials \(F \in \Z[X,Y]\) such that
    \[
        1 \le \deg F \le d,\qquad \norm{F}_\infty \le H.
    \]
    Here \(\norm{F}_\infty\) is the largest absolute coefficient, and the second condition is omitted when \(H = \infty\). 

    We call \(\mathcal F_{d,H}\) the set of \emph{admissible polynomials} at budget \(d,H\). For any admissible polynomial \(F\) the affine curve \(\mathrm{C}_F \coloneqq Z(F)\) is its \emph{admissible (sight) carrier}.
\end{definition}

For a field extension \(K/\Q\), the \(K\)-points of \(\mathrm{C}_F\) are
\[
    \mathrm{C}_F(K)=\{(a,b)\in K^2:F(a,b)=0\}.
\]
An embedding \(\sigma:K\hookrightarrow K'\) sends these points to \(\mathrm{C}_F(K')\) by applying \(\sigma\) to each coordinate. A chosen real embedding places the \(K\)-points in the real plane, where we consider visibility. If a sight arc joins two such points, it lies in \(\mathrm{C}_F(\R)\), but need not consist of \(K\)-points. Different real embeddings may place the endpoints differently relative to the polygon.

If a rational polynomial vanishes at an algebraic target, it vanishes at every Galois conjugate of that target. Choosing a carrier through one point therefore imposes conditions at an entire Galois orbit. Those conjugates need not lie in the polygon and may not even be real.

Write \(\Reg(\mathrm{C}_F)\) for the regular locus, i.e., the nonsingular points of the curve. Define the \emph{regular trace} in \(P\) as
\[
    \mathrm{R}_F(P) \coloneqq \Reg(\mathrm{C}_F)(\R) \cap P
       =\{z \in P : F(z)=0,\ \nabla F(z) \ne 0 \}.
\]
The regular trace is a semialgebraic set with finitely many connected components. Regularity refers to the carrier; its intersection with \(P\) may have endpoints or isolated points.

\begin{figure}[htbp]
    \centering
    \begin{tikzpicture}[
        x=1cm,y=1cm,
        carrier/.style={draw=black!45,line width=.8pt},
        sight/.style={draw=blue!70!black,line width=1.4pt},
        blocked/.style={draw=red!75!black,line width=1.4pt,dashed},
        gallery/.style={draw=black!75,fill=black!5,line width=.8pt},
        dot/.style={circle,fill=black,inner sep=1.6pt},
        every node/.style={font=\small}
        ]
        \foreach \panel/\offset in {a/0,b/4.9}{
          \begin{scope}[shift={(\offset,0)}]
            \path[gallery]
              (0,0)--(4,0)--(4,3)--(2.5,3)--(2.5,1.6)
              --(1.5,1.6)--(1.5,3)--(0,3)--cycle;
            \node at (.25,2.7) {\(P\)};
            \node at (2,-.48) {(\panel)};
          \end{scope}
        }

        \draw[carrier,domain=-.13:4.13,samples=100]
          plot (\x,{.5+5/14*(\x-2)^2});
        \draw[sight,domain=.6:3.4,samples=100]
          plot (\x,{.5+5/14*(\x-2)^2});
        \node[dot,label=below:{\(p\)}] at (.6,1.2) {};
        \node[dot,label=below:{\(q\)}] at (3.4,1.2) {};

        \begin{scope}[shift={(4.9,0)}]
            \draw[carrier,domain=-.13:4.13,samples=100]
              plot (\x,{2.25-15/28*(\x-2)^2});
            \draw[sight,domain=.6:1.5,samples=50]
              plot (\x,{2.25-15/28*(\x-2)^2});
            \draw[blocked,domain=1.5:2.5,samples=50]
              plot (\x,{2.25-15/28*(\x-2)^2});
            \draw[sight,domain=2.5:3.4,samples=50]
              plot (\x,{2.25-15/28*(\x-2)^2});
            \node[dot,label=below:{\(p\)}] at (.6,1.2) {};
            \node[dot,label=below:{\(q\)}] at (3.4,1.2) {};
        \end{scope}
    \end{tikzpicture}

    \caption{Along the depicted carriers, \(p\) sees \(q\) in (a), but not in (b).}
    \label{fig:carrier-sight-arcs}
\end{figure}
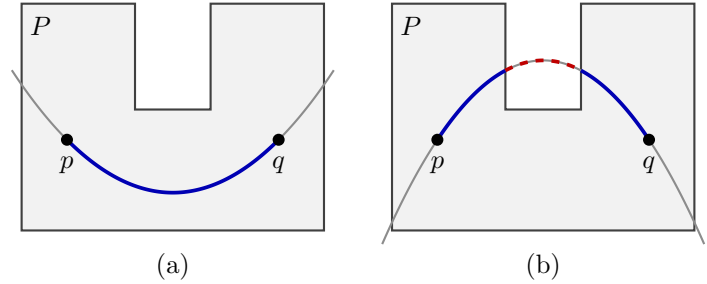

\begin{definition}[Arithmetic visibility]\label{def:visibility}
    Distinct points \(p,q \in P\) satisfy \(p \sim_{d,H} q\) if they belong to the same connected component of \(\mathrm{R}_F(P)\) for some \(F \in \mathcal F_{d,H}\). We separately declare \(p \sim_{d,H}p\) for every \(p \in P\).

    If this relation is satisfied, we say that \(p\) \emph{sees} \(q\) or that \(p\) \emph{guards} \(q\) (see \cref{fig:carrier-sight-arcs}).
\end{definition}

\begin{remark}
    Distinct points in a connected semialgebraic set can be joined by an arc. We call such an arc within one \(\mathrm{R}_F(P)\) a \emph{sight arc}. Thus, for \(p\) to see \(q\) they must lie on a common admissible carrier and must be joined by an arc of that carrier contained in \(P\). The definition prohibits a change of carriers at an intersection and passage through a singular point of that carrier. Self-visibility is separate, so a target on no admissible carrier can still guard itself.
\end{remark}

\begin{definition}[Arithmetic art gallery and guard number]\label{def:guard-number}
    An \emph{arithmetic art gallery} consists of a rational polygon \(P\), a set of targets \(T \subseteq P\), and a set of guards \(G \subseteq P\), equipped with the visibility relations \(\sim_{d,H}\) of \cref{def:visibility}. For a budget \(d,H\), a finite set \(G_0\subseteq G\) is a \emph{cover} of \(T\) if every \(t \in T\) is visible from some \(g \in G_0\). The \emph{guard number} is
    \[
        \gamma_{d,H}(P;T,G)=\min\{\card{G_0} : G_0\subseteq G \text{ is a cover of } T\},
    \]
    with value \(\infty\) if no cover exists. 
    When no confusion can arise, we call either these data or the polygon itself a \emph{gallery}.
\end{definition}

\begin{remark}\label{rem:rational-gallery}
    We assume that \(P\) is rational so that its boundary is defined over the same field as our sight carriers. Its vertices provide rational guard locations, and its edge lines become admissible carriers at sufficiently large coefficient height. These facts let us compare with classical vertex guarding and obtain explicit height bounds. Arithmetic visibility also makes sense in more general polygons, but their boundaries may introduce a further coefficient-field obstruction.

    We leave guard and target locations unrestricted so that the effects of the carrier budget can be distinguished from restrictions on the points themselves. Individual arithmetic results impose coordinate restrictions where needed.
\end{remark}

A guard need not use the same carrier for every target. Thus a bound on the number of targets lying on one carrier does not, by itself, bound the number visible from one guard. If \(T\) is finite and \(T \subseteq G\), self-visibility gives \(\gamma_{d,H}(P;T,G) \le \card{T}\).

For \(B \ge 1, D \in \N\), and \(K\) a number field, define the finite target sets
\[
    T_{K,B}(P) \coloneqq \left \{q\in P(K) : \height(q) \le B \right\},
    \qquad
    T_{D,B}(P) \coloneqq \left \{q\in P(\Qbar) : [\Q(q):\Q] \le D,\ \height(q)\le B \right\}.
\]
We call these sets \emph{target-height slices}: they consist of the points with height at most \(B\) subject to the indicated field or degree restriction. Their finiteness follows from Northcott's theorem \cite[Theorem~1.6.8]{BombieriGubler2006}, applied to the coordinates. A height bound without a field or degree bound is insufficient; the examples in \cref{sec:jump} have bounded target height and arbitrarily large coordinate-field degree. We use \(g \in G\) and \(t \in T\) when their guarding roles matter, and \(p,q\) for general geometric endpoints.

An \emph{obstruction to visibility} will mean a stated geometric or arithmetic condition preventing sight between distinct points at the allowed budgets. The qualifier \emph{degree} or \emph{height} names the source of an arithmetic restriction. The existence of a contained regular arc and coverage by a given number of guards remain separate questions.

\begin{remark}\label{rem:coefficient-field}
    The carrier field is \(\Q\), independently of the target field. If carriers are instead defined over \(K\), every pair of \(K\)-points has a \(K\)-line; degree-one visibility at unbounded coefficient height is then segment visibility. We return to this change in \cref{prob:field-kernel}.

    For a primitive integer coefficient vector, the maximum of its \(\ell\)-adic absolute values is one at every rational prime \(\ell\). Consequently its global projective height equals its largest absolute coefficient. Replacing the present coefficient bound by that global height does not change the rational-carrier model.
\end{remark}

\begin{figure}[htbp]
    \centering
    \begin{tikzpicture}[
        x=1.5cm,y=1.1cm,
        axis/.style={draw=black!25,line width=.5pt,->},
        curve/.style={draw=blue!70!black,line width=1.3pt},
        conjugate/.style={draw=orange!85!black,line width=1.3pt},
        dot/.style={circle,fill=black,inner sep=1.6pt},
        every node/.style={font=\small}
        ]
        \foreach \offset in {0,4.4}{
        \begin{scope}[shift={(\offset,0)}]
            \draw[axis] (-1.5,0)--(1.55,0) node[right] {\(X\)};
            \draw[axis] (0,-1.95)--(0,2) node[above] {\(Y\)};
        \end{scope}
        }

        \draw[curve] (-1.3,{-1.3*sqrt(2)})--(1.3,{1.3*sqrt(2)});
        \draw[conjugate] (-1.3,{1.3*sqrt(2)})--(1.3,{-1.3*sqrt(2)});
        \node[dot] at (0,0) {};
        \node at (0,-2.45) {(a) \(Z(Y^2-2X^2)\)};

        \begin{scope}[shift={(4.4,0)}]
            \draw[curve,domain=-1.49:1.49,samples=250,variable=\t]
            plot ({\t*\t-1},{\t*(\t*\t-1)});
            \node[dot,label=right:{\(p\)}] at (.75,{.75*sqrt(1.75)}) {};
            \node[dot,label=right:{\(q\)}] at (.75,{-.75*sqrt(1.75)}) {};
            \draw[fill=white,draw=red!75!black,line width=1pt]
            (0,0) circle[radius=3.5pt];
            \node at (0,-2.45) {(b) \(Z(Y^2-X^2(X+1))\)};
        \end{scope}
    \end{tikzpicture}

    \caption{(a) Two conjugate lines. (b) An absolutely irreducible cubic; removing its node separates \(p\) and \(q\).}
    \label{fig:irreducibility-singularities}
\end{figure}
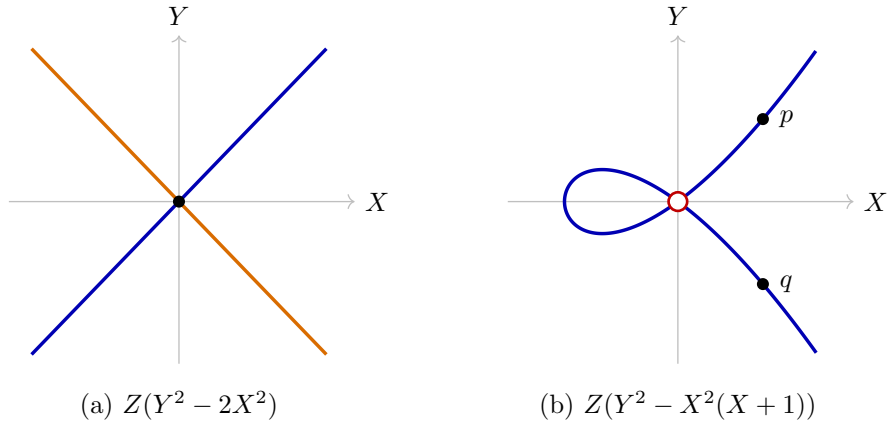

\begin{remark}\label{rem:conditions-example}
    The polynomial \(Y^2-2X^2\) is irreducible over \(\Q\) but its zero locus splits into two conjugate lines over \(\Qbar\). Absolute irreducibility rules out such unions disguised by rational coefficients. The cubic \(Y^2-X^2(X+1)\), by contrast, is absolutely irreducible and has a node at the origin. Its full real locus permits branch switching there, but its regular locus does not. \Cref{ex:nodal-incidence} will show that this distinction changes a guard number at a fixed budget, even after all alternative admissible carriers are considered.
\end{remark}

Degrees and coefficient heights have so far been measured in fixed affine coordinates. How much do these parameters depend on our choice of coordinates? We conclude by showing that degree does not depend on it at all, while coefficient height depends on it only up to a constant.

\begin{proposition}\label{prop:affine}
    Let \(A(z)=Uz+b\), where \(U\in\mathrm{GL}_2(\Q)\) and \(b\in\Q^2\). Then
    \[
        \gamma_{d,\infty}(A(P);A(T),A(G))=\gamma_{d,\infty}(P;T,G).
    \]

    For each \(d\), a constant \(C_A(d)\ge1\) satisfies
    \[
        \gamma_{d,C_A(d)H}(A(P);A(T),A(G))
        \le \gamma_{d,H}(P;T,G).
    \]
    Applying this assertion to \(A^{-1}\) gives the reverse comparison.
\end{proposition}

\begin{proof}
    Substitute \(A^{-1}\) into an admissible polynomial, clear denominators, and take the primitive part. Degree and absolute irreducibility are preserved. The gradient changes by an invertible matrix, and \(A\) maps regular components homeomorphically onto regular components. Expanding the finitely many monomials of degree at most \(d\) bounds the new coefficients by \(C_A(d)\norm{F}_\infty\). Now transport each sight arc with its defining polynomial and then apply the same argument to \(A^{-1}\).
\end{proof}

Thus degree is invariant under rational affine changes, whereas exact coefficient-height levels generally are not. A projective transformation moving the line at infinity into the gallery raises an additional question about which real arcs remain in the allowed region, but that is not covered by this affine comparison.

\section{Computability and finite target sets}\label{sec:finite}

Before turning to particular galleries we observe some basic combinatorial structure underlying finite target sets in galleries. When guards are chosen from the finite target set itself, guarding becomes a domination problem on its visibility graph. At finite budgets, the carrier list is also finite, which will let us compute visibility exactly. To that end, for a finite target set \(T\) define the \emph{visibility graph} \(\operatorname{Vis}_{d,H}(P,T)\) as the graph on \(T\) in which distinct targets are adjacent exactly when they see each other. We recall that a \emph{dominating set} of a graph is a set of vertices such that every vertex lies in it or is adjacent to one of its elements, and the \emph{domination number} is the least size of a dominating set.

\begin{remark}
    The connected components of the regular traces of all admissible carriers define a hypergraph on \(T\): its hyperedges are their nonempty intersections with \(T\). This records which targets lie together on one regular component of one carrier. The visibility graph is the \(2\)-section of this hypergraph, although a clique need not lie on one regular component of a single carrier.
\end{remark}

\begin{proposition}
    \hfill
    \begin{enumerate}
        \item Visibility is reflexive and symmetric, but need not be transitive.
        \item For a finite target set \(T\), the domination number of \(\operatorname{Vis}_{d,H}(P,T)\) is \(\gamma_{d,H}(P;T,T)\).
        \item If \(d\le d'\), \(H\le H'\), \(T'\subseteq T\), and \(G\subseteq G'\), then
        \[
            \gamma_{d',H'}(P;T',G')\le\gamma_{d,H}(P;T,G).
        \]
    \end{enumerate}
\end{proposition}

\begin{figure}[htbp]
    \centering
    \begin{tikzpicture}[x=1cm,y=1cm,font=\small]
        \foreach \shift/\h/\gam in {0/1/2,3.9/2/1,7.8/3/1}{
        \begin{scope}[xshift=\shift cm]
            \coordinate (a) at (0,0); \coordinate (b) at (0,1.7);
            \coordinate (c) at (2.3,1.7); \coordinate (d) at (2.3,0);
            \draw[blue!65!black,very thick] (b)--(a)--(d)--(c);
            \ifnum\h>1 \draw[orange!85!black,very thick] (a)--(c); \fi
            \ifnum\h>2 \draw[purple!80!black,very thick] (b)--(c) (b)--(d); \fi
            \foreach \pt/\lbl/\pos in {a/0/below,b/1/above,c/2/above,d/3/below}
            {\fill (\pt) circle(2pt); \node[\pos=4pt] at (\pt) {\(q_{\lbl}\)};}
            \draw[black,line width=1pt] (a) circle(6pt);
            \ifnum\h=1 \draw[black,line width=1pt] (d) circle(6pt); \fi
            \node[anchor=north,align=center] at (1.15,-.65)
            {\(H=\h\)\\\(\gamma_{1,\h}(\mathcal{S};T,T)=\gam\)};
        \end{scope}
        }
    \end{tikzpicture}

    \caption{The degree-one visibility graphs of \cref{ex:path-four} at carrier heights one, two, and three. Rings indicate minimum guard sets chosen from the targets.}
    \label{fig:path-four}
\end{figure}
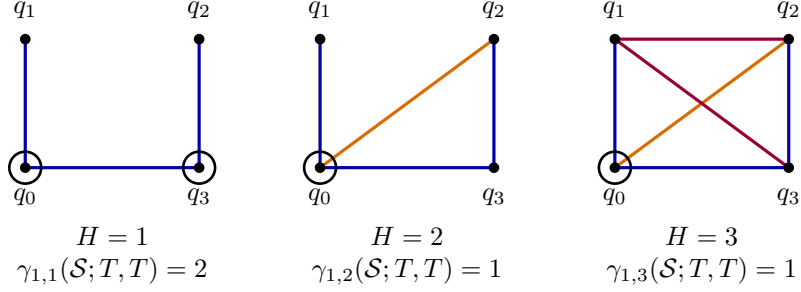

\begin{proof}
    Reflexivity and symmetry are immediate from \cref{def:visibility}, and the \(H=1\) graph in \cref{fig:path-four} shows that transitivity can fail. A target sees exactly itself and its neighbors in the visibility graph, so the covers of \(T\) by guards from \(T\) are exactly the dominating sets. This gives (ii). For (iii), a cover remains a cover when the budget or the guard set grows, and a cover of \(T\) is a cover of every subset of \(T\).
\end{proof}

Arithmetic sight can be restrictive, even at degree one, as the next example demonstrates.

\begin{example}\label{ex:path-four}
    In the square \(\mathcal{S}\), let
    \[
        q_0=(0,0),\quad q_1=(0,\tfrac13),\quad
        q_2=(\tfrac13,\tfrac23),\quad q_3=(\tfrac12,\tfrac12),
        \qquad T=\{q_0,q_1,q_2,q_3\}.
    \]
    Their joining lines have coefficient heights one for \(q_0q_1,q_0q_3,q_2q_3\), two for \(q_0q_2\), and three for \(q_1q_2,q_1q_3\). Thus the degree-one graph is the path \(q_1-q_0-q_3-q_2\) at \(H=1\), and becomes \(K_4\) at \(H=3\).
\end{example}

The following example shows that the regularity requirement in \cref{def:visibility} is consequential.

\begin{figure}[htbp]
    \centering
    \begin{tikzpicture}[x=1.92cm,y=1.33cm,font=\footnotesize]
        \foreach \shift/\panel in {0/a,4/b}{
        \begin{scope}[xshift=\shift cm]
            \draw[gray!55] (-1,-1) rectangle (1,1);
            \begin{scope}
                \clip (-1,-1) rectangle (1,1);
                \draw[gray!65,thin] plot[domain=-1:1,samples=70]
                ({\x*\x-1},{\x*(\x*\x-1)});
                \draw[blue!70!black,thick] plot[domain=1:1.325,samples=45]
                ({\x*\x-1},{\x*(\x*\x-1)});
                \draw[orange!85!black,thick] plot[domain=-1.325:-1,samples=45]
                ({\x*\x-1},{\x*(\x*\x-1)});
            \end{scope}
            \fill (.5625,.703125) circle(1.4pt);
            \fill (.5625,-.703125) circle(1.4pt);
            \node[right=2pt] at (.5625,.703125) {\(a_+\)};
            \node[right=2pt] at (.5625,-.703125) {\(a_-\)};
            \draw[black] (.5625,.703125) circle(3.4pt);
            \ifnum\shift=0
            \draw[fill=white] (0,0) circle(2pt);
            \draw[black] (.5625,-.703125) circle(3.4pt);
            \else
            \fill (0,0) circle(2pt);
            \fi
            \node at (0,-1.45) {(\panel)};
        \end{scope}
        }
    \end{tikzpicture}

    \caption{(a) Two target guards are needed when the node is removed.
    (b) Allowing passage through it permits one guard. Rings indicate guards chosen from \(T\).}
    \label{fig:nodal-practical}
\end{figure}
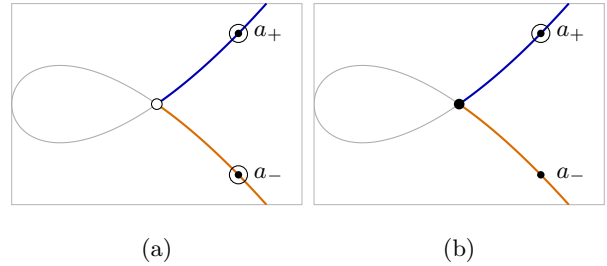

\begin{example}\label{ex:nodal-incidence}
    In \(P=[-1,1]^2\), take \(d=3,H=1\) and
    \[
        F=Y^2-X^2(X+1),\qquad
        T=\{a_+,a_-\},\qquad a_\pm=(\tfrac9{16},\pm\tfrac{45}{64}).
    \]
    The cubic is absolutely irreducible and has only one singular point, the origin. The targets lie on different connected components of the regular trace. Moreover, \(\mathrm{C}_F\) is the only admissible carrier through either target at this budget (we prove this in \cref{sec:nodal-exclusions}). Hence \(\gamma_{3,1}(P;T,T)=2\). Allowing sight to pass through the node would drop the guard number to one at the same budget (\cref{fig:nodal-practical}).
\end{example}

\begin{corollary}\label{cor:forced}
    At a given budget \(d,H\), targets on no admissible regular trace must guard themselves. If such a target does not belong to \(G\), then \(T\) has no cover, and \(\gamma_{d,H}(P;T,G)=\infty\).
\end{corollary}

For finite sets of explicitly represented algebraic guards and targets, finiteness of the carrier list also lets us compute the guard number exactly. For a finite budget \(d,H\), there are \(M_d=\binom{d+2}{2}\) possible monomials, so a bound on the number of admissible carriers is:
\begin{equation}\label{eq:curve-count}
    \card{\mathcal F_{d,H}} \le \frac{(2\lfloor H \rfloor+1)^{M_d}-1}{2}.
\end{equation}

This counts coefficient choices modulo sign, including inadmissible ones, so it is a coarse bound.

\begin{proposition}\label{thm:computable}
    For rational \(P\), finite sets \(T,G\) of explicitly represented real algebraic points, and positive integers \(d,H\), the guard number \(\gamma_{d,H}(P;T,G)\) is exactly computable. In particular, this holds for \(T = G = T_{K,B}(P)\), with \(K\) and its real embedding specified and \(B \ge 1\) explicitly real algebraic.
\end{proposition}

\begin{proof}
There are four finite steps.

\begin{enumerate}[label=\arabic*.,leftmargin=*,nosep]
    \item Enumerate bounded primitive integer polynomials of degrees one through \(d\), and test absolute irreducibility by exact algebraic factorization.
    \item Compute each regular trace's connected components, including polygon boundary points but excluding singular points, by the semialgebraic procedure of \cite[Chapter~16]{BasuPollackRoy2006}.
    \item For each guard, list the targets sharing one such component, and add the guard itself if it is a target.
    \item Test subsets of \(G\) in increasing size until their visible sets cover \(T\); return \(\infty\) if none does.
\end{enumerate}

For the specialization, enumerate the bounded-height target list first, as described in \cref{sec:finite-implementation}.
\end{proof}

\section{Straight-line visibility and degree-constrained guarding}\label{sec:linear}

Our first task is to check that the model does not lose the classical theory. At degree one and sufficiently large height on a rational polygon, it should reduce to rational straight-line guarding. Once this is settled, the interesting questions are what curved sight adds, and what the arithmetic of the targets takes away. We therefore compare straight-line guarding with curved sight in a one-sided comb, and then with guarding all points over a quadratic field, first by point guards and then by edge guards. Finally, alternating notches will show why degree still limits both point and edge guards.

\subsection{Straight-line sight and a one-sided comb}\label{sec:comb}
\begin{proposition}\label{prop:linear}
    Distinct points \(g,t \in P(\Q)\) of projective heights at most \(A,B\) have a primitive joining-line equation of height at most \(2AB\). For \(H \in [2AB, \infty]\), they see each other at degree one exactly when \([g,t] \subseteq P\).
\end{proposition}

\begin{proof}
    Take primitive homogeneous integer coordinates for the points. The coefficients of their line are the primitive part of their cross product, each coordinate bounded by \(2AB\). A connected line trace contains the segment between any two of its points, proving the equivalence.
\end{proof}

This bound makes Fisk's vertex-guard construction available at a common coefficient height. It also relates rational target testing to a real gallery guarded from rational positions.

\begin{corollary}\label{cor:fisk}
    If \(P\) has \(n\) vertices of height at most \(V\), then, for \(B \ge V\),
    \[
        \gamma_{1, 2VB}(P;T_{\Q,B}(P), T_{\Q,B}(P)) \le \lfloor n/3 \rfloor.
    \]
\end{corollary}

\begin{proof}
    Triangulate, three-color, and choose a smallest color class \cite{Fisk1978}. These marked vertices cover the whole real polygon by segments and, as \(P\) has rational vertices with heights \(V \le B\), they belong to the target-height slice. By \cref{prop:linear}, these vertices cover the target-height slice using degree-one carriers of coefficient height at most \(2VB\). There are at most \(\lfloor n/3\rfloor\) such guards.
\end{proof}

Let \(\gamma_{\Q}(P)\) denote the minimum number of rational guards seeing every real point by segments. Then rational targets suffice to determine it.

\begin{proposition}\label{prop:rational-specialization}
    \[
        \gamma_{1,\infty}(P;P(\Q),P(\Q)) = \gamma_{\Q}(P).
    \]
\end{proposition}

\begin{proof}
    The rational vertices provide a finite segment cover of \(P\), and their joining lines to rational targets are rational, so both guard numbers are finite. \emph{A fortiori} a minimum set of rational guards seeing every real point of \(P\) also sees every rational point. Conversely, for a minimum rational guard set seeing all rational targets, approximate any real target by a sequence of rational ones. One guard sees a subsequence and closedness keeps its limiting segment in \(P\). The same guard set therefore covers all real targets, whence the stated guard numbers are equal.
\end{proof}

This need not equal the unrestricted real optimum: an integer-vertex polygon can need four rational guards but only three real guards; see \cite[Theorem~1]{AbrahamsenAdamaszekMiltzow2017}.

The textbook comb \cite[Section~1.4, Figure~1.17]{DevadossORourke2025} suggests that adding teeth requires guards. We use rectangular teeth to check equations directly. The targets below already require the full classical count, yet one central guard sees them along parabolas.

\begin{theorem}\label{thm:central-comb}
    For \(m\ge1\), put \(L=2m+1\), \(a_i=(2i-m+1)L\), and define
    \[
         \mathcal P_m=([-mL,mL]\times[-2m,2m])\cup
         \bigcup_{i=0}^{m-1}[a_i-1/2,a_i+1/2]\times[2m,2m+3/2].
    \]
    Let \(T_m=\{(a_i,L):0\le i<m\}\). The polygon has \(4m+4\) vertices and classical guard number \(m\), already required for \(T_m\). For \(G=\mathcal P_m\) or \(\mathcal P_m(\Q)\),
    \[
         \gamma_{1,H}(\mathcal P_m;T_m,G)=m,\qquad
         \gamma_{2,H}(\mathcal P_m;T_m,G)=1,
         \qquad H\in[4m^2,\infty].
    \]
    The second cover uses \(g=(0,0)\), distinct from all targets.
\end{theorem}

The sight cones below the tooth openings are disjoint. For the curved cover, the arcs \((a_i(2u-u^2),Lu)\), \(0\le u\le1\), bend through the base and enter their teeth vertically. Their defining equations are the primitive parts of \((2i-m+1)Y^2-2(2i-m+1)LY+LX\). The containment and height checks are in \cref{sec:comb-deferred}.

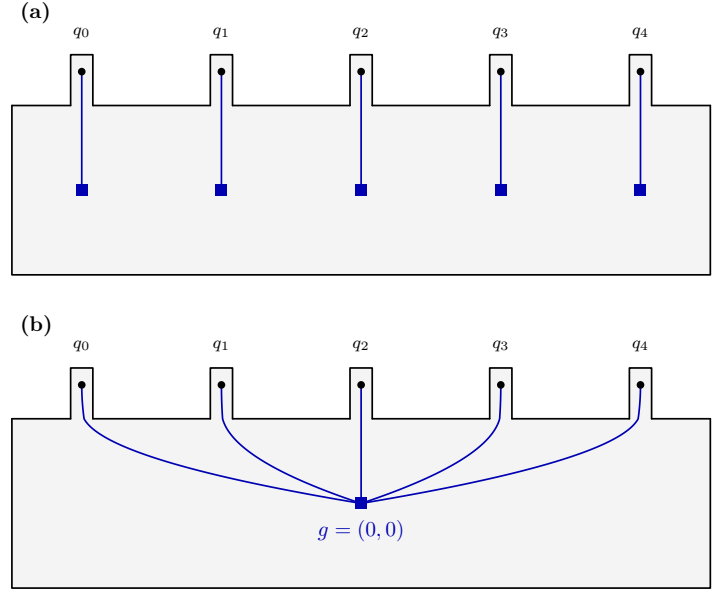
\begin{figure}[t]
    \centering
    \scalebox{0.8}{%
    \begin{tikzpicture}[x=.105cm,y=.14cm,line cap=round,line join=round]
        \foreach \shift/\mode/\heading in
        {0/0/{(a)},-37/1/{(b)}}{
        \begin{scope}[yshift=\shift*.14cm]
            \fill[gray!9] (-55,-10) rectangle (55,10);
            \draw[thick] (-55,10)--(-55,-10)--(55,-10)--(55,10);
            \foreach \a in {-44,-22,0,22,44}{
            \fill[gray!9] (\a-1.75,10) rectangle (\a+1.75,16);
            \draw[thick]
            (\a-1.75,10)--(\a-1.75,16)
            --(\a+1.75,16)--(\a+1.75,10);
            }
            \draw[thick] (-55,10)--(-45.75,10) (-42.25,10)--(-23.75,10)
            (-20.25,10)--(-1.75,10) (1.75,10)--(20.25,10)
            (23.75,10)--(42.25,10) (45.75,10)--(55,10);
            \node[anchor=west,font=\small\bfseries] at (-55,21) {\heading};
            \foreach \a/\i in {-44/0,-22/1,0/2,22/3,44/4}{
            \ifnum\mode=0
            \draw[blue!70!black,thick] (\a,0)--(\a,14);
            \fill[blue!70!black] (\a-.9,-.65) rectangle (\a+.9,.65);
            \else
            \draw[blue!70!black,thick,domain=0:1,samples=100,variable=\u]
            plot ({\a*(2*\u-\u*\u)},{11*\u+3*max(0,11*\u-10)});
            \fi
            \fill[black] (\a,14) circle[radius=1.8pt];
            \node[above,font=\scriptsize] at (\a,17) {\(q_{\i}\)};
            }
            \ifnum\mode=1
            \fill[blue!70!black] (-.9,-.65) rectangle (.9,.65);
            \node[below,font=\small,blue!70!black] at (0,-1) {\(g=(0,0)\)};
            \fi
        \end{scope}
        }
    \end{tikzpicture}
    }

    \caption{One-sided comb with \(m=5\) and \(H=100\). The same polygon and targets are shown in both panels; tooth widths and heights are exaggerated for clarity. Squares are guards; dots are targets. In (a), five straight-line guards cover the whole polygon, and the targets already require five guards even with unrestricted guard locations. In (b), the four parabolic sight arcs and middle vertical segment connect the single central guard to all five targets at degree two. The exact construction is given in \cref{thm:central-comb}.}
    \label{fig:central-comb}
\end{figure}

\subsection{From finite targets to quadratic galleries}
We now replace a finite target set by all \(K\)-points, keeping rational carrier coefficients. At degree one, this arithmetic change can outweigh the polygon's geometry, even in a convex gallery.

\begin{proposition}\label{prop:point-field-linear}
    For a real number field \(K \ne \Q\),
    \[
     \gamma_{1,\infty}(P;P(K),P) = \infty.
    \]
\end{proposition}

\begin{proof}
    Consider a finite guard set and any interior rectangle within the gallery. Choose any point within the rectangle with coordinates \((r,\beta)\), where \(r \in \Q\) is different from every guard's abscissa and \(\beta \in K\) is irrational. Then the only rational line through \((r,\beta)\) is \(X = r\), which contains no guard. So no finite cover exists, whence the guard number is infinite.
\end{proof}

For quadratic fields a small construction reverses this failure. Its geometric hypothesis will also supply our edge covers.

\begin{lemma}\label{lem:convex-connector}
    Suppose \(P\) is convex, \(g \in P(\Q)\), and a real algebraic target \(t \in P\) of coordinate degree \(e\) lies on a rational line. Then a rational absolutely irreducible curve of degree at most \(e\) contains a sight arc between \(g\) and \(t\). If \(t\) is rational, or \(g\) lies on the rational line containing \(t\), their joining segment suffices. Coincident points use reflexivity.
\end{lemma}

\begin{figure}[htbp]
    \centering
    \begin{tikzpicture}[
        x=1cm,y=1cm,
        font=\footnotesize,
        line cap=round,line join=round,
        gallery/.style={draw=gray!70,fill=gray!8,thick},
        reference/.style={draw=black,densely dashed},
        targetline/.style={draw=gray!65,dashed},
        sight/.style={draw=cyan!65!blue,line width=1.2pt}
    ]
        \pgfmathsetmacro{\connectoralpha}{sqrt(2)/2}
        \pgfmathsetmacro{\connectordelta}{\connectoralpha-.25}

        \foreach \offset in {0,10}{
            \begin{scope}[xshift=\offset cm,x=3cm,y=3cm]
                \path[gallery] (0,0) rectangle (1,1);

                \draw[targetline] (0,.5)--(1,.5);
                \node[above left,text=gray!65] at (.98,.5) {\(\ell\)};
                \draw[reference] (0,0)--(\connectoralpha,.5);

                \ifnum\offset=10
                    \draw[sight,domain=0:1,samples=100,variable=\s]
                        plot (
                            {.25+\connectordelta*\s
                             -.25*(8-16*(.25+\connectordelta*\s)^2)/7},
                            {.5*(1-(8-16*(.25+\connectordelta*\s)^2)/7)}
                        );
                \fi

                \fill (0,0) circle[radius=1.8pt];
                \node[below left] at (0,0) {\(g\)};

                \fill (\connectoralpha,.5) circle[radius=1.8pt];
                \node[above right] at (\connectoralpha,.5) {\(t\)};

                \draw[gray!75,fill=white]
                    (.25,.5) circle[radius=1.8pt];
                \node[below] at (.25,.5) {\(b_a\)};
            \end{scope}
        }

        \begin{scope}[xshift=6cm,yshift=1.5cm,x=2cm,y=1.5cm]
            \path[gallery]
                (0,1)--(1,1)--(.5,-1)--(-.5,-1)--cycle;

            \draw[targetline,->] (-.33,0)--(.9,0)
                node[right] {\(u\)};
            \draw[->,gray!65] (0,-1.08)--(0,1.14)
                node[above] {\(w\)};

            \draw[reference] (0,1)--(\connectordelta,0);
            \draw[sight,domain=0:1,samples=100,variable=\s]
                plot (
                    {\connectordelta*\s},
                    {(8-16*(.25+\connectordelta*\s)^2)/7}
                );

            \fill (0,1) circle[radius=1.8pt];
            \node[right=-15pt] at (0,0.9) {\(g'\)};

            \fill (\connectordelta,0) circle[radius=1.8pt];
            \node[below right] at (\connectordelta,0) {\(t'\)};

            \draw[gray!75,fill=white]
                (0,0) circle[radius=1.8pt];
            \node[below left] at (0,0) {\(0\)};
        \end{scope}

        \draw[->,thick] (3.3,1.5)--(4.7,1.5)
        node[midway,above=5pt] {\(\Phi_a^{-1}\)};
        \draw[->,thick] (8.35,1.5)--(9.7,1.5)
        node[midway,above=5pt] {\(\Phi_a\)};

        \node at (1.5,-.5) {(a)};
        \node at (6.5,-.5) {(b)};
        \node at (11.5,-.5) {(c)};
    \end{tikzpicture}

    \caption{The convex connector construction: (a) the gallery setup; (b) its inverse image, with the constructed graph; (c) the graph's image joining the guard and target. Black dashed segments provide the reference for the estimate. Here \(g=(0,0)\), \(t=(\sqrt{2}/2,1/2)\), \(b=(0,1/2)\), \(v=(1,0)\), \(a=1/4\), and \(M(X)=X^2-1/2\).}
    \label{fig:convex-connector}
\end{figure}
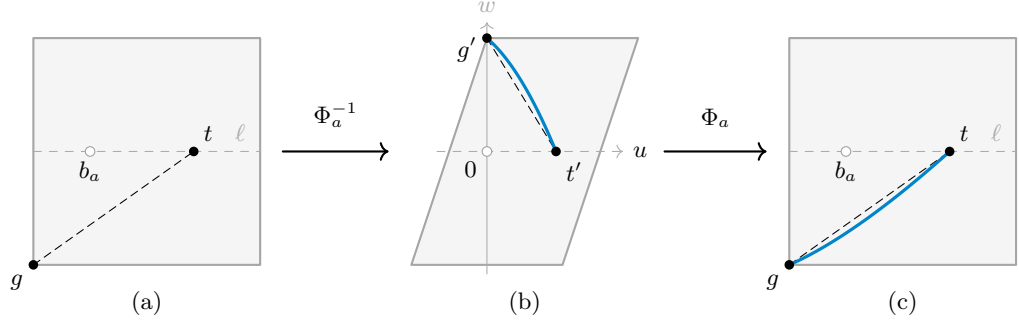

Suppose \(t\) is nonrational and \(g\) lies off its rational line. Write \(t=b+\alpha v\), with \(b,v \in \Q^2\) and \(v \ne 0\), and let \(M\) be the minimal polynomial of \(\alpha\). For rational \(a\) close to \(\alpha\), put \(b_a=b+av\). The affine parameterization
\[
    \Phi_a(u,w)=b_a+uv+w(g-b_a)
\]
sends \((0,1)\) to \(g\) and \((\alpha-a,0)\) to \(t\). To construct a rational polynomial graph through these normalized points, use \(w=M(a+u)/M(a)\): the shifted minimal polynomial vanishes at \(\alpha-a\), and division by \(M(a)\) makes its value at zero equal to one. For rational \(a\) sufficiently close to \(\alpha\), applying \(\Phi_a\) gives a regular connecting arc inside \(P\). The proof, including the containment estimate, is in \cref{sec:edge-deferred}; \cref{fig:convex-connector} illustrates the construction.

\begin{corollary}\label{cor:point-quadratic}
    For a real quadratic field \(K\), every rational point in a convex \(P\) sees all of \(P(K)\) at degree two and unrestricted height. Thus
    \[
        \gamma_{1,\infty}(P;P(K),P(\Q)) = \infty,
        \qquad \gamma_{2,\infty}(P;P(K),P(\Q)) = 1.
    \]
    For any rational simple \(n\)-gon,
    \[
        \gamma_{2,\infty}(P;P(K),P(\Q)) \le \lfloor n/3 \rfloor.
    \]
\end{corollary}

\begin{proof}
    Set \(K=\Q(\sqrt u)\), with \(u>1\) a squarefree integer. Then a \(K\)-point can be written as \(a+\sqrt u\,b\), with \(a,b \in \Q^2\), identifying it on the rational line through \(a\) in direction \(b\). If the polygon is convex, then \cref{prop:point-field-linear} verifies the degree-one guard number and the previous lemma verifies the degree-two guard number.

    Otherwise, triangulate the rational polygon \(P\) using its vertices, three-color the triangulation, and select a smallest color class. In each triangle, apply the lemma with its selected vertex as the guard. This verifies that \(\gamma_{2,\infty}(P;P(K),P(\Q)) \le \lfloor n/3 \rfloor\).
\end{proof}

The proof of \cref{cor:point-quadratic} demonstrates that we can apply \cref{lem:convex-connector} within a rational triangle contained in \(P\), without requiring the whole polygon to be convex. For an algebraic target on a side of that triangle, we can also make the lemma's construction effective. For a side with endpoints \(b,c\), write the target as \(t=b+\alpha(c-b)\). Bounds on the distances between distinct roots of \(\alpha\)'s minimal polynomial control the rational approximation and give an explicit coefficient-height bound.

The classical \emph{kernel} consists of all points that see every point of \(P\) by segments. So a rational interior guard in the kernel spans a triangle contained in \(P\) with every edge of \(P\). The following corollary therefore applies the effective construction to every algebraic boundary target from such a guard.

\begin{corollary}\label{cor:star-boundary-height}
    Let \(P\) be a rational polygon and let \(g \in \operatorname{Ker}(P) \cap (\operatorname{int}P)(\Q)\). Let \(t=b+ \alpha(c-b)\) lie in an edge with rational endpoints \(b,c\), where \(0 < \alpha < 1\) is nonrational and has minimal polynomial of degree \(e\). If \(g\) is off the edge line, there is an admissible carrier of degree \(e\) joining \(g\) to \(t\) inside \(P\), with the explicit coefficient-height bound in \cref{eq:star-height-bound}. This degree is minimal. If \(g\) lies on the edge line, that line suffices at height at most \(2\height(b)\height(c)\). Rational boundary targets use segments of height at most \(2\height(g)\height(t)\).
\end{corollary}

Here the polygon contributes a contained triangle; the target's algebraic degree bounds the degree of the constructed carrier. The proof in \cref{sec:edge-deferred} gives the height bound directly from the minimal polynomial and the triangle's rational coordinates. General polygons and targets will require the arc-lifting construction of \cref{sec:lifting}.

\subsection{Closed edges and the chords of a polygon}\label{sec:edge-guards}
The textbook's Unsolved Problem~5 asks how many edges must be selected to guard a polygon \cite[Section~1.4]{DevadossORourke2025}. 

\begin{definition}[Edge guard visibility]
    Consider a rational polygon \(P\), a target set \(T \subseteq P\), and a budget \(d,H\). A closed edge \(E\) sees \(t \in T\) if there is a \(g \in E\) such that \(g \sim_{d,H} t\). 

    Define \(\eta_{d,H}(P;T)\), the \emph{edge guard number}, as the minimum number of closed edges needed to cover \(T\), with value \(\infty\) if no edge collection forms a cover. Edges are maximal, and their observing points may be real. Let \(\eta_{\R}(P)\) be the classical closed-edge count, i.e., the minimum number of closed edges needed to see every point of \(P\) by segments.
\end{definition}

\begin{proposition}\label{prop:edge-rational}
    If the vertex heights are at most \(V\), then
    \[
        \eta_{1,\infty}(P;P(\Q)) = \eta_{\R}(P).
    \]
    For some \(B_0(P)\) and every \(B \ge B_0(P)\),
    \[
        \eta_{1,2VB}(P;T_{\Q,B}(P)) =
        \eta_{1,\infty}(P;T_{\Q,B}(P)) = \eta_{\R}(P).
    \]
    For every \(B \ge 1\), an edge sees a rational target of height at most \(B\) by a segment exactly when it sees it at degree one and coefficient height \(2VB\).
\end{proposition}

An edge seeing a rational target can use a line through a polygon vertex, which supplies the height bound. Each inadequate edge collection misses a \emph{rational} target, by the closedness of its visibility region. Choosing one missed rational target for each edge collection smaller than a classical minimum cover gives a finite target set that already requires the same number of edges. The proof, including the coefficient-height bound, is in \cref{sec:edge-deferred}.

We can use algebraic boundary targets as a first point of contrast between the theories. A convex gallery has a one-edge classical cover. But a degree-\(e\) point in an edge cannot lie on a rational curve of smaller degree unless that curve contains the entire edge line. Irreducibility then leaves only the line as a carrier. A selected closed edge can see these targets on itself and its two neighbors, using their shared vertices. At degree \(e\), \cref{lem:convex-connector} joins a rational point on any selected edge to every target in a finite boundary slice, at a sufficiently large common coefficient-height budget.

\begin{theorem}\label{thm:edge-boundary}
    Let \(P\) be a convex rational \(n\)-gon and let \(K \subset \R\) have degree \(e \ge 2\). Let \(h_P\) be the largest coefficient height of the primitive integer equations of its edge lines. For some \(B_0(P,K)\), the target set \(T=T_{K,B}(P)\cap\partial P\) satisfies
    \[
     \eta_{d,H}(P;T) = \left \lceil \frac{n}{3} \right \rceil,
     \qquad B \ge B_0(P,K),\quad 1 \le d<e,\quad H \in [h_P,\infty].
    \]
    For every finite \(B \ge 1\), there is a finite \(H_*(P,K,B) \ge h_P\) such that any edge covers \(T\) at degree \(e\) and every \(H \in [H_*,\infty]\). Thus \(\eta_{e,H}(P;T)=1\) when \(T\) is nonempty. The same statements, with possibly different thresholds, hold for the height slice of all real algebraic boundary points whose coordinate-field degree is at most \(e\).
\end{theorem}

The one-in-three count here comes from boundary adjacency, not from the triangulation in Fisk's proof. Height \(h_P\) makes the edge lines available; at degree \(e\), a finite maximum over the convex connections of \cref{lem:convex-connector} makes every edge a possible guard. The proof is in \cref{sec:edge-deferred}.

For interior quadratic targets, the arrangement of edges matters beyond boundary adjacency. Each nonrational quadratic target lies on a unique rational line: the line through it and its conjugate. At degree one, an observing point must lie on this line, and the segment joining it to the target must stay inside \(P\). Which edges can see the target therefore depends on where this segment can reach the boundary. We describe these possibilities using the polygon's chords and obtain an exact edge guard count.

A \emph{proper interior chord} joins non-vertex points on distinct edges and has its open segment in \(\operatorname{int}(P)\). Let \(b(P)\) be the minimum number of polygon edges that can be selected so that every proper interior chord has at least one endpoint on a selected edge. To express this as a graph problem, form the \emph{chord graph} of \(P\): its vertices are the polygon's edges, and two vertices are adjacent when a proper interior chord joins the corresponding edges. A \emph{vertex cover} is a set of graph vertices containing at least one endpoint of every graph edge. Thus selecting polygon edges that meet every proper interior chord is equivalent to choosing a vertex cover of the chord graph, and \(b(P)\) is its minimum size.

\begin{proposition}\label{prop:edge-chords}
    Let \(K\) be a real quadratic field and \(P\) have \(n\) vertices. For some \(B_0(P,K)\),
    \[
        \eta_{1,H}(P;T_{K,B}(P)) = b(P),
        \qquad B \ge B_0(P,K),\quad H \in [2B^2,\infty].
    \]
    Also \(\eta_{1,\infty}(P;P(K)) = b(P) \le n-1\). Any \(n-1\) edges cover \(T_{K,B}(P)\) at height \(2B^2\) for every \(B \ge 1\). The statements also hold for height slices of \emph{all} real points of coordinate degree at most two.
\end{proposition}

For the lower bound, suitable quadratic targets force the selected edges to form a vertex cover of the chord graph. Conversely, such a selection reaches every target along a contained segment of a rational line, including cases where the segment ends at a polygon vertex. The bound \(2\height(t)^2\) on a quadratic target's line height supplies the coefficient budget. The proof is in \cref{sec:edge-deferred}.

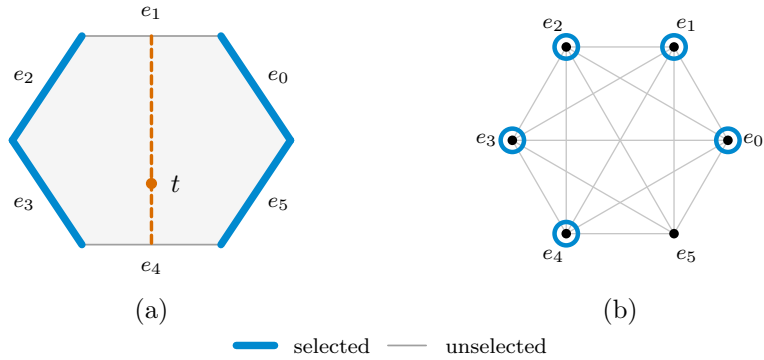
\begin{figure}[htbp]
    \centering
    \begin{tikzpicture}[x=.92cm,y=.92cm,font=\small,line cap=round,line join=round]
        \begin{scope}
            \coordinate (v0) at (2,0); \coordinate (v1) at (1,1.5);
            \coordinate (v2) at (-1,1.5); \coordinate (v3) at (-2,0);
            \coordinate (v4) at (-1,-1.5); \coordinate (v5) at (1,-1.5);
            \fill[gray!8] (v0)--(v1)--(v2)--(v3)--(v4)--(v5)--cycle;
            \draw[gray!75,line width=.65pt] (v0)--(v1)--(v2)--(v3)--(v4)--(v5)--cycle;
            \foreach \i in {0,...,5}{
            \pgfmathtruncatemacro{\j}{mod(\i+1,6)}
            \coordinate (e\i) at ($(v\i)!.5!(v\j)$);
            \node[font=\scriptsize] at ($(0,0)!1.22!(e\i)$) {\(e_{\i}\)};
            }
            \foreach \i in {0,2,3,5}{
            \pgfmathtruncatemacro{\j}{mod(\i+1,6)}
            \draw[cyan!70!blue,line width=2.8pt] (v\i)--(v\j);
            }
            \draw[orange!85!black,densely dashed,line width=1.3pt] (e1)--(e4);
            \fill[orange!85!black] ($(e1)!{sqrt(2)/2}!(e4)$) circle[radius=2.1pt];
            \node[anchor=west] at (.13,-.62) {\(t\)};
            \node[anchor=north] at (0,-2.12) {(a)};
        \end{scope}
        \begin{scope}[xshift=6.2cm]
            \foreach \i in {0,...,5} \coordinate (c\i) at ({60*\i}:1.55);
            \foreach \i in {0,...,5}{
            \foreach \j in {0,...,5}{
            \ifnum\j>\i \draw[gray!45,line width=.5pt] (c\i)--(c\j); \fi
            }
            }
            \foreach \i in {0,...,5}{
            \fill (c\i) circle[radius=1.8pt];
            \node[font=\scriptsize] at ($(0,0)!1.24!(c\i)$) {\(e_{\i}\)};
            }
            \foreach \i in {0,...,4} \draw[cyan!70!blue,line width=1.7pt] (c\i) circle[radius=4.4pt];
            \node[anchor=north] at (0,-2.12) {(b)};
        \end{scope}
        \draw[cyan!70!blue,line width=2.8pt] (1.22,-2.95)--(1.77,-2.95);
        \node[anchor=west,font=\scriptsize] at (1.9,-2.95) {selected};
        \draw[gray!75,line width=.65pt] (3.4,-2.95)--(3.95,-2.95);
        \node[anchor=west,font=\scriptsize] at (4.08,-2.95) {unselected};
    \end{tikzpicture}

    \caption{A convex hexagon has \(b(P)=5\). In (a), an attempted cover by four edges fails: the quadratic target \(t=(0,3/2-3\sqrt{2}/2)\) has \(X=0\) as its only rational line, so none sees the target at degree one. In (b), every edge pair admits a proper chord, giving the chord graph \(K_6\); five circled vertices form a minimum vertex cover.}
    \label{fig:chord-cover}
\end{figure}

\begin{figure}[h]
    \centering
    \begin{tikzpicture}[
        x=.86cm,y=.77cm,
        line cap=round,line join=round,
        bluearc/.style={draw=blue!70!black,line width=1.2pt},
        orangearc/.style={draw=orange!85!black,line width=1.2pt}
    ]
        \pgfmathsetmacro{\notchw}{.035}

        \fill[gray!9] (-1,-2.3) rectangle (11,2.3);
        \foreach \h in {0,2,4,6,8,10}{
            \fill[white]
                (\h-\notchw,-2.31) rectangle (\h+\notchw,.5);
            \draw[thick]
                (\h-\notchw,-2.3)--(\h-\notchw,.5)
                --(\h+\notchw,.5)--(\h+\notchw,-2.3);
        }
        \foreach \h in {1,3,5,7,9}{
            \fill[white]
                (\h-\notchw,-.5) rectangle (\h+\notchw,2.31);
            \draw[thick]
                (\h-\notchw,2.3)--(\h-\notchw,-.5)
                --(\h+\notchw,-.5)--(\h+\notchw,2.3);
        }
        \draw[thick]
            (-1,-2.3)--(-1,2.3)--(1-\notchw,2.3)
            (9+\notchw,2.3)--(11,2.3)--(11,-2.3)
            --(10+\notchw,-2.3)
            (-1,-2.3)--(-\notchw,-2.3);
        \foreach \h in {0,2,4,6,8}{
            \draw[thick]
                (\h+\notchw,-2.3)--(\h+2-\notchw,-2.3);
        }
        \foreach \h in {1,3,5,7}{
            \draw[thick]
                (\h+\notchw,2.3)--(\h+2-\notchw,2.3);
        }

        \draw[densely dashed,gray] (-1,0)--(11,0);
        \node[anchor=east,font=\scriptsize] at (-1.05,0)
            {\(Y=0\)};
        \node[anchor=east,font=\scriptsize] at (-1.05,.5)
            {\(\frac12\)};
        \node[anchor=east,font=\scriptsize] at (-1.05,-.5)
            {\(-\frac12\)};
        \node[anchor=east,font=\scriptsize] at (-1.05,2.3)
            {\(Y_6\)};
        \node[anchor=east,font=\scriptsize] at (-1.05,-2.3)
            {\(-Y_6\)};
        \foreach \x in {-1,11}{
            \draw[white,line width=3pt]
                (\x-.1,1.75)--(\x+.1,1.95);
            \draw[black,thin]
                (\x-.1,1.75)--(\x+.1,1.95);
            \draw[white,line width=3pt]
                (\x-.1,-1.95)--(\x+.1,-1.75);
            \draw[black,thin]
                (\x-.1,-1.95)--(\x+.1,-1.75);
        }

        \draw[bluearc,domain=0:2,samples=70,variable=\u]
            plot (\u,{2*(\u-1)^2-1});
        \draw[orangearc,domain=0:2,samples=70,variable=\u]
            plot ({8+\u},{2*(\u-1)^2-1});

        \draw[bluearc,domain=0:4,samples=160,variable=\u]
            plot ({2+\u},
                {(2*(\u-2)^4-8*(\u-2)^2+3)/3});
        \draw[orangearc,domain=0:4,samples=160,variable=\u]
            plot ({4+\u},
                {(2*(\u-2)^4-8*(\u-2)^2+3)/3});

        \foreach \i in {0,...,5}{
            \fill (2*\i,1) circle[radius=1.8pt];
            \node[above,font=\small] at (2*\i,1.12)
                {\(q_{\i}\)};
        }
        \draw[blue!70!black,very thick]
            (1.88,.88) rectangle (2.12,1.12);
        \draw[orange!85!black,very thick]
            (7.88,.88) rectangle (8.12,1.12);
        \node[above,font=\small\bfseries] at (5,2.65)
            {\(m=6,\quad d=4,\quad H=7683\)};
        \node[below,font=\scriptsize] at (-1,-2.4)
            {\(X=-1\)};
        \node[below,font=\scriptsize] at (11,-2.4)
            {\(X=11\)};
        \foreach \i in {0,...,10}{
            \node[below,font=\tiny] at (\i,-2.4) {\(\i\)};
        }

        \begin{scope}[yshift=-3.65cm]
            \draw[gray!65,thick] (0,0)--(10,0);
            \draw[bluearc] (0,0)--(4,0);
            \draw[orangearc] (6,0)--(10,0);

            \draw[gray!65,thick]
                (0,0) .. controls (1,-.85) and (3,-.85) .. (4,0);
            \draw[bluearc]
                (2,0) .. controls (3,.85) and (5,.85) .. (6,0);
            \draw[orangearc]
                (4,0) .. controls (5,-.85) and (7,-.85) .. (8,0);
            \draw[gray!65,thick]
                (6,0) .. controls (7,.85) and (9,.85) .. (10,0);

            \foreach \i in {0,...,5}{
                \fill (2*\i,0) circle[radius=1.8pt];
                \node[below,font=\scriptsize] at (2*\i,-.08)
                    {\(q_{\i}\)};
            }
            \draw[blue!70!black,very thick]
                (1.88,-.12) rectangle (2.12,.12);
            \draw[orange!85!black,very thick]
                (7.88,-.12) rectangle (8.12,.12);
            \node[below,font=\small] at (5,-1.02) {\(P_6^2\)};
        \end{scope}
    \end{tikzpicture}

    \caption{Degree-four sight in the alternating comb. The squares at \(q_1,q_4\) give a minimum two-guard cover with guards chosen from the targets, using two quadratic and two quartic arcs. Only the carriers used by the selected guards are drawn above. The graph below is \(P_6^2\); colored edges show sight from the selected guards. Every graph edge is available at coefficient height \(7683\), far below the general bound \(H_{6,4}\). The rectangle's height is compressed and notch widths exaggerated.}
    \label{fig:alternating-comb}
\end{figure}
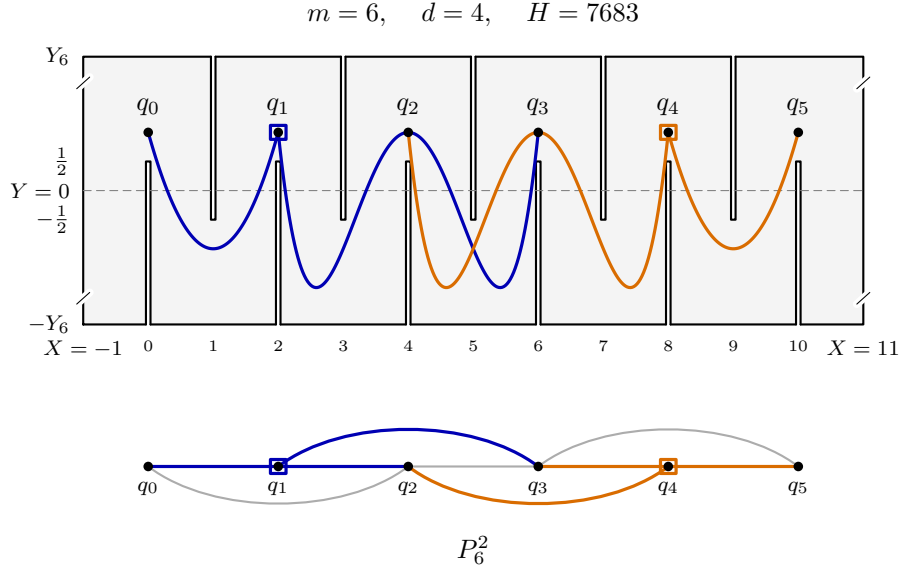

\Cref{fig:chord-cover} shows why the chord condition is necessary, even in a convex polygon.

The dependence on \(B^2\) in \cref{prop:edge-chords} is necessary. In the square, some interior quadratic targets of height at most \(B\) have a unique rational line whose coefficient height is of order \(B^2\). If the coefficient budget excludes that line, even selecting every edge would not cover the target.

\begin{proposition}\label{prop:edge-quadratic-height}
    Let \(K=\Q(\sqrt u)\subset\R\), with \(u>1\) squarefree. Let \(H_{\mathrm{cov}}(B)\) be the least coefficient-height budget at which some collection of edges covers \(T_{K,B}(\mathcal{S})\) at degree one. For all sufficiently large \(B\),
    \[
         \left\lfloor \frac{B^2}{2 \sqrt u} \right\rfloor \le H_{\mathrm{cov}}(B) \le 2B^2.
    \]
    Thus \(H_{\mathrm{cov}}(B)=\Theta_K(B^2)\). This concerns the existence of an edge cover, rather than its size once a cover exists.
\end{proposition}

For the lower bound, the proof in \cref{sec:edge-deferred} constructs an interior target of height at most \(B\) whose only rational line has coefficient height \(\lfloor B^2/(2\sqrt u)\rfloor\). Any smaller coefficient budget leaves this target invisible from every edge.

\begin{corollary}\label{cor:edge-convex}
    Let \(K\) be a real quadratic field and \(P\) a convex rational \(n\)-gon. Then \(b(P)=n-1\), and
    \[
        \eta_{1,\infty}(P;P(K)) = n-1,\quad \eta_{2,\infty}(P;P(K)) = 1.
    \]
    For all sufficiently large \(B\), these same counts hold for \(T_{K,B}(P)\) at every \(H \in [H_2(P,K,B),\infty]\), for some finite \(H_2\ge 2B^2\). Any single edge covers the target set at degree two.
\end{corollary}

\begin{proof}
    Every distinct edge pair of a convex polygon admits a proper chord, so its chord graph is complete, whence \(b(P)=n-1\) by \cref{prop:edge-chords}. By \cref{lem:convex-connector}, a rational midpoint on any edge sees every target in \(P(K)\) at degree two and unrestricted height. Thus any edge covers the whole quadratic gallery. For each edge midpoint and each target in the finite slice, choose a carrier supplied by \cref{lem:convex-connector}. Choose \(H_2\) as the maximum of \(2B^2\) and the coefficient heights of these finitely many carriers.
\end{proof}

In a convex polygon, classical segment visibility requires only one edge. For quadratic targets and rational carriers, the minimum is instead \(n-1\) edges at degree one, but returns to one at degree two. For a general rational polygon, we can use the degree-two vertex guards obtained from the triangulation argument and select an incident edge for each guard, giving the following bound.

\begin{corollary}\label{cor:edge-fisk}
    For a rational simple \(n\)-gon and a real quadratic field \(K\),
    \[
        \eta_{2,\infty}(P;P(K)) \le \lfloor n/3\rfloor.
    \]
    The same edges cover every finite slice at some common finite height.
\end{corollary}

\begin{proof}
    Select one edge incident to each rational vertex guard in \cref{cor:point-quadratic}. Those edges cover \(P(K)\). For a finite target slice, the maximum of the chosen carrier heights gives a common finite coefficient budget.
\end{proof}

The triangulation bound need not attain the classical minimum. We leave open whether every rational polygon \(P\) admits a real quadratic field \(K\) such that
\[
    \eta_{2,\infty}(P;P(K)) \le \eta_{\R}(P).
\]

\subsection{Alternating obstacles and degree bounds}
The one-sided comb gave degree-two sight room to pass through a common base, but alternating notches constrain that freedom. In the construction below, a sight arc between distant targets must repeatedly cross \(Y=0\), and the carrier's degree bounds the number of crossings. The height of the rectangle accommodates the interpolants' amplitude, while the notch width is chosen from bounds on their derivatives. The lower bound uses only alternation.

For \(m \ge 1\), put 
\[
    \Lambda_m = (4m)^{2m},\qquad \delta_m = 1/(4\Lambda_m), \qquad Y_m = (2m+1)\Lambda_m.
\]
From \({[-1,2m-1] \times [-Y_m,Y_m]}\) remove the relatively open notches
\[
    \begin{array}{ll}
         (h-\delta_m,h+\delta_m) \times [-Y_m,1/2),&h=0,2, \ldots, 2m-2,\\
         (h-\delta_m,h+\delta_m) \times (-1/2,Y_m],&h=1,3, \ldots, 2m-3.
    \end{array}
\]
The second family is empty when \(m=1\). Call the remaining polygon \(\mathcal{A}_m\). It has \(n=8m\) vertices. Let \(q_i=(2i,1)\) and \(T_m=\{q_i:0 \le i<m\}\).

\begin{theorem}\label{thm:alternating-comb}
    Let \(m,d \ge 1\), \(s = \min \{\lfloor d/2 \rfloor, m-1\}\), and put \(H_{m,d} = 1\) if \(s = 0\), and \(H_{m,d} = (2s+1)(2s)!(4m)^{2s}\) otherwise. At every \(H \in [1,\infty]\),
    \[
        q_i \sim_{d,H}q_j \Longrightarrow \abs{i-j} \le s, \qquad \gamma_{d,H}(\mathcal{A}_m;T_m,T_m) \ge \left\lceil \frac{m}{2s+1} \right\rceil. 
    \]
    If \(H \in[H_{m,d},\infty]\), then
    \[
        q_i \sim_{d,H}q_j \Longleftrightarrow \abs{i-j} \le s,
    \]
    and
    \[
        \gamma_{d,H}(\mathcal{A}_m;T_m,T_m) = \left\lceil \frac{m}{2s+1} \right\rceil = \Theta(\max\{1,n/d\}).
    \]
\end{theorem}

At \(H \ge H_{m,d}\), the visibility graph joins targets at most \(s\) positions apart, so it is the \(s\)th power of a path. Indeed, \(s=0\) gives no edges, \(s=1\) gives the path, and \(s=2\) also joins targets two steps apart. The proof bounds \(2\abs{i-j}\) distinct crossings of \(Y=0\), then constructs all permitted arcs by Newton interpolation. Note that the crossing bound applies even when the defining polynomial has arbitrary real coefficients. To attain the bound with admissible carriers, we construct polynomial graphs with rational coefficients, then clear denominators and bound the coefficients of their primitive integer equations. Details are in \cref{sec:comb-deferred}. 

\begin{corollary}\label{cor:alternating-free-guards}
    For \(G = \mathcal{A}_m\) or \(\mathcal{A}_m(\Q)\), at every \(H \in [1,\infty]\),
    \[
        \gamma_{d,H}(\mathcal{A}_m;T_m,G) \ge \left\lceil \frac{m}{d+2} \right\rceil.
    \]
    For \(H \in [H_{m,d},\infty]\) the upper bound is \(\lceil m/(2s+1) \rceil\), so this count also has order \(\Theta(\max\{1,n/d\})\).
\end{corollary}

The crossing argument from a guard \(g=(x,y)\) gives
\begin{equation}\label{eq:alternating-reach}
    \abs{x-2i} \le d+1
\end{equation}
for every visible target \(q_i\), regardless of its carrier. The resulting interval contains at most \(d+2\) targets. This proves the lower bound; the preceding target guards give the upper bound.

Below we find that edge guards change the bounds on the guard number but retain the growth order. 

Put \(Q_m=4\Lambda_m\) and
\[
    H^{\mathrm{edge}}_{m,d} = \max \left\{ 2mQ_m+2,4(2s)!(2s+1)(2m+1)^{2s} \prod_{k=1}^{2s}(Q_mk-1) \right\}.
\]
\begin{proposition}\label{prop:edge-alternating-sparse}
At every \(H \in [1,\infty]\), an edge cover of \(T_m\) needs at least \(\lceil m/(d+3) \rceil\) edges. For \(H \in [H^{\mathrm{edge}}_{m,d},\infty]\),
\[
 \left \lceil \frac{m}{d+3} \right\rceil \le \eta_{d,H}(\mathcal{A}_m;T_m) \le \left\lceil \frac{m}{s+1} \right\rceil,
 \quad \eta_{d,H}(\mathcal{A}_m;T_m) = \Theta(\max\{1,n/d\}).
\]
\end{proposition}

\begin{figure}[tbp]
    \centering
    \begin{tikzpicture}[
        x=1.1cm,y=.8cm,font=\small,
        line cap=round,line join=round
        ]
        \def\notchw{.08}

        \fill[gray!8] (-1,-2.1) rectangle (7,2.1);

        \foreach \h in {0,2,4,6}{
        \fill[white]
          (\h-\notchw,-2.11) rectangle (\h+\notchw,.5);
        \draw[gray!75,thick]
          (\h-\notchw,-2.1)--(\h-\notchw,.5)
          --(\h+\notchw,.5)--(\h+\notchw,-2.1);
        }

        \foreach \h in {1,3,5}{
        \fill[white]
          (\h-\notchw,-.5) rectangle (\h+\notchw,2.11);
        \draw[gray!75,thick]
          (\h-\notchw,2.1)--(\h-\notchw,-.5)
          --(\h+\notchw,-.5)--(\h+\notchw,2.1);
        }

        \draw[gray!75,thick]
        (-1,-2.1)--(-1,2.1)--(1-\notchw,2.1)
        (1+\notchw,2.1)--(3-\notchw,2.1)
        (3+\notchw,2.1)--(5-\notchw,2.1)
        (5+\notchw,2.1)--(7,2.1)--(7,-2.1)
        --(6+\notchw,-2.1)
        (6-\notchw,-2.1)--(4+\notchw,-2.1)
        (4-\notchw,-2.1)--(2+\notchw,-2.1)
        (2-\notchw,-2.1)--(\notchw,-2.1)
        (-\notchw,-2.1)--(-1,-2.1);

        \draw[gray,densely dashed] (-1,0)--(7,0);
        \node[anchor=east,font=\scriptsize]
        at (-1.03,0) {\(Y=0\)};

        \pgfmathsetmacro{\quadA}
        {(2+1.5/(1-\notchw))/(2-\notchw)}

        \foreach \i/\clr in 
        {0/cyan!70!blue,2/orange!85!black}{
        \draw[draw=\clr,line width=1.15pt]
          (2*\i-\notchw,.5)--(2*\i+\notchw,.5);

        \draw[draw=\clr,line width=1.15pt]
          (2*\i+\notchw,.5)--(2*\i,1);

        \draw[
          draw=\clr,line width=1.15pt,
          domain=\notchw:2,samples=80,variable=\u
        ]
          plot (
            {2*\i+\u},
            {\quadA*\u*\u+(2-3*\quadA)*\u+2*\quadA-3}
          );

        \fill[fill=\clr]
          (2*\i+\notchw,.5) circle[radius=1.8pt];
        \node[below right,font=\scriptsize,text=\clr]
          at (2*\i+\notchw+.02,.45) {\(v_{\i}\)};
        }

        \foreach \i in {0,...,3}{
        \fill (2*\i,1) circle[radius=1.8pt];
        \node[above] at (2*\i,1.12) {\(q_{\i}\)};
        }

        \node[above] at (3,2.22)
        {\(m=4,\quad d=2,\quad s=1\)};
    \end{tikzpicture}

    \caption{The upper-bound construction in \cref{prop:edge-alternating-sparse}, with blocks \(\{q_0,q_1\}\) and \(\{q_2,q_3\}\).
    The selected bottom-notch tip edges are colored, with their right endpoints \(v_i\) marked by dots. Each \(v_i\) sees \(q_i\) by a segment and \(q_{i+1}\) by a quadratic sight arc. The rectangle's height is compressed and notch widths exaggerated.}
    \label{fig:alternating-edges}
\end{figure}
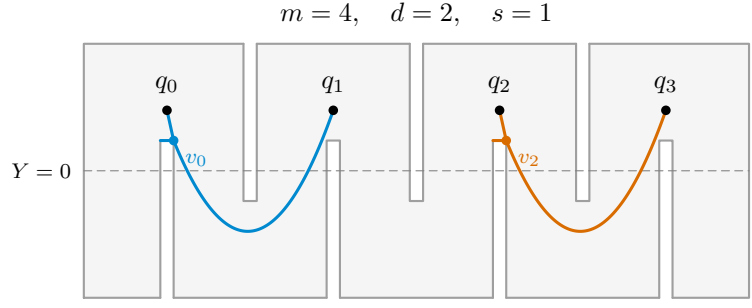

Each edge has horizontal span at most two, so \cref{eq:alternating-reach} limits it to \(d+3\) targets. For the upper bound divide the targets into blocks of at most \(s+1\). A bottom-notch tip edge sees its block from its right endpoint, using a segment and interpolating graphs; see \cref{fig:alternating-edges}. The coefficient check is in \cref{sec:edge-deferred}.

Now let us ask edge guards to cover all rational points, not only \(T_m\). We can still use \(T_m\) to arrive at a lower bound. A different edge collection uses horizontal sight arcs to cover the whole rational gallery.

\begin{theorem}\label{thm:edge-alternating}
Let \(U_B=T_{\Q,B}(\mathcal{A}_m)\). For \(d\ge1\), \(B\ge\max\{1,2m-2\}\), and \(H\in[B,\infty]\),
\[
    \left\lceil \frac{m}{d+3} \right\rceil \le \eta_{d,H}(\mathcal{A}_m;U_B) \le 2m.
\]
The lower bound holds at every coefficient height. The same bounds hold for \(\mathcal{A}_m(\Q)\) at unrestricted height. Thus both counts are \(\Theta_d(n)\) for fixed \(d\) in the stated ranges.
\end{theorem}

\begin{proof}
The slice contains \(T_m\), establishing the lower bound by \cref{prop:edge-alternating-sparse}. Select the left vertical edge of \(\mathcal{A}_m\) and the right vertical face of every notch, \(2m\) edges in all. Every horizontal component has its left end on one of them. For a rational target of ordinate \(a/b\) in lowest terms, its horizontal line \(bY-a=0\) has height at most \(B\), whence we have a cover. The same edges work for all rational targets at unrestricted height.
\end{proof}

\section{Sight around a corner}\label{sec:l-gallery}

We now use an L-shaped gallery to examine how degree, guard placement, and coefficient field affect visibility. Throughout this section let
\[
    \mathcal{L} = ([0,2]\times[0,1])\cup([0,1]\times[0,2]).
\]
First we compare edge guard numbers at degrees one and two with rational carriers. We then fix a guard outside the classical kernel and show that rational conics reach an explicit family of quadratic targets, whereas allowing coefficients in the target field gives a quadratic cover of the whole gallery.

\begin{example}\label{ex:edge-l-shape}
    Suppose \(K\) is a real quadratic field. The two long edges of \(\mathcal{L}\), at the bottom and left, admit proper chords to every other edge. Among the four short edges, the two bounding each arm admit a chord to one another, but no chord joins an edge from one arm to the other. Such a segment would cross \(x=y>1\), outside \(\mathcal{L}\). The largest set of pairwise nonadjacent vertices in the chord graph has size two. Since the complement of a largest pairwise nonadjacent set is a minimum vertex cover, \(b(\mathcal L)=6-2=4\). By \cref{prop:edge-chords},
    \[
        \eta_{1,\infty}(\mathcal{L};\mathcal{L}(K))=4.
    \]
    Every target belongs to one of the two rational rectangles, both of which contain \((0,0)\). Applying \cref{lem:convex-connector} in that rectangle shows that either long edge, observing from \((0,0)\), sees the entire \(K\)-gallery at degree two. Hence
    \[
        \eta_{2,\infty}(\mathcal{L};\mathcal{L}(K))=1.
    \]
    For \(T_{K,B}(\mathcal L)\), the same counts are \(4\) and \(1\) once \(B\) and the common coefficient-height budget are sufficiently large.
\end{example}

This one-edge cover uses a point \((0,0)\) in the classical kernel. Can a rational guard outside of that kernel see around the corner to cover a target set? The following conics answer this for an explicit family without changing the coefficient field.

\begin{proposition}\label{prop:l-rational-conic}
    Fix \(g=(3/2,1/2)\). For a rational nonsquare \(r=a/b\) in lowest terms, with \(b>0\) and \(1/4 < r \le 3/4\), put \(t_r=(1/2,1+\sqrt{r})\). There is a regular sight arc from \(g\) to \(t_r\) inside \(\mathcal{L}\) on an absolutely irreducible rational conic of primitive integer coefficient height at most \(24b\). Degree two is the least possible degree for this pair, even without a coefficient-height bound. At \(r=1/2\) one such primitive equation is
    \[
        8Y^2+16XY-6X-24Y+7 = 0,
    \]
    of height \(24\).
\end{proposition}

\begin{proof}
    Since \(t_r=(1/2,1+\sqrt{r})\), the polynomial \((Y-1)^2-r\) vanishes at the target. Adding a term divisible by \(X-1/2\) preserves this vanishing, so we seek a conic of the form
    \[
        F_r(X,Y)=(X-1/2)L_r(Y)+(Y-1)^2-r,
    \]
    where \(L_r\) is affine linear. To make \(F_r\) vanish at \(g=(3/2,1/2)\), we need \(L_r(1/2)=r-1/4\). We choose \(L_r(Y)=2Y+r-5/4\), giving
    \[
        F_r(X,Y)=(X-1/2)(2Y+r-5/4)+(Y-1)^2-r.
    \]
    This choice makes the connecting arc stay in the two arms, as we now verify. The polynomial \(F_r\) vanishes at \(g\) and \(t_r\). On \(1/2 \le Y \le 1+\sqrt{r}\) its zero set is the graph
    \[
        X=\frac{1}{2}+\frac{r-(Y-1)^2}{2Y+r-5/4},
    \]
    whose denominator is at least \(r-1/4>0\). For \(Y=1/2+s\), with \(0 \le s \le 1/2\), its numerator is \(r-1/4+s-s^2\), between zero and the denominator \(r-1/4+2s\). Thus \(1/2 \le X \le 3/2\), keeping this part in the lower arm. For \(1 \le Y \le 1+\sqrt{r}\), the numerator lies in \([0,r]\) and the denominator is at least \(r+3/4\). Consequently \(0 \le X-1/2 \le r/(r+3/4) \le 1/2\), keeping this part in the upper arm. The graph is embedded, and \(F_{r,X}=2Y+r-5/4>0\) on it.

    With \(c=r-5/4\), the symmetric matrix of the homogenized conic is
    \[
        \begin{pmatrix}
            0 & 1 & c/2\\
            1 & 1 & -3/2\\
            c/2 & -3/2 & 13/8-3r/2
        \end{pmatrix}.
    \]
    Its determinant is \((1-c^2)/4>0\). The projective conic is therefore nonsingular and absolutely irreducible. Multiplying \(F_r\) by \(8b\) gives the integer polynomial
    \[
        16bXY+(8a-10b)X+8bY^2-24bY+13b-12a,
    \]
    whose height is \(24b\). Taking its primitive part can only decrease this bound, and gives the displayed height-\(24\) equation when \(r=1/2\).

    Finally, parameterize the segment by \((1-\lambda)g+\lambda t_r\). For \(1/(1+2\sqrt{r})< \lambda <1 /2\), both coordinates exceed one, so this nonempty portion lies outside \(\mathcal{L}\). Every arc on the line through the two points contains their segment. Degree one is therefore impossible.
\end{proof}

The coefficient bound has quadratic order in the target height. Indeed, \(\height(\sqrt{r})=\sqrt{b}\), and the usual sum and product bounds give
\[
    \frac{\sqrt{b}}{2}\le\height(t_r)\le4\sqrt{b}.
\]
Thus the sufficient bound \(24b\) is \(O(\height(t_r)^2)\). This bounds the height of the constructed carrier, not the minimum height of a conic.

\begin{figure}[tb]
    \centering
    \begin{tikzpicture}[x=2.05cm,y=2.05cm,line cap=round,line join=round]
        \foreach \shift/\panel in {0/a,5.6/b}{
        \begin{scope}[xshift=\shift cm]
            \fill[gray!8] (0,0)--(2,0)--(2,1)--(1,1)--(1,2)--(0,2)--cycle;
            \draw[thick] (0,0)--(2,0)--(2,1)--(1,1)--(1,2)--(0,2)--cycle;
            \draw[red!60!black,densely dashed] (1.5,.5)--(.5,{1+sqrt(.5)});
            \fill (1.485,.485) rectangle (1.515,.515);
            \node[below right,font=\small] at (1.5,.5) {\(g\)};
            \fill (.5,{1+sqrt(.5)}) circle[radius=2pt];
            \node[above left,font=\small] at (.5,{1+sqrt(.5)}) {\(t_{1/2}\)};
            \node[below,font=\small] at (1,-.12) {(\panel)};
        \end{scope}
        }
        \draw[cyan!65!blue,very thick,domain=.5:1.7071067812,samples=90,variable=\y]
        plot ({.5+(.5-(\y-1)*(\y-1))/(2*\y-.75)},\y);
        \begin{scope}[xshift=5.6cm]
            \draw[gray!65,densely dotted] (1.5,.5)--(.5,.5)--(.5,{1+sqrt(.5)});
            \draw[orange!85!black,very thick,domain=0:1,samples=80,variable=\s]
            plot ({1.5-2*\s+\s*\s},{.5+(.5+sqrt(.5))*\s*\s});
            \draw[fill=white] (.5,.5) circle[radius=2pt];
            \node[below left,font=\small] at (.5,.5) {\(c\)};
        \end{scope}
    \end{tikzpicture}

    \caption{The same guard and quadratic target in \(\mathcal{L}\). (a) The rational conic \(Z(8Y^2+16XY-6X-24Y+7)\) gives a contained sight arc of height \(24\). (b) A quadratic B\'ezier curve gives another contained arc, here with coefficients in \(\Q(\sqrt{2})\). This construction works for every target over its coefficient field, as shown in \cref{ex:field-kernel-l}. The dashed segment leaves the polygon.}
    \label{fig:l-rational-conic}
\end{figure}
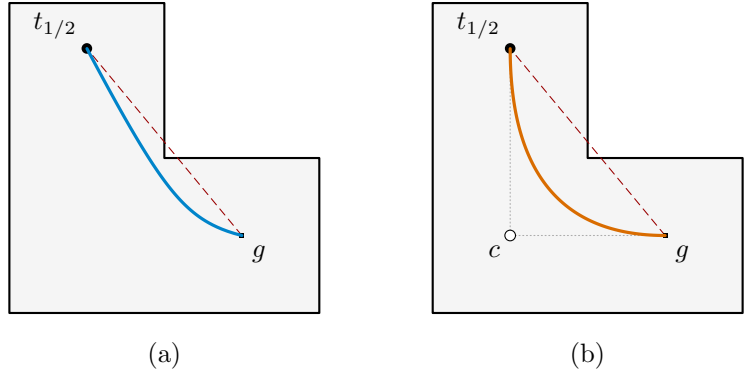

For comparison, now allow the carrier coefficients to lie in a real number field \(K\). In this comparison we retain rational guards and take \(H=\infty\), isolating the effect of changing the carrier field.

\begin{definition}[\(K\)-sight]\label{def:k-sight}
    Points \(p,q \in P\) satisfy \(p \sim_{d,H}^K q\) when the points coincide or lie in one connected regular trace component of a sight carrier defined by an absolutely irreducible \(F \in K[X,Y]\). This polynomial must have degree at most \(d\) and coefficient height at most \(H\), measured by the absolute multiplicative Weil height \(\height([\operatorname{coeff}(F)])\) of the coefficient vector. If this relation is satisfied, we say that \(p\) sees or guards \(q\) over \(K\). We omit \(H\) when \(H = \infty\).

When \(K=\Q\), this agrees with \cref{def:visibility}: every rational defining polynomial has a primitive integer representative whose coefficient height, measured as above, equals its largest absolute coefficient; see \cref{rem:coefficient-field}.
\end{definition}

\begin{definition}[Arithmetic kernel]\label{def:arith-kernel}
    The \emph{arithmetic kernel} of the gallery \((P;T,G)\) at budget \(d,H\) over \(K\) is
    \[
        \operatorname{Ker}_{d,H}(P,K;T,G)=\{g\in G : g\sim^K_{d,H}t \text{ for every } t\in T\}.
    \]
    We omit \(K\) when \(K=\Q\), omit \(H\) when \(H=\infty\), and omit \(T,G\) when \(T=P(K)\) and \(G=P(\Q)\). For nonempty \(T\), the kernel is nonempty exactly when one guard from \(G\) covers \(T\).
\end{definition}

No primitivity condition is needed since the coefficient height is unchanged by multiplying the defining polynomial by a nonzero element of \(K\). At degree one and unrestricted coefficient height,
\[
    \operatorname{Ker}_1(P,K) = \operatorname{Ker}(P) \cap \Q^2.
\]
A rational guard in the classical kernel sees every \(K\)-point by a \(K\)-line. Conversely, a guard in \(\operatorname{Ker}_1(P,K)\) sees every rational target; density of rational targets and closedness of \(P\) imply that it sees every real target.

\begin{example}\label{ex:field-kernel-l}
    For every real number field \(K\), the guard \(g=(3/2,1/2)\) satisfies
    \[
        g \notin \operatorname{Ker}_1(\mathcal{L},K),\qquad g \in \operatorname{Ker}_2(\mathcal{L},K).
    \]
    The segment to \((1,2)\) leaves \(\mathcal{L}\). Lower-arm targets are visible by lines over \(K\). For an upper-arm target \(t=(u,v)\), with \(0\le u\le1\) and \(1<v\le2\), use the quadratic B\'ezier curve with control point \(c=(1/2,1/2)\):
    \[
        \beta_t(s)=(1-s)^2g+2s(1-s)c+s^2t,\qquad 0 \le s \le 1.
    \]
    Its coordinates are \(X=3/2-2s+(u+1/2)s^2\) and \(Y=1/2+(v-1/2)s^2\). The nonnegative coefficients sum to one, so the curve stays in \([0,3/2]\times[1/2,2]\). For \(s \le 1/2\) we have \(Y \le 7/8<1\). For \(s \ge 1/2\), we have \(X\le3/2-2s+(3/2)s^2 \le 1\). These two parts lie in the lower and upper arms, respectively.

    Put \(a=v-1/2>0\) and \(b=u+1/2\). Its carrier is
    \[
        Z\!\left( [a(X-3/2)-b(Y-1/2)]^2-4a(Y-1/2) \right).
    \]
    The invertible affine coordinates \(U=a(X-3/2)-b(Y-1/2)\), \(V=Y-1/2\) give \(U^2-4aV=0\), an absolutely irreducible nonsingular parabola over \(K\). Along the arc \(U=-2as\), so it is embedded and regular.
\end{example}

In the last example the coefficients depend on the target and may lie outside \(\Q\). It therefore covers every \(K\)-point at degree two in the \(K\)-coefficient model, whereas \cref{prop:l-rational-conic} gives rational carriers and explicit heights for a particular family. Whether the same guard sees every quadratic point of the upper arm on a rational conic remains a further question; see \cref{prob:field-kernel}.

At a fixed degree, then, the choice of carrier field changes which targets can be reached. To isolate the arithmetic restrictions, we return to our rational carrier model and work in the unit square, where every joining segment stays inside the gallery.

\section{How much coefficient height does guarding cost?}\label{sec:scarcity}

Classically, one guard anywhere in the square sees the whole gallery. At small coefficient height, arithmetic visibility can require many guards even for a finite set of rational targets. We ask how much coefficient height is needed to recover one-guard coverage of the square's finite rational target-height slices.

Throughout this section the target field is \(\Q\), so we suppress it:
\[
    T_B(P)\coloneqq T_{\Q,B}(P),\qquad T_B\coloneqq T_B(\mathcal{S}).
\]

We first count the targets, then bound how many lie on a single curve. Since a fixed budget admits only finitely many carriers, we can compare the total number of targets with the number lying on any admissible carrier.

\begin{lemma}\label{lem:square-count}
    As \(B \to \infty\),
    \[
        \card{T_B}=\frac{B^3}{3\zeta(3)}+O(B^2),
    \]
    where \(\zeta\) is the Riemann zeta function.
\end{lemma}

\begin{proof}
    A target has a unique expression \(q=(a/c,b/c)\), where \(0\le a,b\le c\le B\), \(c\ge1\), and \(\gcd(a,b,c)=1\). These triples lie in a square-based cone in \(\R^3\), truncated by \(c\le B\); the slice at height \(c\) has \((c+1)^2\) integer pairs \((a,b)\). Without the coprimality restriction their number is
    \[
        N(B)\coloneqq\sum_{1\le c\le B}(c+1)^2=B^3/3+O(B^2).
    \]
    The identity \(\mathbf1_{\gcd(a,b,c)=1}=\sum_{j\mid a,b,c}\mu(j)\), where \(\mu\) is the M\"obius function and the indicator equals one if \(\gcd(a,b,c)=1\) and zero otherwise, gives
    \[
    \card{T_B}
     =\sum_{j\le B} \mu(j) N(B/j)
     =\frac{B^3}{3} \sum_{j \le B} \frac{\mu(j)}{j^3}
       +O\left(B^2 \sum_{j\le B} \frac1{j^2} \right).
    \]
    Now \(\sum_{j \ge 1} \mu(j)/j^3= 1/\zeta(3)\), with tail \(\sum_{j>B} \mu(j)/j^3= O(B^{-2})\). Multiplying this tail by \(B^3/3\) contributes only \(O(B)\), while the displayed error is \(O(B^2)\). This proves the estimate.
\end{proof}

\begin{lemma}\label{lem:occupancy}
    If \(0 \ne F \in \Q[X,Y]\) is irreducible of degree at most \(d\), then
    \[
        \card{\{q \in T_B : F(q)=0\}} = O(dB^2),
    \]
    with an absolute implied constant.
\end{lemma}

\begin{proof}
    The number of rational coordinates in \([0,1]\) of height at most \(B\) is \(1+\sum_{m \le B} \varphi(m)=O(B^2)\), where \(\varphi\) is Euler's totient function. Unless \(\mathrm{C}_F\) is a vertical line, \(F(x,Y)\) is nonzero for every rational \(x\), since otherwise \(X-x\) divides \(F\). There are then at most \(d\) choices of \(y\) for each \(x\). If \(F\) depends only on \(X\), it either has no rational root or defines a vertical line, on which there are \(O(B^2)\) possible \(y\)-coordinates.
\end{proof}

For geometrically integral projective plane curves of degree \(r\), Walsh's theorem gives the sharper uniform estimate \(O_r(B^{2/r})\) \cite[Theorem~1.1]{Walsh2015}. It builds on earlier determinant-method bounds for rational points \cite{HeathBrown2002} and for integral points \cite{BombieriPila1989}. The elementary bound is enough here because our family includes lines, which already contain \(\Theta(B^2)\) targets.

\begin{theorem}[Sparsity at a fixed carrier budget]\label{thm:sparsity}
    For \(d \ge 1\) and finite \(H \ge 1\),
    \[
        \card{T_B}-\gamma_{d,H}(\mathcal{S};T_B,T_B)=\Theta_{d,H}(B^2).
    \]
    Consequently,
    \[
        \frac{\gamma_{d,H}(\mathcal{S};T_B,T_B)}{\card{T_B}} \longrightarrow 1, \qquad B \to \infty.
    \]
\end{theorem}

\begin{proof}
By \cref{eq:curve-count}, the number of admissible carriers is a finite constant depending only on the carrier degree and coefficient height, \emph{not} the target height. Hence, by \cref{lem:occupancy}, the union of their regular traces contains \(O_{d,H}(B^2)\) targets. It follows that \(\card{T_B}-O_{d,H}(B^2)\) is a lower bound on the number of guards required, by \cref{cor:forced}. This proves the upper bound 
\(\card{T_B}-\gamma_{d,H}(\mathcal{S};T_B,T_B) \le O_{d,H}(B^2)\).

For the matching lower bound, the line \(Z(X)\) is admissible and contains \(1+\sum_{m \le B} \varphi(m)= \Theta(B^2)\) targets. One of them guards all targets on that side of the square; every remaining target may guard itself. The ratio assertion follows from \cref{lem:square-count} via
\[
    \frac{\gamma_{d,H}(\mathcal{S};T_B,T_B)}{\card{T_B}}
    = 1 - \frac{\card{T_B} - \gamma_{d,H}(\mathcal{S};T_B,T_B)}{\card{T_B}}
    = 1 - \Theta_{d,H}(B^{-1}) \longrightarrow 1.
\]
\end{proof}

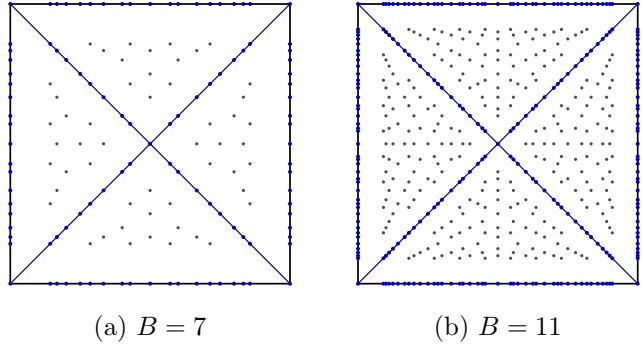
\begin{figure}[tbp]
    \centering
    \begin{tikzpicture}[x=3.7cm,y=3.7cm]
        \begin{scope}[xshift=0cm]
            \draw[black,line width=.65pt] (0,0) rectangle (1,1);
            \draw[blue!50!black,line width=.45pt] (0,0)--(1,1) (0,1)--(1,0);
            \foreach \ax/\by/\cz in {
              0/0/1,0/1/1,1/0/1,1/1/1,0/1/2,1/0/2,1/1/2,1/2/2,2/1/2,0/1/3,
              0/2/3,1/0/3,1/1/3,1/2/3,1/3/3,2/0/3,2/1/3,2/2/3,2/3/3,3/1/3,
              3/2/3,0/1/4,0/3/4,1/0/4,1/1/4,1/3/4,1/4/4,3/0/4,3/1/4,3/3/4,
              3/4/4,4/1/4,4/3/4,0/1/5,0/2/5,0/3/5,0/4/5,1/0/5,1/1/5,1/4/5,
              1/5/5,2/0/5,2/2/5,2/3/5,2/5/5,3/0/5,3/2/5,3/3/5,3/5/5,4/0/5,
              4/1/5,4/4/5,4/5/5,5/1/5,5/2/5,5/3/5,5/4/5,0/1/6,0/5/6,1/0/6,
              1/1/6,1/5/6,1/6/6,5/0/6,5/1/6,5/5/6,5/6/6,6/1/6,6/5/6,0/1/7,
              0/2/7,0/3/7,0/4/7,0/5/7,0/6/7,1/0/7,1/1/7,1/6/7,1/7/7,2/0/7,
              2/2/7,2/5/7,2/7/7,3/0/7,3/3/7,3/4/7,3/7/7,4/0/7,4/3/7,4/4/7,
              4/7/7,5/0/7,5/2/7,5/5/7,5/7/7,6/0/7,6/1/7,6/6/7,6/7/7,7/1/7,
              7/2/7,7/3/7,7/4/7,7/5/7,7/6/7
            } {\fill[blue!75!black] ({\ax/\cz},{\by/\cz}) circle (.75pt);}
            \foreach \ax/\by/\cz in {
              1/2/4,2/1/4,2/3/4,3/2/4,1/2/5,1/3/5,2/1/5,2/4/5,3/1/5,3/4/5,
              4/2/5,4/3/5,1/2/6,1/3/6,1/4/6,2/1/6,2/3/6,2/5/6,3/1/6,3/2/6,
              3/4/6,3/5/6,4/1/6,4/3/6,4/5/6,5/2/6,5/3/6,5/4/6,1/2/7,1/3/7,
              1/4/7,1/5/7,2/1/7,2/3/7,2/4/7,2/6/7,3/1/7,3/2/7,3/5/7,3/6/7,
              4/1/7,4/2/7,4/5/7,4/6/7,5/1/7,5/3/7,5/4/7,5/6/7,6/2/7,6/3/7,
              6/4/7,6/5/7
            } {\fill[black!70] ({\ax/\cz},{\by/\cz}) circle (.6pt);}
            \node[font=\small,anchor=north] at (.5,-.075) {(a) \(B=7\)};
        \end{scope}
        \begin{scope}[xshift=4.6cm]
            \draw[black,line width=.65pt] (0,0) rectangle (1,1);
            \draw[blue!50!black,line width=.45pt] (0,0)--(1,1) (0,1)--(1,0);
            \foreach \ax/\by/\cz in {
              0/0/1,0/1/1,1/0/1,1/1/1,0/1/2,1/0/2,1/1/2,1/2/2,2/1/2,0/1/3,
              0/2/3,1/0/3,1/1/3,1/2/3,1/3/3,2/0/3,2/1/3,2/2/3,2/3/3,3/1/3,
              3/2/3,0/1/4,0/3/4,1/0/4,1/1/4,1/3/4,1/4/4,3/0/4,3/1/4,3/3/4,
              3/4/4,4/1/4,4/3/4,0/1/5,0/2/5,0/3/5,0/4/5,1/0/5,1/1/5,1/4/5,
              1/5/5,2/0/5,2/2/5,2/3/5,2/5/5,3/0/5,3/2/5,3/3/5,3/5/5,4/0/5,
              4/1/5,4/4/5,4/5/5,5/1/5,5/2/5,5/3/5,5/4/5,0/1/6,0/5/6,1/0/6,
              1/1/6,1/5/6,1/6/6,5/0/6,5/1/6,5/5/6,5/6/6,6/1/6,6/5/6,0/1/7,
              0/2/7,0/3/7,0/4/7,0/5/7,0/6/7,1/0/7,1/1/7,1/6/7,1/7/7,2/0/7,
              2/2/7,2/5/7,2/7/7,3/0/7,3/3/7,3/4/7,3/7/7,4/0/7,4/3/7,4/4/7,
              4/7/7,5/0/7,5/2/7,5/5/7,5/7/7,6/0/7,6/1/7,6/6/7,6/7/7,7/1/7,
              7/2/7,7/3/7,7/4/7,7/5/7,7/6/7,0/1/8,0/3/8,0/5/8,0/7/8,1/0/8,
              1/1/8,1/7/8,1/8/8,3/0/8,3/3/8,3/5/8,3/8/8,5/0/8,5/3/8,5/5/8,
              5/8/8,7/0/8,7/1/8,7/7/8,7/8/8,8/1/8,8/3/8,8/5/8,8/7/8,0/1/9,
              0/2/9,0/4/9,0/5/9,0/7/9,0/8/9,1/0/9,1/1/9,1/8/9,1/9/9,2/0/9,
              2/2/9,2/7/9,2/9/9,4/0/9,4/4/9,4/5/9,4/9/9,5/0/9,5/4/9,5/5/9,
              5/9/9,7/0/9,7/2/9,7/7/9,7/9/9,8/0/9,8/1/9,8/8/9,8/9/9,9/1/9,
              9/2/9,9/4/9,9/5/9,9/7/9,9/8/9,0/1/10,0/3/10,0/7/10,0/9/10,1/0/10,
              1/1/10,1/9/10,1/10/10,3/0/10,3/3/10,3/7/10,3/10/10,7/0/10,7/3/10,7/7/10,
              7/10/10,9/0/10,9/1/10,9/9/10,9/10/10,10/1/10,10/3/10,10/7/10,10/9/10,0/1/11,
              0/2/11,0/3/11,0/4/11,0/5/11,0/6/11,0/7/11,0/8/11,0/9/11,0/10/11,1/0/11,
              1/1/11,1/10/11,1/11/11,2/0/11,2/2/11,2/9/11,2/11/11,3/0/11,3/3/11,3/8/11,
              3/11/11,4/0/11,4/4/11,4/7/11,4/11/11,5/0/11,5/5/11,5/6/11,5/11/11,6/0/11,
              6/5/11,6/6/11,6/11/11,7/0/11,7/4/11,7/7/11,7/11/11,8/0/11,8/3/11,8/8/11,
              8/11/11,9/0/11,9/2/11,9/9/11,9/11/11,10/0/11,10/1/11,10/10/11,10/11/11,11/1/11,
              11/2/11,11/3/11,11/4/11,11/5/11,11/6/11,11/7/11,11/8/11,11/9/11,11/10/11
            } {\fill[blue!75!black] ({\ax/\cz},{\by/\cz}) circle (.75pt);}
            \foreach \ax/\by/\cz in {
              1/2/4,2/1/4,2/3/4,3/2/4,1/2/5,1/3/5,2/1/5,2/4/5,3/1/5,3/4/5,
              4/2/5,4/3/5,1/2/6,1/3/6,1/4/6,2/1/6,2/3/6,2/5/6,3/1/6,3/2/6,
              3/4/6,3/5/6,4/1/6,4/3/6,4/5/6,5/2/6,5/3/6,5/4/6,1/2/7,1/3/7,
              1/4/7,1/5/7,2/1/7,2/3/7,2/4/7,2/6/7,3/1/7,3/2/7,3/5/7,3/6/7,
              4/1/7,4/2/7,4/5/7,4/6/7,5/1/7,5/3/7,5/4/7,5/6/7,6/2/7,6/3/7,
              6/4/7,6/5/7,1/2/8,1/3/8,1/4/8,1/5/8,1/6/8,2/1/8,2/3/8,2/5/8,
              2/7/8,3/1/8,3/2/8,3/4/8,3/6/8,3/7/8,4/1/8,4/3/8,4/5/8,4/7/8,
              5/1/8,5/2/8,5/4/8,5/6/8,5/7/8,6/1/8,6/3/8,6/5/8,6/7/8,7/2/8,
              7/3/8,7/4/8,7/5/8,7/6/8,1/2/9,1/3/9,1/4/9,1/5/9,1/6/9,1/7/9,
              2/1/9,2/3/9,2/4/9,2/5/9,2/6/9,2/8/9,3/1/9,3/2/9,3/4/9,3/5/9,
              3/7/9,3/8/9,4/1/9,4/2/9,4/3/9,4/6/9,4/7/9,4/8/9,5/1/9,5/2/9,
              5/3/9,5/6/9,5/7/9,5/8/9,6/1/9,6/2/9,6/4/9,6/5/9,6/7/9,6/8/9,
              7/1/9,7/3/9,7/4/9,7/5/9,7/6/9,7/8/9,8/2/9,8/3/9,8/4/9,8/5/9,
              8/6/9,8/7/9,1/2/10,1/3/10,1/4/10,1/5/10,1/6/10,1/7/10,1/8/10,2/1/10,
              2/3/10,2/5/10,2/7/10,2/9/10,3/1/10,3/2/10,3/4/10,3/5/10,3/6/10,3/8/10,
              3/9/10,4/1/10,4/3/10,4/5/10,4/7/10,4/9/10,5/1/10,5/2/10,5/3/10,5/4/10,
              5/6/10,5/7/10,5/8/10,5/9/10,6/1/10,6/3/10,6/5/10,6/7/10,6/9/10,7/1/10,
              7/2/10,7/4/10,7/5/10,7/6/10,7/8/10,7/9/10,8/1/10,8/3/10,8/5/10,8/7/10,
              8/9/10,9/2/10,9/3/10,9/4/10,9/5/10,9/6/10,9/7/10,9/8/10,1/2/11,1/3/11,
              1/4/11,1/5/11,1/6/11,1/7/11,1/8/11,1/9/11,2/1/11,2/3/11,2/4/11,2/5/11,
              2/6/11,2/7/11,2/8/11,2/10/11,3/1/11,3/2/11,3/4/11,3/5/11,3/6/11,3/7/11,
              3/9/11,3/10/11,4/1/11,4/2/11,4/3/11,4/5/11,4/6/11,4/8/11,4/9/11,4/10/11,
              5/1/11,5/2/11,5/3/11,5/4/11,5/7/11,5/8/11,5/9/11,5/10/11,6/1/11,6/2/11,
              6/3/11,6/4/11,6/7/11,6/8/11,6/9/11,6/10/11,7/1/11,7/2/11,7/3/11,7/5/11,
              7/6/11,7/8/11,7/9/11,7/10/11,8/1/11,8/2/11,8/4/11,8/5/11,8/6/11,8/7/11,
              8/9/11,8/10/11,9/1/11,9/3/11,9/4/11,9/5/11,9/6/11,9/7/11,9/8/11,9/10/11,
              10/2/11,10/3/11,10/4/11,10/5/11,10/6/11,10/7/11,10/8/11,10/9/11
            } {\fill[black!70] ({\ax/\cz},{\by/\cz}) circle (.6pt);}
            \node[font=\small,anchor=north] at (.5,-.075) {(b) \(B=11\)};
        \end{scope}
    \end{tikzpicture}

    \caption{Rational targets in the same square at \(d=H=1\). Blue points lie on the four sides or the two diagonals, which together contain all admissible regular traces at this budget. Black points lie on no admissible carrier and must guard themselves. At \(B=7\) there are \(105\) blue and \(52\) black points; at \(B=11\) there are \(249\) blue and \(268\) black points.}
    \label{fig:scarcity-targets}
\end{figure}

With a larger coefficient budget, two separate questions arise. The first is whether every target lies on an admissible carrier, and the second is whether a few guards can use those carriers to see every target.

\begin{proposition}\label{prop:minkowski-web}
    There is an absolute constant \(C>0\) such that every rational target in \(T_B\) lies on an admissible line of coefficient height at most \(C\sqrt{B}\), for every \(B \ge 1\).
\end{proposition}

\begin{proof}
    For a primitive representative \(\widetilde q=(a,b,c)\) of a rational target, the integer coefficient vectors of lines through \(q\) form the rank-two lattice
    \[
        \Lambda_q=\{u \in \Z^3 : u \cdot \widetilde{q} = 0\}.
    \]
    By \cref{lem:line-uniform}, its covolume in its real span is \(\norm{\widetilde q}_2 \le \sqrt{3}B\). Minkowski's convex-body theorem gives a nonzero \(u \in \Lambda_q\) with \(\norm{u}_\infty \le C \sqrt{B}\). With \(u=(u_0,u_1,u_2)\), the coefficients \(u_0,u_1\) cannot both vanish: since \(c>0\), the relation \(u \cdot \widetilde{q}=0\) would then give \(u_2=0\), contrary to \(u \ne 0\). Thus \(u_0X+u_1Y+u_2=0\) defines a line through \(q\). Taking its primitive part preserves the coefficient bound.
\end{proof}

The proposition says that the cost of finding a line through every target is \(O(\sqrt{B})\). The following theorem says that covering them with boundedly many degree-one guards costs \(\Theta(B)\).

\begin{theorem}[Height scale for degree-one guarding]\label{thm:linear-threshold}
    Let \(P\) be a rational polygon with \(n\) vertices. For every positive integer \(k\), if
    \[
        \gamma_{1,H(B)}(P;T_B(P),T_B(P))\le k
    \]
    for all sufficiently large \(B\), then a constant \(c_{P,k}>0\) satisfies \(H(B)\ge c_{P,k}B\) for all sufficiently large \(B\). Conversely, there is a constant \(C_P>0\) such that
    \[
        \gamma_{1,C_PB}(P;T_B(P),T_B(P)) \le \left\lfloor \frac{n}{3} \right\rfloor
    \]
    for all sufficiently large \(B\). Thus the coefficient-height scale for a number of degree-one rational guards bounded independently of \(B\) is \(\Theta_P(B)\).
\end{theorem}

The proof, in \cref{sec:lattice-proofs}, counts the targets visible from one guard uniformly in its position. Suppose first that \(H \le B\). A rational line of coefficient height \(h \le B\) contains \(O(B^2/h)\) targets, and the lines through a rational guard with heights in \((h/2,h]\) number \(O(h^2)\). Summing over the bands with upper endpoints \(H, H/2, H/4, \ldots\) therefore bounds one guard's coverage by \(O(B^2H)\). Since the polygon contains \(\Theta_P(B^3)\) targets, covering them with a bounded number of guards requires \(H\) to be at least a positive constant times \(B\). If \(H>B\), this lower bound already holds. The matching upper bound follows from Fisk's vertex guards (\cref{cor:fisk}) and the joining-line height bound (\cref{prop:linear}).

In the square the threshold is attained in the simplest possible way, by a single guard at a corner.

\begin{proposition}\label{prop:square-collapse}
    For every \(B \ge 1\),
    \[
        \gamma_{1,B}(\mathcal{S};T_B,T_B) = 1.
    \]
\end{proposition}

\begin{proof}
    Use the origin as guard. A target \(t=(a/c,b/c) \ne (0,0)\) has the joining line \(Z(bX-aY)\). After division by its content \(\gcd(a,b)\), its coefficient height is at most \(B\), and the segment remains in the square. The origin sees itself, so this gives a one-guard cover. Since \(T_B\) is nonempty, the guard number is exactly one.
\end{proof}

\section{Algebraic relation degree and visibility}\label{sec:relation}

Increasing the height budget helped us guard rational targets without increasing the degree. For algebraic targets this approach can fail. An algebraic point may lie on no rational curve of the permitted degree, however large we allow its coefficients to be. In that case, increasing the degree is necessary for sight from another point. We study the least degree of a nonzero rational polynomial vanishing at an algebraic point and establish a sharp bound in terms of the target-field degree.

For \(q=(a,b)\in\Qbar^2\), consider its evaluation ideal
\[
    \mathfrak m_q\coloneqq \ker \left(\Q[X,Y] \xrightarrow {\operatorname{ev}_q}\Q(a,b) \right).
\]
The rational obstruction index of \([1:a:b]\) is the least degree of a rational projective hypersurface containing it \cite[Definition~2.1(i)]{AmorosoViada2012}. Homogenizing and dehomogenizing equations gives its affine form
\begin{equation}\label{eq:relation-degree}
    \omega_\Q(q)\coloneqq \omega_\Q([1:a:b]) = \min_{0 \ne F \in \mathfrak m_q} \deg F.
\end{equation}
We use \emph{(algebraic) relation degree} as an affine name for this standard invariant. The minimum is over the entire ideal, not just a chosen set of generators. The coefficient field \(\Q\) matters: allowing all algebraic coefficients would give linear equations through every point.

\begin{proposition}\label{prop:relation-obstruction}
    If \(\omega_\Q(q)>d\), no carrier of degree at most \(d\) passes through \(q\), regardless of coefficient height. A target at \(q\) can therefore be covered only by a guard at \(q\).
\end{proposition}

\begin{proof}
    An admissible carrier would supply a nonzero element of \(\mathfrak m_q\) of degree at most \(d\), contradicting \cref{eq:relation-degree}. Thus \(q\) must guard itself. If \(q\) is not an allowed guard location, no cover exists.
\end{proof}

For a fixed number field \(K\), how large can the relation degree be as \(q\) varies in \(K^2\)? The field degree \(e=[K:\Q]\) imposes an  upper bound: there are \(M_d=\binom{d+2}{2}\) monomials of total degree at most \(d\), and their evaluations at \(q\) lie in the \(e\)-dimensional \(\Q\)-vector space \(K\). Indeed, evaluation is a linear map
\[
    \operatorname{ev}_q: \Q[X,Y]_{\le d} \longrightarrow K,\qquad \dim_\Q \Q[X,Y]_{\le d} = M_d.
\]
If \(M_d>e\), its kernel is nonzero, so every \(q \in K^2\) satisfies a relation of degree at most \(d\). The next theorem says that this dimension bound is sharp. It concerns existence of equations through points, not yet existence of sight arcs.

\begin{theorem}[Sharp field-degree threshold]\label{thm:field-threshold}
    Let \(K \subset \R\) be a number field of degree \(e\), and let \(d \ge 1\).
    \begin{enumerate}
        \item If \(e < M_d\), then \(\omega_\Q(q) \le d\) for every \(q \in K^2\).
        \item If \(e \ge M_d\), choose a \(\Q\)-basis of \(K\) and identify \(K^2\) with \(\Q^{2e}\). The points with \(\omega_\Q(q)>d\) are the rational points of a nonempty Zariski-open subset of \(\A^{2e}_\Q\). Under the chosen real embedding, every nonempty open subset of \(\R^2\) contains infinitely many such points.
    \end{enumerate}
\end{theorem}

\begin{proof}
    The first assertion is the dimension argument above. For the second, the expected-dimension theorem of Asgarli, Love, and Yip \cite[Theorem~1.2]{AsgarliLoveYip2025} supplies a point of \(\PP^2(K)\) at which no nonzero rational degree-\(d\) form vanishes when \(e \ge M_d\). In particular, the form \(X_0^d\) does not vanish there, so the point already lies in the standard chart \(X_0 \ne 0\). Dehomogenization gives an affine point \(q\) with \(\omega_\Q(q)>d\). A lower-degree affine relation homogenizes, after multiplication by a power of \(X_0\), to a degree-\(d\) form vanishing at the projective point.

    In the chosen basis of \(K\), the evaluation matrix has polynomial entries in the \(2e\) rational coordinates of \(q\). Injectivity is the nonvanishing of at least one maximal minor. The resulting Zariski-open set is nonempty by the point just found.

    The real linear map \(\R^{2e} \to \R^2\) induced by the chosen embedding is surjective. The inverse image of a nonempty open subset is open and meets the maximal-rank locus, whose algebraic complement has empty interior. Rational coefficient vectors are dense there. There are infinitely many such vectors, and their images are distinct because the chosen basis embeds \(K^2\) into \(\R^2\).
\end{proof}

\begin{remark}\label{rem:relation-not-sight}
    The dimension argument above gives a nonzero rational polynomial \(F\) vanishing at the target \(q\). Siegel's lemma can bound the coefficient height of a primitive integer equation \(F\). If \(F\) is reducible, choose a primitive integer irreducible factor \(G\) that still vanishes at \(q\). Gelfond's lemma bounds the coefficient height of \(G\) in terms of the degree and height of \(F\) \cite[Lemma~1.6.11]{BombieriGubler2006}. This bound may exceed the allowed budget \(H\) because factors need not have smaller coefficients than the original polynomial.

    Even \emph{if} \(G\) fits the budget, it need not be absolutely irreducible or nonsingular at \(q\), as in \cref{rem:conditions-example}. And even when these requirements hold, a sight arc must still join the guard to \(q\) inside \(P\). A polynomial relation through \(q\) is therefore necessary for sight from another point, but not sufficient.
\end{remark}

\begin{corollary}\label{cor:no-uniform}
    Let \(P\) be a rational polygon, \(d \ge 1\), and let \(K \subset \R\) have degree at least \(M_d\). For every positive integer \(N\), there is a \(B\) such that, for every \(H \in [1,\infty]\), 
    \[
        \gamma_{d,H}(P;T_{K,B}(P),T_{K,B}(P)) \ge N.
    \]
This holds \emph{even} when \(P\) is convex.
\end{corollary}

\begin{proof}
    By \cref{thm:field-threshold}, choose \(N\) distinct interior points with relation degree greater than \(d\) and set \(B\) as their maximum height. Each selected point belongs to \(T_{K,B}(P)\) and requires a guard at its own location by \cref{prop:relation-obstruction}. Every cover therefore contains these \(N\) points, regardless of \(H\). Thus the same \(B\) proves the lower bound at every coefficient-height budget.
\end{proof}

Increasing height alone therefore cannot restore sight when the permitted degree is below a target's relation degree. The next section shows that, once the degree may also grow, a rational interior guard can see every target of a finite algebraic set.

\section{Arithmetic interpolation and arc lifting}\label{sec:lifting}

If we allow both carrier bounds \(d\) and \(H\) to grow, then a rational carrier can contain any finite collection of algebraic points. Visibility requires more: the guard and target must be joined by a regular arc of that carrier contained in the polygon.

The first proposition is a geometric consequence of the Shape Lemma for a finite reduced set. Cox, Little, and O'Shea formulate the relevant coordinate description as an exercise \cite[Chapter~2, \S4, Exercise~16]{CoxLittleOShea2005}. We give the interpolation argument in the form needed here.

\begin{proposition}\label{prop:interpolation}
    Every finite subset of \(\Qbar^2\) is contained in a nonsingular, absolutely irreducible affine plane curve over \(\Q\).
\end{proposition}

\begin{proof}
    For an empty initial set any rational line suffices. Let \(\Omega\) be the finite Galois-stable set consisting of the given points and their conjugates. Choose \(c \in \Q\) so that \(U=X+cY\) separates its points, meaning that it takes distinct values at distinct points. Use the new coordinates \((U,V)=(X+cY,Y)\), keeping the second coordinate unchanged so that we can interpolate a graph \(V=r(U)\). For each pair of distinct points at most one value of \(c\) is forbidden. By Lagrange interpolation we obtain a unique polynomial \(r \in \Qbar[U]\), of degree less than \(\card{\Omega}\), with
    \[
        r(U(z)) = V(z),\qquad z \in \Omega.
    \]
    Uniqueness and Galois equivariance of the data imply \(r \in \Q[U]\). The curve \(Z(V-r(U))\) contains \(\Omega\). It is absolutely irreducible and nonsingular. In the coordinates \(U,V\) its equation is monic and linear in \(V\), with \(V\)-derivative one. Returning to the original coordinates, clear denominators and remove the content to obtain a primitive integer equation.
\end{proof}

The real curve just constructed is the graph of a polynomial in the coordinate \(U\), so its real points are connected. Its regular trace inside \(P\) need not be connected. Interpolation therefore settles incidence but does not ensure that sight arcs joining the points stay inside \(P\).

\begin{theorem}[Arithmetic arc lifting]\label{thm:lifting}
    Let \(P\) be a rational polygon, let \(g \in (\operatorname{int}P)(\Q)\), and let \(t \in P(\Qbar)\). There is a primitive absolutely irreducible \(F \in \Z[X,Y]\) such that \(g\) and \(t\) lie on the same connected component of \(\mathrm{R}_F(P)\). The projective closure of \(\mathrm{C}_F\) may be chosen smooth.
\end{theorem}

\paragraph{How the construction works.}
An arbitrary approximation could move an algebraic endpoint off the curve. We therefore begin with a rational polynomial \(F_0\) whose values and first derivatives are fixed at the endpoints and their Galois conjugates. The resulting local branches supply the ends of a smooth reference arc. At a boundary endpoint its tangent direction points strictly into the polygon.

Next we construct a defining function \(f\) for a compact regular real curve containing this arc, arranging that \(f=F_0\) near the fixed points. Form \(A_0\) as a sum of squares of rational polynomials vanishing at the fixed points, using the evaluation kernel described in \cref{sec:effective-preparation}. The polynomial \(A_0\) vanishes to second order there, so a correction of the form \(A_0\beta\) can improve the approximation without changing the first-order data. Control in \(\mathscr{C}^1\), rather than merely uniform control of function values, preserves the regular arc and its inward direction at the boundary.

Finally, a small rational perturbation through the fixed points makes the projective curve smooth. This Bertini step need not fix the derivatives exactly: it keeps them sufficiently close for the arc to persist. The common proof of this theorem and its effective refinement is in \cref{sec:effective-proof}. The figures follow the same polygon and endpoints through the construction; see \cref{fig:lifting-assembly} for the reference curve.

The earlier \cref{lem:convex-connector} solves a special case explicitly: convexity and a rational line through the target give a carrier of degree at most the target's degree. Here the polygon may be nonconvex and the target need not lie on a rational line.

The qualitative conclusion also holds for an algebraic interior guard. For distinct endpoints, let \(S\) be the union of their Galois orbits, with repeated points included only once, and put \(N=\card{S}\). Let \(\mathcal{V}\) be the rational vector space of polynomials of degree at most \(2N-1\) vanishing throughout \(S\). After extending scalars to \(\R\), the interpolation construction in the proof of \cref{lem:equivariant-jets} allows the gradients at the real points of \(S\) to be assigned independently.

Next choose an oriented polygonal path from \(t\) to \(g\) with no self-intersections, lying in \(\operatorname{int}P\) except possibly at \(t\), and avoiding the other real points of \(S\). If \(t \in \partial P\), require its initial tangent to point strictly into \(P\). Choose the endpoint covectors as the clockwise normals to the oriented endpoint tangents, and choose nonzero covectors at the other real points of \(S\). 

Then express an interpolant with these gradients in a rational basis of \(\mathcal{V}\), and approximate its basis coefficients by rational numbers. The resulting rational polynomial \(F_0\) still vanishes throughout \(S\), and its gradients remain nonzero and sufficiently close to preserve the endpoint orientations and the strict inward condition. Its first jets are automatically Galois-equivariant. 

Now adjust the path's short endpoint segments to follow the actual tangent directions of \(F_0=0\) at \(t\) and \(g\), choosing their orientations close to the original ones and reconnecting them locally to the remaining path. With the coefficient approximation sufficiently close, these adjustments preserve containment in \(P\), avoidance of the other real points of \(S\), and the absence of self-intersections. If \(t \in \partial P\), its initial tangent remains strictly inward. This construction remains valid when \(g\) and \(t\) belong to the same Galois orbit. The reference-function construction, polynomial approximation, and smooth projective completion then proceed with this interpolation set.

\paragraph{Relation to the literature.}
The ideal-generator approximation strategy belongs to the Tognoli and Akbulut--King tradition discussed by Ghiloni--Savi \cite[\S3 and \S3.1]{GhiloniSavi2025}. Their rational-coefficient approximation \cite[Lemma~3.3(ii)--(iii)]{GhiloniSavi2025} expresses the function as a sum of polynomial ideal generators with smooth coefficients, then approximates those coefficients by rational polynomials. We use the same strategy in a finite-point setting, but arrange exact agreement with an interpolant near the points. Since \(f=F_0\) near the fixed real points, the quotient \((f-F_0)/A_0\) extends by zero to a \(\mathscr{C}^2\) function. The correction lies in the squared ideal, thus preserving first jets. The regular arc's containment at a polygon boundary and the rational smooth projective completion are further requirements, not conclusions of that lemma.

Fernando studies polynomial and Nash paths through specified control points \cite[Main Theorems~1.8--1.9]{Fernando2025}. Our additional requirements concern a rational implicit carrier, regularity at the endpoints and along the chosen arc, and a smooth projective closure. In particular, absolute irreducibility of the image of a rational polynomial parameterization does not itself give these regularity and projective conditions.

Effective Bertini--Northcott results control arithmetic invariants without requiring a connecting real arc to remain inside the polygon \cite[Theorem~1.3]{PazukiWidmer2021}. We now bound the approximation and the final perturbation while retaining that requirement.

\subsection{Effective bounds}
The degree and coefficient-height bounds obtained from this construction are computable, but can be very large. A narrow passage or an endpoint close to a corner can require much finer approximation than the number of sides alone suggests. Below we compute a reference arc and its margins inside the polygon, then turn those margins into degree and height bounds.

A rational polynomial \(F_0=P_0/D_0\), with \(\norm{P_0}_\infty =C_0\), fixes the endpoint jets. A polynomial \(A_0\) with vanishing first jets has sum of absolute coefficient values at most \(A_*\) and \(\mathscr{C}^1\) norm at most \(A_C\) on \([-R,R]^2\). The \(\mathscr{C}^2\) factor \(u=(f-F_0)/A_0\) has \(\mathscr{C}^2\) norm at most \(M\). The tolerance \(\epsilon_0\) is small enough to keep the approximating arc regular and inside \(P\), including at a boundary target. Finally \(L\) controls the sizes of the evaluation matrices used to choose polynomial bases. \Cref{sec:effective-proof} computes all these quantities from the input; the formulas below turn them into carrier budgets.

\begin{theorem}[Effective arithmetic arc lifting]\label{thm:effective-lifting}
    Let \(P\) be a rational polygon, let \(g \in (\operatorname{int}P)(\Q)\), and let \(t \in P(\Qbar)\). From exact algebraic descriptions of these data one can compute integers \(m,H \ge 1\) and a primitive polynomial \(F \in \Z[X,Y]\) of degree at most \(m\) and coefficient height at most \(H\), with smooth absolutely irreducible projective closure, such that \(g\) and \(t\) lie on the same component of \(\mathrm{R}_F(P)\).

    For distinct endpoints, the construction gives the following explicit bounds. Put \(N=[\Q(t_x,t_y):\Q]+1\). The preparation in \cref{sec:effective-proof} computes positive integers \(R,M,A_*,D_0,C_0,L\) and positive rational numbers \(\epsilon_0,A_C\). Define
    \[
        \varepsilon = \frac{\epsilon_0}{4A_C},\qquad
        n = \max\left\{ 1,\left\lceil \left( \frac{8RM}{\varepsilon} \right)^2 \right\rceil \right\},
        \qquad Q_u = \max\left\{ 1,\left\lceil \frac{n+1}{\varepsilon} \right\rceil \right\}.
    \]
    Let
    \[
        d_u = Q_u(2R)^{2n},\qquad H_u = (Q_uM+1)[4(1+R)^2]^n,
    \]
    and put
    \[
        m = 2N+2n,\qquad D_G = D_0d_u,\qquad C_G = d_uC_0+D_0A_*H_u.
    \]
    With \(J_m = \binom{m+2}{2}\), \(k = J_m-N\), \(b_m = N!L^{mN}\), and \(\delta = 3(m-1)^2\), one may take
    \[
        Q_{\Gamma} = 1+\left\lceil \frac{4\delta k m J_mb_m R^m}{\epsilon_0} \right\rceil, \qquad H = Q_{\Gamma}C_G+D_G\delta k J_mb_m.
    \]
    All the preparation data are obtained by terminating calculations from \(P,g,t\). If \(g=t\), the primitive equation of the vertical rational line through \(g\) suffices.
\end{theorem}

The proof in \cref{sec:effective-proof} combines an explicit rational polynomial approximation with a small rational perturbation that results in a smooth projective curve. Using a tensor-product Bernstein polynomial with rounded sample values, we construct a rational approximation with explicit degree and height bounds. Our construction directly bounds the errors in function values and first derivatives throughout the square \([-R,R]^2\). Let \(S\) consist of \(g\) and the Galois conjugates of \(t\). The discriminant is not identically zero on the linear system of homogeneous degree-\(m\) polynomials vanishing at \(S\), and its degree \(3(m-1)^2\) gives an explicit bound on the number of values needed along each coordinate of a finite parameter grid.

Note that these are only sufficient, not optimal, bounds. For a finite target set, the construction gives common computable budgets by taking the largest degree bound and the largest height bound. Given how coarse this approach is, the much smaller budgets achieved by the conics in \cref{prop:l-rational-conic} illustrate why direct constructions remain useful. Sharp bounds in terms of the polygon and target data remain open; see \cref{prob:effective-lifting}.

\begin{corollary}[Eventual one-guard stabilization]\label{cor:stabilization}
    Let \(P\) be a rational polygon, let \(T\) be a nonempty finite subset of \(P(\Qbar)\), and let \(g \in (\operatorname{int}P)(\Q)\). There are computable finite \(d_0,H_0\) such that
    \[
        \gamma_{d,H}(P;T,T\cup\{g\})=1,\qquad d \ge d_0,\ H \in [H_0,\infty].
    \]
    If \(g \in T\), the same conclusion holds with candidate set \(T\). In particular this applies to \(T_{K,B}(P)\) or \(T_{D,B}(P)\) as soon as \(B \ge \height(g)\).
\end{corollary}

\begin{proof}
    For each \(t \in T \setminus\{g\}\), apply \cref{thm:effective-lifting} with target \(t\), obtaining an equation \(F_t\). Take \(d_0,H_0\) to be the respective maxima of one and the degrees and coefficient heights of these finitely many equations. Then \(g\) sees every target; it sees itself by convention. Since \(T\) is nonempty, the guard number is therefore one.
\end{proof}

\begin{corollary}\label{cor:edge-eventual}
    Let \(P\) be a rational polygon and let \(T\) be a nonempty finite subset of \(P(\Qbar)\). There are computable \(d_*,H_*\) such that any edge of \(P\) guards \(T\) whenever \(d \ge d_*\) and \(H \in [H_*,\infty]\). In particular, \(\eta_{d,H}(P;T)=1\) throughout this range. 
\end{corollary}

To obtain this edge version, attach a small rational triangle outside \(P\) along a short segment of the chosen edge, and put the rational guard in its interior. Apply \cref{thm:effective-lifting} in the enlarged polygon, then retain the part of each arc after its last entry into \(P\). That entry lies on the selected edge. Taking maxima over the finitely many edges gives common budgets. The full argument is in \cref{sec:edge-deferred}.

We now compare this eventual one-guard coverage with straight-line guarding.

\begin{proposition}\label{prop:empty-kernel}
    If the classical kernel \(\operatorname{Ker}(P)\) is empty, then for all sufficiently large \(B\) no single straight-line point guard anywhere in \(P\) covers \(T_{\Q,B}(P)\). Nevertheless, a rational interior guard sees each such finite target set arithmetically at some finite budget.
\end{proposition}

\begin{proof}
    If arbitrarily large target slices admitted a single straight-line guard, a limit of these guards would see all of \(P\). To show this, choose \(B_j\) tending to infinity and corresponding guards \(g_j\). Compactness gives a subsequence converging to \(g \in P\). Fix a rational target \(t\). It belongs to every sufficiently large slice, so \([g_j,t] \subseteq P\) eventually. For each \(\lambda \in [0,1]\), closedness keeps the limit of \((1-\lambda)g_j+\lambda t\) in \(P\). Thus \([g,t] \subseteq P\). This same guard sees every rational target. Now approximate a real target by rational targets in \(P\) and apply closedness once more: \(g\) sees every point of \(P\), contradicting its empty kernel. The arithmetic assertion follows from \cref{cor:stabilization}, with one rational interior guard and budgets that may depend on the slice.
\end{proof}

These budgets depend on the target set. A finite set has a common sufficient budget, whereas a rational interior guard can see \emph{all} algebraic targets only with bounds allowed to depend on the target. At finite \(d,H\), the finitely many carriers cannot even contain every rational point in \(P\). At fixed degree and unrestricted height, algebraic targets of larger relation degree still cannot be reached. Removing both bounds, on the other hand, makes every rational interior point see every algebraic target by \cref{thm:lifting}. These facts motivate the questions about fixed-degree kernels in \cref{prob:field-kernel}; simply removing budgets does not retain the geometric information carried by the classical kernel.

\section{Sight at the next degree}\label{sec:jump}

Arc lifting gives sufficient budgets, but does not determine how quickly visibility might change with increases in degree. Can the guard number change all at once when the degree increases by just one? To explore this possibility, let us return to the unit square \(\mathcal{S}\) and consider the curve
\[
    \mathrm{C} \coloneqq Z(Y-X^{d+1}).
\]
Its points over any field \(K/\Q\) are \((a,a^{d+1})\), with \(a \in K\). Its real trace in the square is a connected graph, and hence supplies a sight arc between any two of its points. Polynomial-graph visibility, including monomial graphs, also occurs in \cite{GoinsHarrisKubikMbirika2018}; there the question is which lattice points are first reached, not which regular arcs have bounded rational defining equations.

We will choose algebraic targets on \(\mathrm{C}\) that satisfy no nonzero rational polynomial equation of degree at most \(d\). At those degrees, each target must therefore guard itself. At degree \(d+1\), however, \(\mathrm{C}\) becomes an admissible carrier of coefficient height one, and any target sees all the others along it. To obtain arbitrarily many targets without increasing their projective heights, we construct their coordinates from a cyclotomic unit and take all conjugates. The relation-degree criterion of \cref{sec:relation} explains why the smaller degree budget prevents sight; the defining polynomial of \(\mathrm{C}\) shows exactly when sight becomes possible. Here those bounds are exact, rather than the coarse, but sufficient, bounds of the general construction.

\begin{theorem}[Degree jump and guard number collapse]\label{thm:main-jump}
    For every \(d \ge 1\) and \(M \ge 1\), there exists a finite Galois-stable target set
    \[
        T \subset(\operatorname{int}\mathcal{S})(\Qbar)
    \]
    such that
    \begin{enumerate}[label=(\roman*)]
        \item \(\card{T}>M\) and every \(t \in T\) has projective height exactly \(4^{d+1}\);
        \item
        \[
         \gamma_{d,\infty}(\mathcal{S};T,T) = \card{T};
        \]
        \item
        \[
         \gamma_{d+1,1}(\mathcal{S};T,T) = 1.
        \]
    \end{enumerate}
    In other words, the arithmetic kernel changes from \(\operatorname{Ker}_{d}(\mathcal{S};T,T)=\varnothing\) to \(\operatorname{Ker}_{d+1,1}(\mathcal{S};T,T)=T\). An increase of one in carrier degree therefore changes an arbitrarily large guard number to one in the same convex integer polygon.
\end{theorem}

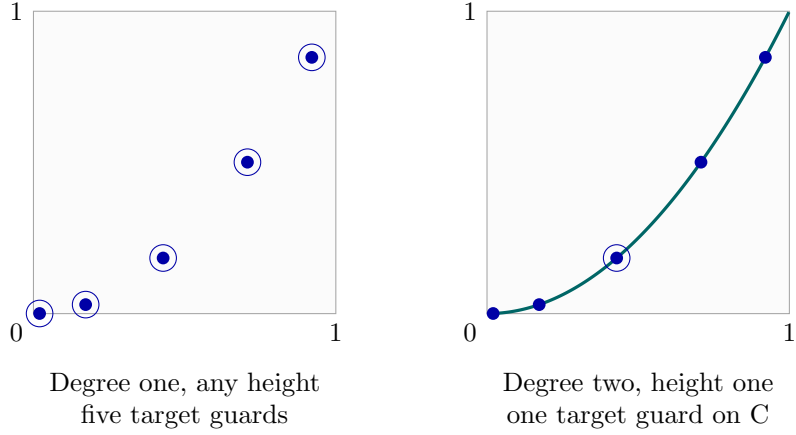
\begin{figure}[tbp]
    \centering
    \begin{tikzpicture}[font=\small]
        \foreach \shift/\deg in {0/1,6/2} {
        \begin{scope}[shift={(\shift,0)},x=4cm,y=4cm]
            \fill[gray!3] (0,0) rectangle(1,1);
            \draw[gray!70] (0,0) rectangle(1,1);
            \ifnum\deg=2
            \draw[teal!80!black,very thick,domain=0:1,samples=100] plot(\x,{\x*\x});
            \fi
            \foreach \k in {1,...,5} {
            \pgfmathsetmacro{\a}{(1+cos(360*\k/11))/2}
            \fill[blue!65!black] (\a,{\a*\a}) circle(2.4pt);
            \ifnum\deg=1
            \draw[blue!65!black] (\a,{\a*\a}) circle(5pt);
            \fi
            \ifnum\deg=2
            \ifnum\k=3 \draw[blue!65!black] (\a,{\a*\a}) circle(5pt); \fi
            \fi
            }
            \node[below left] at (0,0) {\(0\)};
            \node[below] at (1,0) {\(1\)};
            \node[left] at (0,1) {\(1\)};
            \ifnum\deg=1
            \node[align=center,text width=4.6cm,anchor=north] at (.5,-.16)
            {Degree one, any height\\five target guards};
            \else
            \node[align=center,text width=4.6cm,anchor=north] at (.5,-.16)
            {Degree two, height one\\one target guard on \(\mathrm{C}\)};
            \fi
        \end{scope}
        }
    \end{tikzpicture}

    \caption{The same five algebraic targets at two degree budgets. At degree one, each must guard itself, even with unrestricted coefficient height. At degree two, any one target sees all five along \(\mathrm{C}=Z(Y-X^2)\), whose equation has coefficient height one. Here \(\ell=11\), so the targets have degree five and projective height sixteen.}
    \label{fig:degree-jump}
\end{figure}

\begin{proof}
    \emph{Choose the coordinates.} Choose an odd prime \(\ell\) such that
    \[
        e \coloneqq \frac{\ell-1}{2} > \max\{d(d+1),M\}.
    \]
    Let \(\zeta_\ell\) be a primitive \(\ell\)th root of unity and set
    \[
        \theta \coloneqq 2+\zeta_\ell+\zeta_\ell^{-1}.
    \]
    The maximal real subfield of \(\Q(\zeta_\ell)\) has degree \(e\), and \(\theta\) is a primitive element.  Its \(e\) conjugates are
    \[
        \theta_k = 2+2\cos(2\pi k/\ell),\qquad 1 \le k \le e,
    \]
    and are distinct elements of \((0,4)\). Now set
    \[
        \alpha = \theta/4.
    \]
    Its conjugates \(\alpha_\sigma\) are distinct elements of \((0,1)\); \cref{fig:jump-conjugates} illustrates this construction.

    \begin{figure}[tbp]
        \centering
        \begin{tikzpicture}[scale=1.15,every node/.style={font=\small}]
            \draw[black!55] (0,0) circle (1.25);
            \draw[->,black!55] (-1.65,0)--(1.7,0);
            \draw[->,black!55] (0,-1.55)--(0,1.65);
            \foreach \k in {1,...,10} {
            \fill[black] ({1.25*cos(360*\k/11)},{1.25*sin(360*\k/11)}) circle (1.4pt);
            }
            \foreach \k in {1,...,5} {
            \draw[draw=linkcyan,dashed] ({1.25*cos(360*\k/11)},{1.25*sin(360*\k/11)})
              --({1.25*cos(360*\k/11)},{-1.25*sin(360*\k/11)});
            \fill[linkcyan] ({1.25*cos(360*\k/11)},0) circle (2pt);
            }
            \node[above] at (0,1.72) {\(\zeta_{11}^{\,k}\)};
            \node[below] at (0,-1.67) {complex conjugate pairs};
            \draw[->,thick] (2.03,.35)--(3.15,.35);
            \node[align=center,font=\scriptsize] at (2.58,.73) {\(z\longmapsto(1+\operatorname{Re}z)/2\)};
            \draw[->,black!55] (3.3,0)--(6.6,0);
            \draw (3.45,.06)--(3.45,-.06) node[below=4pt] {\(0\)};
            \draw (6.45,.06)--(6.45,-.06) node[below=4pt] {\(1\)};
            \foreach \k in {1,...,5} {
            \pgfmathsetmacro{\a}{(1+cos(360*\k/11))/2}
            \fill[linkcyan] ({3.45+3*\a},0) circle (2pt);
            }
            \node[align=center] at (4.9,-1.05)
            {five real conjugates\\in \((0,1)\)};
        \end{tikzpicture}

    \caption{For \(\ell=11\), conjugate pairs of nontrivial roots of unity give the five real coordinates \(\alpha_k=(1+\cos(2\pi k/11))/2\), \(1\le k\le5\), in \((0,1)\), each of degree five.}
    \label{fig:jump-conjugates}
    \end{figure}
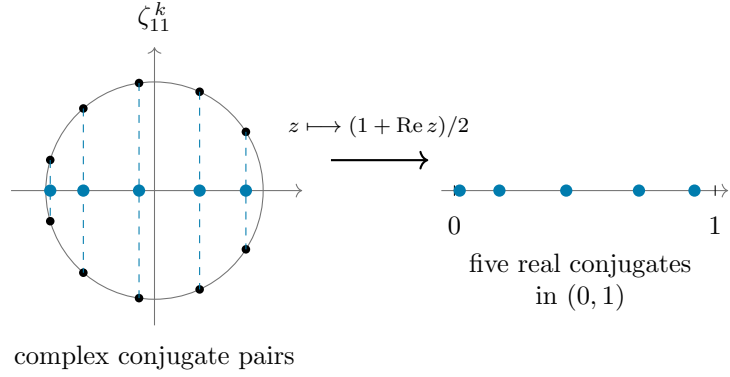

    For every embedding \(\sigma : \Q(\theta) \hookrightarrow \R\), define
    \[
        t_\sigma = \left( \alpha_\sigma,\alpha_\sigma^{d+1} \right)
    \]
    and let \(T\) be the set of these \(e\) points. It is a single Galois orbit and lies in the interior of \(\mathcal{S}\). By the strict inequality on \(e\), we have \(\card{T}=e>M\) and \(d(d+1)<e\); the latter is needed in the degree comparison below.

    \emph{Determine the relation degree of each target.} Suppose \(0 \ne F \in \Q[X,Y]\) has total degree at most \(d\) and vanishes at one \(t_\sigma\).  Substitute
    \[
        X = U,\qquad Y = U^{d+1}
    \]
    to obtain
    \[
        G(U) = F \left(U,U^{d+1}\right) \in \Q[U].
    \]
    Each monomial \(X^iY^j\) with \(i+j \le d\) becomes \(U^{i+(d+1)j}\), so
    \[
        \deg G \le d(d+1) < e.
    \]
    Since \(G(\alpha_\sigma)=0\) and \(\alpha_\sigma\) has degree \(e\), it follows that \(G\) is the zero polynomial.

    Thus \(F\) vanishes identically on \(\mathrm{C}\), not merely at the chosen target. Algebraically, the kernel of the substitution homomorphism
    \[
        \Q[X,Y] \longrightarrow \Q[U],\qquad X \longmapsto U,\quad Y \longmapsto U^{d+1},
    \]
    is the principal ideal generated by \(Y-X^{d+1}\). Indeed, division by this monic polynomial in \(Y\) leaves the remainder \(F(X,X^{d+1})\), which has just been shown to be zero. Thus \(F\) would be a nonzero multiple of a polynomial of total degree \(d+1\), contrary to \(\deg F \le d\).  Since the defining equation of \(\mathrm{C}\) has degree \(d+1\), this proves \(\omega_\Q(t_\sigma)=d+1\). By \cref{prop:relation-obstruction}, every target requires its own guard, and reflexivity gives
    \[
        \gamma_{d,\infty}(\mathcal{S};T,T) = e > M.
    \]

    \emph{One more degree makes the same curve available to every target.} At degree \(d+1\), consider
    \[
        F_{d+1}(X,Y) = Y-X^{d+1}.
    \]
    It is primitive, absolutely irreducible because it is monic and linear in \(Y\), and everywhere regular because
    \[
        \frac{\partial F_{d+1}}{\partial Y} = 1.
    \]
    Its trace in \(\mathcal{S}\) is the connected graph
    \[
        Y = X^{d+1},\qquad 0 \le X \le 1,
    \]
   which contains every target. Any one target therefore guards all of \(T\) at coefficient height \(1\). Since \(T\) is nonempty, \(\gamma_{d+1,1}(\mathcal{S};T,T)=1\).

    \emph{Keep height independent of orbit size.} The numerator of \(\alpha\) was chosen to be a unit. Since \(\theta=(1+\zeta_\ell)(1+\zeta_\ell^{-1})\) and \(\Phi_\ell(-1)=1\), the full cyclotomic norm of \(\theta\) is \(1\). The norm from the maximal real subfield is therefore \(1\) as well: its square is the full norm and all conjugates of \(\theta\) are positive. Thus \(\theta\) is an algebraic unit.  At every nonarchimedean place of \(\Q(\theta)\) one therefore has \(\abs{\theta}_v=1\), while every archimedean conjugate of \(\alpha=\theta/4\) lies in \((0,1)\). By the product formula \(\height(\alpha)=4\). Moreover, place by place,
    \[
        \max\{1,\abs{\alpha}_v,\abs{\alpha}_v^{d+1}\} = \max\{1,\abs{\alpha}_v\}^{d+1}.
    \]
    Consequently,
    \[
        \height([1:\alpha:\alpha^{d+1}]) = \height(\alpha)^{d+1}=4^{d+1},
    \]
    independently of \(\ell\), \(e\), and \(M\).

    By Galois invariance of the Weil height, every target in \(T\) therefore has height \(4^{d+1}\). Together with \(\card{T} = e > M\), this proves (i) and completes the proof.
\end{proof}

\begin{remark}\label{rem:jump-comparison}
The first coordinate has such large algebraic degree that substituting \(Y=X^{d+1}\) turns any putative lower-degree relation at the target into an identity on \(\mathrm{C}\). But a nonzero equation of this curve has degree at least \(d+1\). Galois orbits that impose independent conditions on low-degree forms are known~\cite{AsgarliGhiocaReichstein2024,AsgarliLoveYip2025}. Neither those results nor \cref{thm:field-threshold} is used in the proof above; they serve only as comparison. The proof uses the degree and height of a real cyclotomic unit and the calculation above of the kernel of the substitution map.
\end{remark}

We have now seen both sources of difficulty. A polygon can prevent a regular arc from remaining inside it; the alternating combs make this limitation persist over real coefficients at bounded degree. Arithmetic can instead exclude a suitable equation even in a square. Increasing height removes the scarcity caused by a finite carrier list, but not a target's relation-degree barrier. Arc lifting overcomes both barriers for each finite algebraic target set, and the present example shows that a single increase in degree can already be enough, with no increase in target height as the set grows.

\section{Questions suggested by the examples}\label{sec:questions}

The examples above leave three main questions: which visibility graphs occur, when one guard sees a whole arithmetic gallery, and how small the sight budgets can be. We keep these questions separate from possible changes to the sight rule. Of these, \cref{prob:height-growth,prob:effective-lifting}, on the height scale at higher degree and on sharper effective bounds, respectively, are direct extensions of results developed above.

\subsection{Visibility graphs}
\begin{problem}[Realizing arithmetic visibility graphs]\label{prob:target-recognition}
    For \(d \ge 1\), \(H \in [1,\infty]\), and a real number field \(K\), characterize the graphs \(\operatorname{Vis}_{d,H}(P,T)\) arising from finite target sets \(T \subset P(K)\) in rational simple polygons \(P\). How does fixing the polygon or requiring the targets to be its vertices change the question?
\end{problem}

Taking \(T=V(P)\), \(d=1\), and \(H=\infty\) gives classical vertex visibility recognition restricted to rational realizations \cite{TOPPVisibilityRecognition}. One may additionally prescribe the cyclic vertex order. A finite coefficient budget can even exclude a polygon side's line. For general finite targets, the alternating combs realize powers of paths on interior target sets. When \(P\) is fixed, bounding target height makes realization a finite search over subsets of \(T_{K,B}(P)\); finding the least height of a realization remains a further issue.

\subsection{Whole-gallery kernels}
At fixed degree and coefficient height, infinitely many rational targets lie on no admissible carrier. Removing both bounds lets every rational interior guard see every algebraic target. Removing only the coefficient-height bound need not restore one-guard coverage: with rational carriers, the degree bound can force the whole-gallery kernel to be empty even in a convex polygon.

Suppose \(d<[K:\Q]\). By rationally translating a primitive element of \(K\), choose \(\alpha\) in the horizontal interval of a small rectangle inside \(P\). Choose distinct rational \(y_0,y_1\) in its vertical interval. A rational polynomial \(F\) of degree at most \(d\) vanishing at \((\alpha,y_i)\) must satisfy \(F(X,y_i) \equiv 0\), since \(\alpha\) has degree greater than \(d\). Absolute irreducibility therefore forces its carrier to be the horizontal line \(Y=y_i\). These lines are disjoint, so no guard sees both targets.

Instead consider \(\operatorname{Ker}_d(P,K)\) from \cref{def:arith-kernel}: rational guards see all \(K\)-points along \(K\)-coefficient carriers of degree at most \(d\), without a height bound. At degree one this is \(\operatorname{Ker}(P) \cap \Q^2\); \cref{ex:field-kernel-l} shows its enlargement at degree two.

\begin{problem}[Fixed-degree kernels over a number field]\label{prob:field-kernel}
    Which rational polygons have nonempty \(\operatorname{Ker}_d(P,K)\)? For fixed \(K\), can the number of vertices or reflex vertices bound a sufficient degree? How does this kernel compare with the set of rational guards that see every real point of \(P\) along real-coefficient carriers of degree at most \(d\)?
\end{problem}

For rational rather than \(K\)-coefficient conics, does the guard in \cref{prop:l-rational-conic} see every quadratic target in the upper arm of \(\mathcal{L}\), not just the displayed family?

\subsection{Degree and height}
For rational targets, degree one requires and admits coefficient height of order \(B\) (\cref{thm:linear-threshold,prop:square-collapse}). Quadratic targets and degree-one edge guards instead require order \(B^2\), even for a cover by all edges (\cref{prop:edge-quadratic-height}). These different scales make the degree-two point-guard question worth isolating. A further question, on edge guarding at degree two, is raised at the end of \cref{sec:edge-guards}.

\begin{problem}[Height needed for boundedly many guards]\label{prob:height-growth}
    For \(d \ge 2\), find the least growth of \(H(B)\) keeping \(\gamma_{d,H(B)}(\mathcal{S};T_{\Q,B}(\mathcal{S}),T_{\Q,B}(\mathcal{S}))\) bounded. How do other polygons and target fields change the answer?
\end{problem}

Over fields with \([K:\Q] \ge \binom{d+2}{2}\), relation degree prevents boundedly many guards at every height (\cref{thm:field-threshold}).

\begin{problem}[Sharper degree and height bounds]\label{prob:effective-lifting}
    Sharpen \cref{thm:effective-lifting}. Can bounds use polygon vertex data, endpoint degrees and heights, and intrinsic arc geometry rather than the constructed interpolation and cutoff constants? Include boundary distance away from an endpoint and an inward tangent angle at it. What degree--height tradeoffs are possible?
\end{problem}

The general theorem gives coarse formulas; the conics in \cref{sec:l-gallery} and the construction in \cref{cor:star-boundary-height} inside a triangle contained in \(P\) show how direct geometry can produce smaller degree bounds. For the alternating comb, the least height recovering every permitted adjacency may likewise be much smaller than \(H_{m,d}\). Recovering the whole graph and attaining its minimum guard number need not happen at the same height.

\subsection{Changing the sight rule}
\cref{ex:nodal-incidence} shows that singularities are not a cosmetic choice. Allowing passage through the node connects points lying on different regular components. A useful comparison should distinguish branch switching from following a single branch through a singularity.

\begin{problem}[Singularities and admissible sight arcs]\label{prob:singular-traces}
    For the same arithmetic gallery, compare guard numbers for absolutely irreducible carriers under four rules:
    \begin{enumerate}[nosep]
        \item the entire affine carrier is nonsingular;
        \item the sight arc avoids singularities, as here;
        \item the sight arc lifts continuously to the real normalization;
        \item connected full traces permit branch switching at singularities.
    \end{enumerate}
    How large can the differences be, and what budget increase replaces singular sight by regular sight?
\end{problem}

The nodal example distinguishes the second and fourth rules at the same degree and coefficient height. The third requires compatible real preimages and does not permit a jump between the two preimages of a node.

\paragraph{Further directions.}
Pairwise visibility does not record whether several targets lie on one regular component of one carrier. Keeping those component subsets as hyperedges suggests a finer realization question than the graph problem above.

Changing the coefficient field, as in the L-shaped example, brings other absolute values into the height. With a specified number field and bounds on carrier degree and coefficient height, the carrier list is finite. The same holds if only the degree of the coefficient field is bounded. A primitive representative with unit-ideal content need not exist over its ring of integers, so projective coefficient vectors are the natural replacement for primitive integer equations. 

Finite-place conditions can matter even before changing fields. Modulo three, \(Z(Y^2-X^2-3)\) splits into lines, whereas \(Z(Y^2-X^2-3X-1)\) is a smooth conic; both equations have degree two and height three. Joining endpoints by carriers with specified congruences or reduction behavior would add another arithmetic constraint. Congruences require a chosen normalization or a projective formulation.

Projective regions provide a different extension. The two affine branches of a hyperbola belong to one smooth connected real projective conic. To ask which visibility distinctions survive a rational projective transformation, one must specify the region, whether it includes points at infinity, and the coefficient-height convention.

\clearpage
\appendix
\crefalias{section}{appendix}
\crefalias{subsection}{appendix}
\part*{Appendices}
\addcontentsline{toc}{part}{Appendices}
The appendices provide the finite component tests, deferred geometric proofs, uniform lattice counts, and the common construction for qualitative and effective arc lifting.

\section{Details of the finite tests}\label{sec:finite-deferred}

\subsection{The nodal carrier exclusion}\label{sec:nodal-exclusions}
The parameterization \((U^2-1,U(U^2-1))\) gives the two targets of \cref{ex:nodal-incidence} at \(U=\pm5/4\). On the regular cubic its inverse is \(U=Y/X\); the excluded node has parameters \(\pm1\). The targets therefore lie on different regular components. Their two arms meet at the node inside the square. Absolute irreducibility follows because \(X^2(X+1)\) is not a square in \(\C(X)\).

Write \(A=\sum_{i+j\le3}c_{ij}X^iY^j\), with \(c_{ij} \in \{-1,0,1\}\), and put \(z=U^2-1\). In the monic basis
\[
    U^3z^3,\ U^2z^3,\ Uz^3,\ z^3,\ Uz^2,\ z^2,\ Uz,\ z,\ 1,
\]
of distinct degrees \(9,8,7,6,5,4,3,2,0\), the coefficients of \(A(z,Uz)\) are, respectively,
\[
    c_{03},\ c_{12},\ c_{21},\ c_{30}+c_{02},\ c_{11},\
    c_{20}+c_{02},\ c_{01},\ c_{10},\ c_{00}.
\]
A nonzero substituted polynomial has leading coefficient of absolute value at most two, so the rational-root theorem excludes a root with denominator four. Vanishing at either target therefore makes the substitution identically zero. Comparing coefficients gives \(A=c_{02}F\), so a nonzero admissible \(A\) is \(F\) or \(-F\). Thus the only admissible carrier through either target is the nodal cubic itself, whose regular components separate the two targets. The targets therefore do not see one another.

\subsection{The bounded-height target list}\label{sec:finite-implementation}
Write \(e=[K:\Q]\) and choose an integer \(A \ge 2^eB^e\). A target coordinate \(\alpha\) has degree \(r \le e\) and height at most \(B\). The elementary symmetric expressions in its conjugates bound its primitive minimal polynomial by
\[
    \norm{m_\alpha}_\infty \le 2^r \height(\alpha)^r \le A.
\]
Enumerate these bounded integer polynomials, collect their roots in \(K\), and filter coordinate pairs by membership in \(P\) and projective height at most \(B\). This gives precisely \(T_{K,B}(P)\).

\section{Proofs for the comb constructions}\label{sec:comb-deferred}

\begin{proof}[Proof of \cref{thm:central-comb}]
    The segment from a base guard \((x,y)\) to \((a_i,L)\) must cross the tooth opening at \(Y=2m\). Its two bounding rays give \(\abs{x-a_i}\le(L-y)/2\le L-1/2\). These intervals have centers \(2L\) apart and lengths \(2L-1\), so are disjoint. A guard above the base is in one tooth and cannot see another target across a gap. Thus \(m\) guards are necessary. The guards \((a_i,0)\) see their entire teeth and the base, attaining the classical bound; their vertical lines have height at most \(4m^2\).

    For the curved cover put \(r_i=2i-m+1\) and take \(\alpha_i(u)=(a_i(2u-u^2),Lu)\), \(0\le u\le1\). Until \(Lu=2m\) it stays in the base. Afterwards \(u\ge1-1/L\), so
    \[
        \abs{x(\alpha_i(u))-a_i} = \abs{a_i}(1-u)^2\le \abs{r_i}/L<1/2.
    \]
    It stays in its tooth, and its increasing ordinate makes it embedded. Its equation \(r_iY^2-2r_iLY+LX\) solves for \(X\) as a polynomial in \(Y\), with nonzero constant \(X\)-derivative. Its primitive part is nonsingular and absolutely irreducible, of height at most \(4m^2\); if \(r_i=0\), it is \(X\). These arcs join \((0,0)\) to every target, so one guard suffices.
\end{proof}

\begin{proof}[Proofs of \cref{thm:alternating-comb,cor:alternating-free-guards}]
    An injectively parameterized sight arc from \(q_i\) to \(q_j\), \(i<j\), first meets \(X=2i,2i+1,\ldots,2j\) in that order, even if it backtracks. Its ordinates alternate between at least \(1/2\) and at most \(-1/2\). The \(2(j-i)\) intervening parameter intervals give distinct intersections with \(Y=0\). For a carrier \(\mathrm{C}_F\), the polynomial \(F(X,0)\) is nonzero, since otherwise irreducibility would make the carrier \(Y=0\), which contains no target. Hence \(2(j-i) \le d\), giving the target-guard lower bound.

    For the reverse direction put \(r=2(j-i) \le 2s\) and use
    \[
        p_{ij}(X) = \sum_{k=0}^r (-2)^k\binom{X-2i}{k}.
    \]
    The binomial theorem gives \(p_{ij}(2i+u)=(-1)^u\) at integral nodes. On \([2i,2j]\) every falling-product factor has absolute value at most \(2m\), so
    \[
        \abs{p_{ij}} \le \sum_{k=0}^r \frac{(4m)^k}{k!} \le (r+1)(4m)^r<Y_m,
    \]
    and differentiation gives the explicit sum
    \[
        \abs{p'_{ij}} \le 2\sum_{h=0}^{r-1} \frac{(4m)^h}{h!} \le 2r(4m)^{r-1} \le \Lambda_m.
    \]
    The graph fits the outer rectangle and changes by at most \(1/4\) in a notch half-width. It therefore stays above \(1/2\) at even notch bands and below \(-1/2\) at odd ones. Clearing denominators gives
    \[
        r!Y-\sum_{k=0}^r (-2)^k \frac{r!}{k!} \prod_{\nu=0}^{k-1} (X-2i-\nu).
    \]
    Each product has coefficient \(\ell^1\)-norm at most \((2m)^k\), giving height at most \((r+1)r!(4m)^r \le H_{m,d}\). The primitive equation is a nonsingular absolutely irreducible graph. A middle target covers each consecutive block of at most \(2s+1\), giving the exact count; \(s=0\) uses reflexivity. If \(s=m-1\) the answer is one, and otherwise \(d \le 2s+1 \le d+1\), proving the uniform growth rate.

    For a general guard \((x,y)\), let \(N\) count the notch centers \(0,1,\ldots,2m-2\) in the closed interval between \(x\) and \(2i\). For every \(-1 \le x \le 2m-1\), we have \(N \ge \abs{x-2i}\). First encounters with these notch lines give at least \(N-1\) crossings of \(Y=0\), so \(\abs{x-2i} \le d+1\). The interval \([x-d-1,x+d+1]\) contains at most \(d+2\) targets. This proves the free-guard lower bound; the preceding target guards give its upper bound.
\end{proof}

\section{Proofs for convex connections and edge guards}\label{sec:edge-deferred}

\begin{proof}[Proof of \cref{lem:convex-connector}]
    If \(t\) is rational or \(g\) lies on its rational line, the segment from \(g\) to \(t\) suffices. Otherwise write \(t=b+\alpha v\) on its rational line and factor the minimal polynomial over \(\R\) as \(M(X)=(X-\alpha)R(X)\), with \(R(\alpha) \ne 0\). For rational \(a\) close to \(\alpha\), the affine map \(\Phi_a(u,w)=b+av+uv+w(g-b-av)\) is invertible. The equation \(M(a+u)-M(a)w=0\) is a nonsingular absolutely irreducible graph. Put \(\delta=\alpha-a\) and use
    \[
        u(s) = \delta s,\qquad w(s)=\frac{M(a+\delta s)}{M(a)} = (1-s) \frac{R(a+\delta s)}{R(a)},\qquad 0 \le s \le 1.
    \]
    On an interval about \(\alpha\), choose \(\abs{R} \ge c>0\) and \(\abs{R'} \le L_R\). The mean value theorem gives
    \[
        \abs{w(s)-(1-s)} = (1-s)\frac{\abs{R(a+\delta s)-R(a)}}{\abs{R(a)}} \le (L_R/c)\abs{\delta}s(1-s).
    \]
    The difference between \(\Phi_a(u(s),w(s))\) and \((1-s)g+st\) is \([w(s)-(1-s)](g-b-av)\).

    Write the convex polygon using its defining inequalities \(\lambda_j \ge 0\). No edge line contains both \(g,t\) as any two distinct rational lines intersect only at a rational point. Thus \(\lambda_j(g)+\lambda_j(t)>0\), and on the segment
    \[
        (1-s)\lambda_j(g)+s\lambda_j(t) \ge s(1-s)[\lambda_j(g)+\lambda_j(t)].
    \]
    Choose \(a\) so close that the error in every inequality is less than half this positive coefficient times \(s(1-s)\). The whole arc stays in \(P\), including near its endpoints. The strictly monotone coordinate \(u(s)\) makes it embedded. Clearing denominators and taking the primitive part preserves all these properties.
\end{proof}

\begin{proof}[Proof of \cref{prop:edge-rational}]
    For a rational target \(t\) visible from a closed edge \(e\), either \(t \in e\) and reflexivity applies, or a visible line through \(t\) contains another polygon vertex. To prove the latter, orient \(e\) and choose the first position \(g\) on it seeing \(t\); visible positions form a nonempty closed set. If \(g\) is an endpoint, it is the required vertex. Otherwise suppose the line contains no vertex other than possibly \(t\). The segment's interior cannot touch the boundary: contact away from a vertex would put an edge line, and thus its vertices, on the line.

    It remains to check that a small movement of \(g\) preserves containment near both endpoints. For each \(z \in \{g,t\}\) choose \(r_z>0\) so that inside \(B(z,r_z)\) the polygon agrees with its local cone of inward directions. At an edge interior this cone is a half-plane; at a vertex it is the sector between the incident edges. The direction from \(z\) toward the other endpoint lies in the interior of that cone, since a boundary direction would put an incident edge line on the sight line. At \(g\) the boundary edge itself is locally the fixed boundary of a half-plane. Thus these endpoint portions remain inside \(P\) when \(g\) moves sufficiently little along \(e\). The rest of the segment is compact and has positive distance from the boundary, so it also remains inside after the same movement. A small backward movement of \(g\) contradicts its first position. The resulting line through \(t\) and a rational vertex has height at most \(2VB\). Conversely an arc on a line contains its joining segment.

    An edge's segment-visibility region is closed: compactness supplies a convergent subsequence of observing points for converging visible targets, and their limiting segments stay in \(P\). The finite union for an edge collection is closed too. If it is not all of \(P\), its relatively open complement contains a rational target. This proves equality with the classical minimum number of edge guards. Choose one missed rational target for each collection of fewer than \(\eta_{\R}(P)\) edges, and let \(B_0\) bound their finitely many heights. Those collections fail for every \(B \ge B_0\); an optimal classical collection works at height \(2VB\).
\end{proof}

\begin{proof}[Proof of \cref{thm:edge-boundary}]
    Choose a point \(t_i\) of coordinate-field degree \(e\) in the interior of each edge \(E_i\). Rational translates of a primitive element of \(K\) give degree-\(e\) parameters in any open interval; a rational parameterization of an edge preserves their generated field. Let \(B_0\) bound the heights of these finitely many targets. Restricting a rational polynomial of degree below \(e\) to the edge line gives a univariate polynomial vanishing at a degree-\(e\) parameter. It is identically zero, so the edge line divides the equation. Absolute irreducibility makes the line the entire carrier. Convexity and maximal edges imply that its trace in \(P\) is precisely \(E_i\).

    A selected edge sees only the targets lying on it or on either adjacent edge. Each selected edge therefore covers at most three of the \(n\) targets, giving \(\lceil n/3 \rceil\). Conversely, partition the cyclic edge list into groups of at most three consecutive edges and select a middle edge from each group. At height \(h_P\), these selected edges cover every boundary point using its edge line, including from a shared vertex when necessary.

    At degree \(e\), choose a rational midpoint on every edge. Each boundary target lies on a rational line and has coordinate-field degree at most \(e\), so \cref{lem:convex-connector} joins it to every midpoint. Taking the maximum of \(h_P\) and the finitely many carrier heights gives \(H_*\). For the slice allowing all coordinate fields of degree at most \(e\), retain the same degree-\(e\) lower-bound targets and use Northcott finiteness for the upper bound.
\end{proof}

\begin{proof}[Proof of \cref{prop:edge-chords}]
    First, a nonrational quadratic target has exactly one rational line. The relations among \(1,x(t),y(t)\) form a one-dimensional rational space. The line joins \(t\) to its distinct conjugate. At each finite place the cross-product bound has constant one, and at archimedean places constant two. Multiplying normalized local bounds gives line height at most \(2\height(t)^2\). Its projective coefficient point is rational, so this is the height of its primitive integer equation.

    If an edge collection meets all proper chords, it meets every connected component of \(P \cap \ell\), for every real line \(\ell\). Indeed, approximate a point on that line by interior points and the line by lines avoiding vertices and edge directions. Their components are proper chords. Selected endpoints have a subsequence converging on one selected closed edge; the limiting segments stay in \(P\) and on the original line. This includes components through vertices and on the boundary.

    Every pair of edges joined by a proper interior chord admits one with rational endpoints: truncate short endpoint portions, use the middle segment's positive distance from the boundary, and preserve the transverse endpoint intersections under small rational perturbations. Choose one such \([a,b]\) per pair and put \(t=(1-\xi)a+\xi b\) for a fixed irrational \(\xi \in K \cap (0,1)\). Its unique rational line has precisely this component and its two boundary endpoints. An edge cover must therefore select at least one edge from each such pair, requiring \(b(P)\) edges. Let \(B_0\) bound these finitely many targets' heights. This lower bound applies to all slices once \(B \ge B_0\).

    Conversely, a minimum chord cover meets the line component of every quadratic target. Its line height is at most \(2B^2\); rational targets may use their vertical lines of height at most \(B\). This gives the finite-height and whole-gallery upper bounds. Any \(n-1\) edges meet every proper chord. For the degree-at-most-two version keep the same lower-bound targets from one field and apply the upper argument to each target individually.
\end{proof}

\begin{proof}[Proof of \cref{cor:star-boundary-height}]
    Since \(g\) sees every point of the edge, the triangle \(\operatorname{conv}(g,b,c)\) lies in \(P\). The rational and collinear cases use the stated lines and \cref{prop:linear}. Suppose now that the triangle is nondegenerate. Take its rational affine coordinates in which \(b=(0,0)\), \(c=(1,0)\), and \(g=(0,1)\). Write these coordinate functions in the original plane as
    \[
        x = \frac{L_x(X,Y)}D,\qquad y = \frac{L_y(X,Y)}D,
    \]
    where \(L_x,L_y\) are integral affine polynomials and \(D\) is a positive integer. Put
    \[
        W = \max\{D,\norm{L_x}_1,\norm{L_y}_1\}.
    \]
    Let \(M \in \Z[X]\) be the primitive minimal polynomial of \(\alpha\), with coefficient height \(A\). Define
    \[
        \rho=[2A(A+1)]^{-e^2},\qquad \epsilon_{\alpha}=\frac1{A(A+2)^{e-1}},
    \]
    then put
    \[
        c_0=(\rho/2)^{e-1},\qquad L_0=A(e-1)(A+2)^{e-2}, 
        \qquad Q=1+\left\lceil \max\left\{ \frac{4}{\rho},\frac{4}{\epsilon_{\alpha}},\frac{4L_0}{c_0\epsilon_{\alpha}} \right\} \right\rceil.
    \]
    The resulting carrier satisfies
    \begin{equation}\label{eq:star-height-bound}
        \norm{F}_\infty \le 2A(e+1)(2QW)^e.
    \end{equation}

    \emph{Choose a rational approximation with room at both endpoints.} All roots of \(M\) have modulus at most \(A+1\). If two roots are at distance \(\sigma\), then given that the integral discriminant of \(M\) is nonzero, we have 
    \[
        1 \le A^{2e-2}\sigma^2[2(A+1)]^{e(e-1)-2}.
    \]
    Hence their distance is at least \(\rho\). Factoring the nonzero integers \(M(0),M(1)\) also gives \(\alpha,1-\alpha \ge \epsilon_{\alpha}\). Choose \(a=\lfloor Q\alpha\rfloor/Q\), so \(0<a<\alpha<1\) and \(\delta \coloneqq \alpha-a<1/Q\). In the factorization \(M(T)=(T-\alpha)R(T)\), every other root stays at least \(\rho/2\) from \([a,\alpha]\). Thus \(\abs{R} \ge c_0\) and \(\abs{R'} \le L_0\) on this interval.

    \emph{Keep the connecting arc in the triangle.} Use the equation
    \[
        M(x+ay)-M(a)y = 0.
    \]
    The invertible change \((x,y) \mapsto (x+ay,y)\) makes it a polynomial graph, with nonzero constant derivative in its second coordinate. It is therefore nonsingular and absolutely irreducible, of degree \(e\). Its arc from \(g\) to \(t\) is
    \[
        y(s)=(1-s)\frac{R(a+\delta s)}{R(a)},\qquad x(s)=a+\delta s-a y(s),\qquad 0 \le s \le 1.
    \]
    Set \(E(s)=y(s)-(1-s)\). By the mean value theorem,
    \[
        \abs{E(s)} \le (L_0/c_0) \delta s(1-s) \le (\epsilon_{\alpha}/4)s(1-s).
    \]
    Consequently \(y\ge0\), and the other barycentric coordinates satisfy
    \[
        x = \alpha s-aE \ge s(\alpha-\epsilon_{\alpha}/4),\qquad 1-x-y = (1-\alpha)s-(1-a)E \ge s(1-\alpha-\epsilon_{\alpha}/4).
    \]
    All three are nonnegative, so the arc stays in the triangle. The coordinate \(x+ay=a+\delta s\) is strictly increasing, making the arc embedded.

    \emph{Bound the coefficients and check the degree.} Write \(a=k/Q\) and \(M(T)=\sum_jm_jT^j\). Multiplication by \(Q^eD^e\) clears the equation's denominators. Its first term becomes
    \[
        \sum_{j=0}^e m_jQ^{e-j}D^{e-j}(QL_x+kL_y)^j,
    \]
    of coefficient \(\ell^1\)-norm at most \(A(e+1)(2QW)^e\). The second term is \(Q^eM(a)L_yD^{e-1}\), integral and of norm at most \(A(e+1)(QW)^e\). Taking the primitive part gives \cref{eq:star-height-bound}. Finally, an equation of degree less than \(e\) through \(t\) vanishes identically on the edge line, since its restriction has \(\alpha\) as a root. An absolutely irreducible such carrier would be that line, which misses \(g\). Thus degree \(e\) is least possible when \(g\) is off the line.
\end{proof}

\begin{proof}[Proof of \cref{prop:edge-quadratic-height}]
    Write \(\alpha=\sqrt u\). For \(k\ge2\), put
    \[
        a_k = \lfloor (k-1)\alpha \rfloor,\qquad b_k = \lfloor k\alpha \rfloor,
    \]
    and let
    \[
        t_k = ((k-1)\alpha-a_k,k\alpha-b_k).
    \] 
    Both coordinates are strictly between zero and one. The unique rational line through \(t_k\) has equation
    \[
        kX-(k-1)Y+c_k = 0,\qquad c_k = ka_k-(k-1)b_k.
    \]
    Its coefficients are primitive because \(\gcd(k,k-1)=1\), and \(c_k=(k-1)y(t_k)-kx(t_k)\) lies strictly between \(-k\) and \(k-1\). Its coefficient height is exactly \(k\).

    The target coordinates are algebraic integers, so every finite-place factor of \(\height(t_k)\) is one. At the chosen real embedding 
    \[
        \max\{1,\abs{x(t_k)},\abs{y(t_k)}\}
    \]
    equals one, and at the other it is \(k\alpha+b_k\), larger than \((k-1)\alpha+a_k\) and one. Therefore
    \[
        \height(t_k)^2 = k \alpha+\lfloor k \alpha \rfloor < 2\alpha k.
    \]
    For sufficiently large \(B\), choose \(k=\lfloor B^2/(2\alpha) \rfloor \ge 2\). Then \(t_k \in T_{K,B}(\mathcal{S})\). Below coefficient height \(k\), its only rational line is unavailable. Since \(t_k\) is interior, no edge covers it by reflexivity either. This proves the lower bound on \(H_{\mathrm{cov}}(B)\); \cref{prop:edge-chords} gives the upper bound.
\end{proof}

\begin{proof}[Proof of \cref{prop:edge-alternating-sparse}]
    Every edge has horizontal span at most two. By \cref{eq:alternating-reach}, its visible targets occupy an interval of length at most \(2d+4\), containing at most \(d+3\) targets.

    Partition the targets into blocks of at most \(s+1\). For a block starting at \(q_i\), select the bottom-notch tip edge at \(2i\) and its right endpoint \(v_i=(2i+\delta_m,1/2)\). Its segment to \(q_i\) lies above the notch and has equation \(Q_mX+2Y-(2iQ_m+2)=0\). For a later \(q_j\) in the block put \(r=2(j-i)\le2s\) and interpolate
    \[
        z_0 = 2i+\delta_m,\quad y_0 = 1/2;\qquad z_k = 2i+k,\quad y_k = (-1)^k,\quad 1 \le k \le r.
    \]
    After translation, the nodes are \(\delta_m,1,\ldots,r\). The first Lagrange denominator has magnitude \(\prod_{k=1}^r(k-\delta_m) \ge r!(1-\delta_m)^r \ge r!(1-r\delta_m) \ge r!/2\). The \(k\)th remaining one is \((k-\delta_m)(k-1)!(r-k)! \ge k!(r-k)!/2 \ge 1/2\). On the interpolation interval, all numerator factors have magnitude at most \(2m\), giving
    \[
        \abs{p_{ij}} \le 2(r+1)(2m)^r < Y_m,\qquad \abs{p'_{ij}} \le 2r(r+1)(2m)^{r-1} \le \Lambda_m.
    \]
    The graph starts on the right side of its initial notch. At each later notch center it has value \(\pm1\), and the derivative bound retains the required sign throughout its band. It therefore stays in \(\mathcal{A}_m\).

    A common coefficient denominator is \(D_r=2r!\prod_{k=1}^r(Q_mk-1)\) for these factors clear the first Lagrange term, while the remaining terms have denominators dividing \((Q_mk-1)(k-1)!(r-k)!\) after cancellation. The numerator coefficient norms give \(\norm{p_{ij}}_1 \le 2(r+1)(2m+1)^r\). Consequently
    \[
        \norm{D_r(Y-p_{ij}(X))}_\infty \le 4r!(r+1)(2m+1)^r\prod_{k=1}^r(Q_mk-1) \le H^{\mathrm{edge}}_{m,d}.
    \]
    Its primitive part is a nonsingular absolutely irreducible graph. This covers each block, including \(s=0\) by the initial segments. The growth rate follows from \(d+3\le4d\) and \(s+1 \ge d/2\), unless \(s=m-1\), when the upper bound is one.
\end{proof}

\begin{proof}[Proof of \cref{cor:edge-eventual}]
    On the chosen edge take a positive-length rational-endpoint segment \(J\) away from its ends. A thin rational triangle \(\Delta\) can be attached outside \(P\) with \(\Delta\cap P=J\), making \(P'=P \cup \Delta\) simple. Choose \(g \in (\operatorname{int}\Delta)(\Q)\). For each \(t \in T \setminus J\), apply \cref{thm:effective-lifting} to \(P',g,t\). This gives a primitive polynomial \(F_t\) with computable degree and coefficient height whose regular trace in \(P'\) contains an arc from \(g\) to \(t\). Targets on \(J\) use reflexivity. For any other target, the final portion is outside \(\Delta\) and inside \(P\). Let \(s_*\) be the supremum of the parameters at which the arc lies outside \(P\). Closedness and continuity put the point at \(s_*\) on \(J\) and the remaining arc in \(P\).

    Restricting the arc does not change its carrier polynomial. Repeat this construction for each edge, and take \(d_*,H_* \ge 1\) to bound, respectively, the degrees and coefficient heights of all the polynomials obtained. These bounds are computable and give budgets that are valid for every edge.
\end{proof}

\section{Uniform counts for rational points and lines}\label{sec:lattice-proofs}

\cref{thm:linear-threshold} compares the number of targets with the number that a bounded set of guards can reach. The estimates must be uniform since both the guards and their sight arcs may change with \(B\). We again work over \(\Q\) and write \(T_B(P)=T_{\Q,B}(P)\). We continue to use the sup norm on primitive homogeneous coordinates.

\begin{lemma}\label{lem:line-uniform}
    Let \(u=(a,b,c) \in \Z^3\) be primitive, and put \(h=\norm{u}_\infty\).  The lattice
    \[
        L_u=\{v \in \Z^3 : u\cdot v = 0\}
    \]
    has Euclidean covolume \(\norm{u}_2\), the area of a fundamental parallelogram in \(u^\perp\). If \((a,b) \ne (0,0)\) and \(1 \le h \le B\), the line \(aX+bY+c=0\) contains \(O(B^2/h)\) rational points of projective height at most \(B\), with an absolute implied constant.
\end{lemma}

\begin{proof}
    Primitivity gives \(w \in \Z^3\) with \(u \cdot w=1\).  If \(v_1,v_2\) form a basis of \(L_u\), then \(v_1,v_2,w\) form a basis of \(\Z^3\).  Their determinant is \(\pm1\).  Since \(v_1 \times v_2\) is parallel to \(u\), its scalar product with \(w\) shows that \(v_1 \times v_2=\pm u\).  This proves the covolume identity and
    \[
        h \le \operatorname{covol}(L_u) = \norm{u}_2 \le \sqrt{3}\, h.
    \]

    We now bound the number of integer vectors in \(L_u\) whose coordinates have absolute value at most \(B\). This also bounds the number of rational points on the line of projective height at most \(B\). Permute the three coordinate positions so that the coefficient of the coordinate being eliminated has absolute value \(h\).  Denote the other coefficients by \(a',b'\).  The remaining two integer coordinates \(x,y\) must satisfy
    \[
        a'x+b'y \equiv 0\pmod h,\qquad \abs{x},\abs{y} \le B.
    \]
    Put \(g=\gcd(a',h)\).  Primitivity implies \(\gcd(b',g)=1\), so \(g\) divides \(y\).  There are at most \(2B/g+1\) possibilities for \(y\).  For each one, \(x\) lies in a single residue class modulo \(h/g\), giving at most \(2Bg/h+1\) choices.  The eliminated coordinate is then uniquely determined; ignoring its size only enlarges the count.  Thus
    \[
        \card{(L_u \cap [-B,B]^3)} \le (2B/g+1)(2Bg/h+1) = O(B^2/h+B+1).
    \]
    Because \(h \le B\), this is \(O(B^2/h)\).  Primitive homogeneous representatives of the desired rational points form a subset, which proves the assertion.
\end{proof}

\begin{lemma}\label{lem:guard-lines}
    For any rational affine point \(g\) and \(h \ge 1\), the number of primitive integral line equations through \(g\), up to sign, with coefficient height in \((h/2,h]\) is \(O(h^2)\), uniformly in \(g\).
\end{lemma}

\begin{proof}
    Write \(g=(a/c,b/c)\) with \(c \ne 0\).  A coefficient vector \((u,v,w)\) satisfies
    \[
        au+bv+cw = 0.
    \]
    The integers \(u,v\) determine \(w\) uniquely.  There are at most \((2\lfloor h \rfloor+1)^2\) possibilities with \(\abs{u},\abs{v} \le h\). Imposing the bound on \(w\), primitivity, the interval condition, and division by sign can only reduce this count.
\end{proof}

\begin{proof}[Proof of \cref{thm:linear-threshold}]
    The upper bound follows from \cref{cor:fisk}, with the vertex heights absorbed in a constant depending on \(P\).

    For the lower bound, first suppose \(1\le H\le B\). Split the rational lines in the pencil through \(g\), of coefficient height at most \(H\), into dyadic height bands: those with heights in \((h/2,h]\), for \(h=H,H/2,H/4,\ldots\).  Discard empty bands; a nonempty one has \(h \ge 1\).  In each band, \cref{lem:line-uniform,lem:guard-lines} bound the number of targets on all its lines by
    \[
        O(h^2)\,O(B^2/h) = O(B^2h).
    \]
    Here an individual line has height greater than \(h/2\), so the first lemma gives \(O(B^2/h)\) with an absolute constant.  Since
    \[
        \sum_{j\ge0} \frac{H}{2^j} = 2H,
    \]
    summing the band estimates gives \(O(B^2H)\) targets visible from one guard.  The guard itself contributes at most one further target, absorbed by the same bound.  Counting whole lines ignores polygonal occlusion and therefore remains a valid upper bound. Consequently \(k\) guards cover at most \(O(kB^2H)\) targets.

    To bound the total number of targets below, consider homogeneous triples \((x,y,z)\) with \(z>0\), \((x/z,y/z)\in P\), and \(\max\{\abs{x},\abs{y},z\}\le B\).  They lie in the \(B\)-dilate of a fixed positive-volume region that is a finite union of rational polytopes: triangulate \(P\) and take the truncated cone over each triangle. Their intersections lie in lower-dimensional faces, so counting integer points in the union gives \(v_P B^3+O_P(B^2)\), with \(v_P>0\). Removing common divisors by M\"obius inversion then gives
    \[
        \card{T_B(P)} = c_P B^3+O_P(B^2),\qquad c_P = \frac{v_P}{\zeta(3)} > 0.
    \]
    Indeed, the error sums to \(O_P(B^2\sum_{r\ge1}r^{-2})\), and extending the main sum from \(r\le B\) to all positive integers adds only \(O_P(B)\). For sufficiently large \(B\), there are therefore at least \((c_P/2)B^3\) targets.  If \(C\) is the absolute constant in the one-guard estimate, a cover by \(k\) guards must satisfy
    \[
        CkB^2H \ge \frac{c_P}{2}B^3,\qquad \text{hence} \qquad H \ge \frac{c_P}{2Ck}B.
    \]
    When \(H>B\), the same conclusion holds after replacing the last constant by its minimum with \(1\).  This is the required positive constant depending only on \(P\) and \(k\).
\end{proof}

\section{Arc lifting: the common construction and its bounds}\label{sec:relative-proofs}\label{sec:bertini-proofs}\label{sec:effective-proof}

We prove \cref{thm:effective-lifting}; since its conclusion includes that of \cref{thm:lifting}, the qualitative theorem follows \emph{a fortiori}. The geometric task is to build a regular reference arc with the correct endpoint germs, and the arithmetic task is to replace its defining function by a rational polynomial and then choose a smooth projective member without losing that arc. The first jets fix the endpoints during approximation. At the last stage only their values remain fixed, while their derivatives remain close.

\colorlet{liftblue}{cyan!65!blue}
\colorlet{liftorange}{orange!85!black}
\colorlet{liftpurple}{purple!70!black}
\colorlet{liftgreen}{green!50!black}

\tikzset{
    liftboundary/.style={
        fill=gray!8,draw=black,thick,
        line cap=round,line join=round
    },
    liftarc/.style={
        draw=liftblue,very thick,
        line cap=round,line join=round
    },
    liftnew/.style={
        draw=liftorange,very thick,densely dashed,
        line cap=round,line join=round
    },
    liftloop/.style={
        draw=liftblue,line width=.75pt,
        line cap=round,line join=round
    },
    liftgerm/.style={
        draw=liftpurple,line width=2pt,
        line cap=round,line join=round
    },
    liftpoint/.style={
        circle,fill=black,inner sep=0pt,
        minimum size=4pt
    },
    liftguard/.style={
        rectangle,fill=black,inner sep=0pt,
        minimum size=4pt
    },
    liftlabel/.style={}
}

\newcommand{\LiftingPolygon}{%
    \path[liftboundary]
        (0,0)--(4,0)--(4,3)--(2.4,3)
        --(2.4,1.75)--(0,1.75)--cycle;
}

\newcommand{\LiftingArc}[1][]{%
    \draw[liftarc,#1]
        (.65,.55)--(1.05,.55)
        .. controls (1.8,.55) and (3.25,.7) .. (3.25,1.85)
        --(3.25,3);
}

\newcommand{\LiftingLoop}[1][]{%
    \path[liftloop,#1]
        (.45,.55)--(1.05,.55)
        .. controls (1.8,.55) and (3.25,.7) .. (3.25,1.85)
        --(3.25,3.16)
        .. controls (3.25,3.48) and (3.65,3.48) .. (3.65,3.15)
        --(3.65,1.85)
        .. controls (3.65,.18) and (1.5,.13) .. (.45,.13)
        .. controls (.15,.13) and (.15,.55) .. (.45,.55)
        --cycle;
}

\newcommand{\LiftingAuxLoop}[1][]{%
    \path[liftloop,#1]
        (.7,2.9)--(.7,3.1)
        .. controls (.7,3.42) and (1.3,3.42) .. (1.3,3.1)
        --(1.3,2.9)
        .. controls (1.3,2.58) and (.7,2.58) .. (.7,2.9)
        --cycle;
}

\newcommand{\LiftingMarks}{%
    \node[liftguard,label={[liftlabel]below left:\(g\)}]
        at (.65,.55) {};
    \node[liftpoint,label={[liftlabel]above left:\(t\)}]
        at (3.25,3) {};
    \node[liftpoint,label={[liftlabel]below left:\(t'\)}]
        at (.7,3) {};
}

\newcommand{\LiftingGerms}{%
    \draw[liftgerm] (.43,.55)--(.97,.55);
    \draw[liftgerm] (3.25,2.76)--(3.25,3.24);
    \draw[liftgerm] (.7,2.91)--(.7,3.09);
}

\subsection{Interpolation and geometric preparation}\label{sec:effective-preparation}
The first jet is \(j_p^1F = (F(p),dF_p)\). A jet \((0,\nu_p)\), with \(\nu_p \ne 0\), specifies a regular curve through \(p\) with tangent line \(\ker\nu_p\), not its full local germ.

\begin{lemma}\label{lem:equivariant-jets}
    A finite Galois-stable subset \(S \subset \Qbar^2\), with Galois-equivariant algebraic values and differentials, admits a rational polynomial interpolating those first jets. If \(N=\card{S}\), degree at most \(2N-1\) suffices.
\end{lemma}

\begin{proof}
    Over a finite Galois extension \(E/\Q\) containing the data, the ideals \(\mathfrak{m}_p=(X-p_x,Y-p_y)\) are distinct and their squares are comaximal. By Sunzi's theorem, prescribing a first jet at each \(p \in S\) therefore amounts to choosing a residue class modulo every \(\mathfrak{m}_p^2\), and such classes can always be matched by a single polynomial. We construct one explicitly in order to bound its degree. Fix \(p \in S\) and choose, for each other point, a linear polynomial vanishing there but not at \(p\). The product of their squares has degree \(2N-2\) and zero first jets at every other point. Since it is nonzero at \(p\), an affine linear factor can be chosen to give the prescribed value and differential there. The resulting polynomial has degree at most \(2N-1\). Sum these polynomials over \(p \in S\), then average over \(\operatorname{Gal}(E/\Q)\). Galois-equivariance preserves the jets and makes the coefficients rational.
\end{proof}

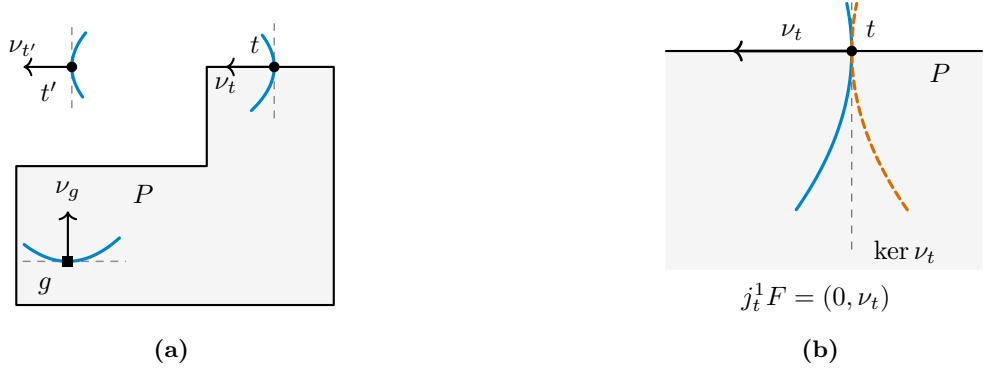
\begin{figure}[htbp]
    \centering
    \begin{subfigure}[t]{.46\linewidth}
        \centering
        \begin{tikzpicture}[scale=1.05,every node/.style={font=\small}]
            \LiftingPolygon
            \node at (1.6,1.4) {\(P\)};
            \draw[black!45,line width=.55pt,dash pattern=on 3pt off 2.5pt] (.08,.55)--(1.38,.55);
            \draw[black!45,line width=.55pt,dash pattern=on 3pt off 2.5pt] (3.25,2.35)--(3.25,3.55);
            \draw[black!45,line width=.55pt,dash pattern=on 3pt off 2.5pt] (.7,2.48)--(.7,3.52);
            \draw[liftarc] plot[variable=\s,domain=-.55:.65,samples=55]
            ({.65+\s},{.55+.7*\s*\s});
            \draw[liftarc] plot[variable=\s,domain=-.55:.4,samples=55]
            ({3.25-.95*\s*\s},{3+\s});
            \draw[liftarc] plot[variable=\s,domain=-.38:.43,samples=55]
            ({.7+.95*\s*\s},{3+\s});
            \draw[->,line width=.8pt] (.65,.55)--(.65,1.18) node[above] {\(\nu_g\)};
            \draw[->,line width=.8pt] (3.25,3)--(2.64,3) node[below] {\(\nu_t\)};
            \draw[->,line width=.8pt] (.7,3)--(.09,3) node[above] {\(\nu_{t'}\)};
            \LiftingMarks
            \path[use as bounding box] (-.25,-.2) rectangle (4.2,3.65);
        \end{tikzpicture}
        \caption{}\label{fig:lifting-jets-a}
    \end{subfigure}\hfill
    \begin{subfigure}[t]{.50\linewidth}
        \centering
        \begin{tikzpicture}[scale=1.05,every node/.style={font=\small}]
            \fill[gray!8] (0,0) rectangle (4,2.75);
            \draw[black,thick] (0,2.75)--(4,2.75);
            \node[anchor=west] at (3.2,2.48) {\(P\)};
            \draw[black!65,dashed] (2.35,.25)--(2.35,3.45);
            \draw[liftarc] plot[variable=\s,domain=-2:.6,samples=50]
            ({2.35-.175*\s*\s},{2.75+\s});
            \draw[liftnew] plot[variable=\s,domain=-2:.6,samples=50]
            ({2.35+.175*\s*\s},{2.75+\s});
            \draw[->,line width=1pt] (2.35,2.75)--(.85,2.75) node[midway,above] {\(\nu_t\)};
            \node[liftpoint,label={above right:\(t\)}] at (2.35,2.75) {};
            \node[anchor=west] at (2.5,.22) {\(\ker\nu_t\)};
            \node at (1.9,-.35) {\(j_t^1F=(0,\nu_t)\)};
            \path[use as bounding box] (-.2,-.55) rectangle (4.15,3.45);
        \end{tikzpicture}
        \caption{}
        \label{fig:lifting-jets-b}
    \end{subfigure}

    \caption{Local data for arithmetic arc lifting. (a) The endpoints \(g,t\) and another real conjugate \(t'\) have fixed zero values and nonzero covectors, shown as normal arrows. Separate solid blue local branches are compatible with the prescribed first jets and dashed gray lines show their tangents. (b) Two branches can share the same first jet without sharing a germ.}
    \label{fig:lifting-jets}
\end{figure}

Assume \(g \ne t\), put \(K=\Q(t_x,t_y)\), \(e=[K:\Q]\), and let \(S\) consist of \(g\) and the \(e\) conjugates of \(t\). Then \(N=e+1\). Choose nonzero rational tangent vectors at \(t\) and \(g\), taking the vector at a boundary target strictly into \(P\). Then take short straight segments from these points in the chosen directions. Connect their other endpoints using rational intermediate waypoints so that the resulting polygonal path from \(t\) to \(g\) has no self-intersections, lies in \(\operatorname{int}P\) except possibly at \(t\), and avoids the other real points of \(S\). Orient this path from \(t\) to \(g\) and choose rational normal covectors \(\nu_t,\nu_g\) whose coefficient vectors are obtained by rotating the path's oriented endpoint tangent vectors clockwise through \(\pi/2\). At every conjugate \(p\) of \(t\), set \(\nu_p=\nu_t\), then assign the first jet \((0,\nu_p)\) at every \(p \in S\). These data are Galois-equivariant.

For the coefficient bounds, choose an integral primitive element \(\theta\) and write
\[
    t = (t_1(\theta)/D,t_2(\theta)/D),\qquad g = (g_1/D_g,g_2/D_g),
\]
where \(t_1,t_2 \in \Z[X]\), \(\deg t_1,\deg t_2<e\), and \(D,D_g>0\). Let \(M_\theta\) be the multiplication matrix of \(\theta\) in the power basis. Choose an integer
\[
    L \ge \max\{1,D,D_g,\abs{g_1},\abs{g_2}, \norm{t_1(M_\theta)}_{\mathrm{row}}, \norm{t_2(M_\theta)}_{\mathrm{row}}\},
\]
where the matrix norm is the maximum absolute row sum. For \(j \ge N-1\), let \(\mathcal{E}_j\) be the integral evaluation matrix obtained by scaling the target evaluation equations by \(D^j\) and the guard evaluation equation by \(D_g^j\). This matrix has rank \(N\), with entries bounded by \(L^j\). For example, the target block of the column corresponding to \(X^aY^b\) is \(D^{j-a-b}t_1(M_\theta)^at_2(M_\theta)^b(1,0,\ldots,0)^{\mathsf{T}}\). To verify the rank, extend scalars to a finite Galois extension \(E/\Q\) containing the coordinates of every point of \(S\). The embeddings of \(K\) into \(E\) give \(K \otimes_{\Q}E \cong E^e\). Under this identification, the target block in power-basis coordinates becomes evaluation at the \(e\) conjugates of \(t\). This is an invertible change of rows, and the scaling factors \(D^j,D_g^j\) are nonzero, so the rank is unchanged. For each \(p \in S\), choose a line through each other point that does not pass through \(p\). The product of their defining linear polynomials has degree \(N-1\), vanishes at every other point, and is nonzero at \(p\). Rescale it to have value one at \(p\). These polynomials give all the standard basis vectors under evaluation at \(S\), so the evaluation map has rank \(N\) over \(E\), and therefore over \(\Q\). 

Now put \(J_j = \binom{j+2}{2}\), the number of columns of \(\mathcal{E}_j\). Choose \(N\) pivot columns whose \(N \times N\) submatrix has nonzero determinant. For each free coordinate, set that coordinate equal to the pivot determinant and the other free coordinates to zero, then solve for the pivot coordinates by Cramer's rule. The resulting integral polynomials \(B_i^{(j)}\) form a rational basis of the evaluation kernel, with
\[
    \norm{B_i^{(j)}}_\infty \le b_j \coloneqq N!L^{jN},\qquad \norm{B_i^{(j)}}_1 \le J_jb_j.
\]
Here \(\norm{\cdot}_1\) is the sum of the absolute values of the coefficients. Multiply all the covectors by a common positive integer to make their coefficients integral, retaining the notation \(\nu_p\), and let \(c_\nu\) bound the absolute values of their coefficients. The first-jet matrix in degree \(m_0 \coloneqq 2N-1\) has rank \(3N\) by \cref{lem:equivariant-jets}. In power-basis coordinates, multiply the target value and derivative equations, including their right-hand sides, by \(D^{m_0}\), and the guard equations by \(D_g^{m_0}\). The resulting matrix and right-hand side have integral entries of absolute value at most \(\max(m_0,c_\nu)L^{m_0}\). Applying Cramer's rule to the first jets gives \(F_0=P_0/D_0\) of degree at most \(m_0\), satisfying \(F_0(p)=0\) and \(d(F_0)_p=\nu_p\) for every \(p \in S\). Set \(C_0 = \norm{P_0}_\infty\). These integers satisfy
\[
    \max\{D_0,C_0\} \le (3N)!\,[\max(m_0,c_\nu)L^{m_0}]^{3N}.
\]
The construction uses their actual integer values.

Use the degree-\(N\) basis, with \(k_N=J_N-N\), and set
\[
    A_0 = \sum_{i=1}^{k_N}(B_i^{(N)})^2,\qquad A_* = k_N(J_Nb_N)^2.
\]
The ideal defined by \(S\) is \(I(S) \coloneqq \{F \in \Q[X,Y] \colon F(p)=0 \text{ for all } p \in S\}\). The polynomial \(A_0\) lies in \(I(S)^2\), has degree at most \(2N\), coefficient sum at most \(A_*\), and vanishing first jets at every point of \(S\). Moreover, its real zero set is precisely \(S \cap \R^2\). Indeed, every \(B_i^{(N)}\) vanishes on \(S\), so \(A_0\) vanishes at its real points. Conversely, for a real point \(w \notin S\) and for each \(p \in S\), choose a linear polynomial over \(E\) vanishing at \(p\) but not at \(w\), and multiply these polynomials. Their product has degree \(N\) and belongs to the evaluation kernel after extending scalars to \(E\), so it is an \(E\)-linear combination of the \(B_i^{(N)}\). Since the product is nonzero at \(w\), at least one \(B_i^{(N)}(w)\) is nonzero. Hence \(A_0(w)>0\).

\begin{lemma}[Computable reference function]\label{lem:computable-f}
    From \(P,g,t,F_0\) as above one can construct a semialgebraic \(\mathscr{C}^2\) function \(f\) on \(\R^2\) with compact regular zero set, agreeing with \(F_0\) on neighborhoods of the points of \(S \cap \R^2\). Its zero set contains an embedded semialgebraic \(\mathscr{C}^3\) arc \(a \colon[0,1] \to P\), with \(a(0)=t\), \(a(1)=g\), and \(a((0,1]) \subset \operatorname{int}P\). If \(t \in \partial P\), then \(a'(0)\) points strictly into \(P\). The construction computes finite semialgebraic descriptions of \(f\) and \(a\) over the real algebraic numbers.
\end{lemma}

\begin{proof}
    \emph{Join the actual endpoint germs.} At each real point of \(S\), at least one partial derivative of \(F_0\) is nonzero. The implicit function theorem therefore expresses the zero set of \(F_0\) in a sufficiently small box as the graph of a real-analytic function of the other coordinate. This function is algebraic because it satisfies the polynomial equation \(F_0=0\). Around each point, search for a box with dyadic rational endpoints on which the chosen partial derivative has a constant nonzero sign and \(F_0\) has opposite signs on the two sides perpendicular to that coordinate direction. Every segment joining those sides parallel to that direction then contains exactly one zero of \(F_0\). These conditions can be checked by exact semialgebraic tests, and the implicit function theorem guarantees that suitable boxes exist. Choose the boxes small enough to be pairwise disjoint. At \(t\) and \(g\), retain the halves of the local branches directed toward the chosen polygonal path; when \(t\) lies on the boundary, its retained half enters \(\operatorname{int}P\).

    Use the polygonal path chosen above, replacing its straight endpoint portions near \(t\) and \(g\) by the local branches of \(F_0=0\). Each branch has the same tangent direction as the segment it replaces. The smoothing construction below joins these branches to the remaining polygonal path while smoothing its corners. Choose the neighborhoods in which these modifications occur to be disjoint and small enough to ensure that the arc is embedded in \(P\) and avoids other points of \(S\).

    We round each corner by interpolating between the tangent vectors of its two adjacent segments using the following \(\mathscr{C}^5\) transition function:
    \[
        \chi(u) = 462u^6-1980u^7+3465u^8-3080u^9+1386u^{10}-252u^{11}, \quad 0 \le u \le 1,
    \]
    where we set \(\chi(u)=0\) for \(u<0\) and \(\chi(u)=1\) for \(u>1\). The derivatives of \(\chi\) of orders \(1,\ldots,5\) vanish at both endpoints, so the extension is \(\mathscr{C}^5\). On \([0,1]\), we have \(\chi'(u)=2772u^5(1-u)^5\), and \(\int_0^1\chi(u)\, du=1/2\). At a corner \(z\), let \(v_-\) and \(v_+\) be the constant tangent vectors of affine parameterizations of the incoming and outgoing segments, both oriented along the path. For sufficiently small \(h>0\), replace the portions from \(z-hv_-\) to \(z\) and from \(z\) to \(z+hv_+\) by
    \[
        \gamma(s)=z-hv_- + sv_- + (v_+ - v_-)\int_0^s \chi(u/(2h))\, du \quad 0 \le s \le 2h.
    \]
    The integral identity gives \(\gamma(0)=z-hv_-\) and \(\gamma(2h)=z+hv_+\). The endpoint values and vanishing derivatives of \(\chi\) ensure that \(\gamma\) matches the adjacent segments in derivatives through order five. Its tangent vector is
    \[
        \gamma'(s)=(1-\chi(s/(2h)))v_- + \chi(s/(2h))v_+.
    \]
    The vectors \(v_-\) and \(v_+\) are not oppositely directed, so there is a vector \(v_0\) having positive scalar product with both. Thus \(\langle \gamma(s),v_0 \rangle\) is strictly increasing, which shows that \(\gamma'\) never vanishes and that \(\gamma\) has no self-intersections.

    To join a retained branch of \(F_0=0\) to a straight segment, write both as graphs over their common tangent direction. On a short overlap interval, interpolate their graph functions with weights \(1-\chi\) and \(\chi\), rescaling the interval to \([0,1]\). Keep the original pieces outside this interval. The interpolated curve remains a graph, so it is regular and embedded, and the vanishing endpoint derivatives of \(\chi\) make the derivatives match through order five at both joins. With the modification neighborhoods chosen above, the resulting arc is piecewise analytic and embedded in \(P\), algebraic on each piece, and agrees with the branches of \(F_0=0\) near \(t\) and \(g\).

    Use the affine parameter on each straight segment also for its adjoining graph interpolation, so that the parameterizations match in derivatives through order five at the joins. Concatenate the parameterized pieces in order from \(t\) to \(g\), then rescale the full parameter interval affinely to \([0,1]\). Denote the resulting map by \(a \colon [0,1] \to P\). The matching derivatives and nonvanishing tangent vectors make \(a\) a regular \(\mathscr{C}^5\) parameterization, and the absence of self-intersections makes it an embedding. Its finitely many semialgebraic pieces give it a computable semialgebraic description. We have \(a(0)=t\), \(a(1)=g\), and \(a((0,1])\subset\operatorname{int}P\). Near its endpoints, its image agrees with the retained branches of \(F_0=0\); at a boundary target, \(a'(0)\) points strictly into \(P\).

    \emph{Complete the arc and match defining functions.} Extend the arc \(a\) slightly beyond \(t\) and \(g\), following the local branches of \(F_0=0\) through those points. Form a second arc by moving each point of the extended arc a small fixed distance in the normal direction on one side. Join the corresponding endpoints with rounded connections to obtain a simple closed \(\mathscr{C}^4\) curve bounding a disk. Place the second arc on the side where \(F_0<0\) near \(t\) and \(g\), so that \(F_0\) increases when crossing the boundary from inside the disk to outside near either point. Our earlier choice of the endpoint covector signs makes these two requirements compatible.

    At each other real point of \(S\), take a short piece of the local branch of \(F_0=0\) and repeat this construction, again placing the enclosed disk on the side where \(F_0<0\) near that point. Choose the added arcs and loops close enough to the original pieces that all the enclosed disks are disjoint. Around each loop, take a thin \emph{tubular neighborhood}: a neighborhood formed by short normal segments extending to both sides of the loop. Choose the tubes narrow enough that they are disjoint from one another and that each avoids all the other disks.

    Here and in the next lemma, each constant is found by the same search. A positive lower bound is obtained by halving a dyadic candidate, starting from one, until a condition holds, and an upper bound by increasing an integer candidate until a condition holds. Each condition is a first-order statement about semialgebraic data over the real algebraic numbers, so it is decidable by quantifier elimination over real closed fields; see \cite[Chapter~2]{BasuPollackRoy2006}. Each search terminates because the condition holds for all sufficiently small, or sufficiently large, candidates, which we verify in each case.

    For the width of each normal strip, the search finds a positive dyadic width on which the normal coordinate map is injective and its Jacobian is nonzero throughout its domain. Construct the rounded connections in these coordinates using the transition function \(\chi\), and check that the resulting closed curve has a nonvanishing tangent and no self-intersections. All the parameterizations, coordinate domains, and their images have finite semialgebraic descriptions. For a regular \(\mathscr{C}^4\) curve, the normal coordinate map is \(\mathscr{C}^3\), and at the curve, the two columns of its Jacobian are a tangent vector and a normal vector, so they are linearly independent. The inverse function theorem and compactness then give a uniform positive width on which the Jacobian is nonzero. Injectivity also holds after shrinking the width: otherwise, distinct normal coordinates with the same image in arbitrarily thin strips would have base points converging to the same point of the embedded curve, contradicting local injectivity there. Thus sufficiently narrow strips satisfy both conditions. This is the normal-coordinate argument in \cite[Theorem~6.24]{Lee2013}, which uses only the finite differentiability available here.

    On each tubular neighborhood, use signed distance to its loop, negative in the region bounded by the loop and positive outside that region. Use the transition function \(\chi\) to construct a semialgebraic \(\mathscr{C}^2\) function that agrees with signed distance on a smaller tubular neighborhood and is constant outside the larger tubular neighborhood, that is, \(-1\) in the bounded region and \(1\) in the unbounded region. In the intervening regions, interpolate from signed distance to \(-1\) inside the loop and to \(1\) outside it, preserving the sign. Multiply these functions to obtain \(f_Z\). By the disjointness conditions above, all the other factors are equal to \(1\) near each loop, so \(f_Z\) agrees there with signed distance. Its zero set is therefore exactly the union of the loops, and its differential is nonzero on every loop.

    Near each point of \(S \cap \R^2\), the loop agrees with the local analytic branch of \(F_0=0\). In local analytic normal coordinates with the branch given by \(y=0\), we have \(f_Z=y\). Since \(F_0\) also vanishes to first order along the branch, it has the form \(F_0=yq\), where \(q\) is analytic and nonzero and our choice of signs makes \(q\) positive. Thus \(F_0/f_Z\) extends to a positive analytic function across the branch. This extension is semialgebraic and computable since it is the quotient away from the branch, and on the branch its value is the ratio of the nonzero normal derivatives of \(F_0\) and \(f_Z\).

    Now choose disjoint open balls centered at the points of \(S \cap \R^2\), small enough that the extended quotients are defined and positive on neighborhoods of their closures. For each center \(p\), use \(\chi\) to construct a radial cutoff function \(\chi_p\), i.e., a function depending only on distance from \(p\), taking values in \([0,1]\), equal to \(1\) on a smaller concentric ball, and equal to \(0\) near the boundary of the larger ball. Write \(q_p\) for the extended quotient on that ball and set the multiplier there equal to \((1-\chi_p)+\chi_p q_p\). This is positive, agrees with \(q_p\) on the smaller ball, and equals \(1\) near the boundary of the larger ball. Setting the multiplier equal to \(1\) outside the balls therefore gives a positive semialgebraic \(\mathscr{C}^2\) function on \(\R^2\). Multiply \(f_Z\) by this function to obtain \(f\). Then \(f=F_0\) near every real point of \(S\), while its zero set remains the same compact union of regular loops and contains the constructed arc from \(t\) to \(g\) in \(P\). In particular, \(f(a(s))=0\) for \(0 \le s \le 1\).
\end{proof}

\begin{figure}[tbp]
    \centering
    \begin{subfigure}[t]{.32\linewidth}
        \centering
        \begin{tikzpicture}[scale=.90,every node/.style={font=\scriptsize}]
            \LiftingPolygon
            \draw[black!45,dashed,line width=1pt] (.65,.55)--(1.05,.55)
            .. controls (1.8,.55) and (3.04,.7) .. (3.04,1.85)
            .. controls (3.04,2.2) and (3.25,2.65) .. (3.25,3);
            \LiftingArc
            \LiftingGerms
            \LiftingMarks
            \node at (1.65,1.45) {\(P\)};
            \node[anchor=east] at (2.4,.87) {\(a([0,1])\)};
            \node[align=center,text=liftpurple] at (1.1,2.12)
            {germs of \(\mathrm{C}_{F_0}(\R)\)};
            \path[use as bounding box] (-.28,-.22) rectangle (4.22,3.72);
        \end{tikzpicture}
        \caption{}
        \label{fig:lifting-assembly-a}
    \end{subfigure}\hfill
    \begin{subfigure}[t]{.32\linewidth}
        \centering
        \begin{tikzpicture}[scale=.90,every node/.style={font=\scriptsize}]
            \LiftingPolygon
            \LiftingLoop
            \LiftingAuxLoop
            \LiftingArc
            \LiftingGerms
            \LiftingMarks
            \node[anchor=west] at (3.78,1.94) {\(Z_0\)};
            \node at (1.65,1.45) {\(P\)};
            \node at (1.42,3.18) {\(Z_1\)};
            \path[use as bounding box] (-.28,-.22) rectangle (4.22,3.72);
        \end{tikzpicture}
        \caption{}
        \label{fig:lifting-assembly-b}
    \end{subfigure}\hfill
    \begin{subfigure}[t]{.32\linewidth}
        \centering
        \begin{tikzpicture}[scale=.90,every node/.style={font=\scriptsize}]
            \LiftingPolygon
            \LiftingLoop[fill=liftblue!15]
            \LiftingAuxLoop[fill=liftblue!15]
            \LiftingArc
            \LiftingGerms
            \LiftingMarks
            \node at (1.05,.32) {\(-\)};
            \node at (1,3) {\(-\)};
            \node at (2.04,2.35) {\(+\)};
            \node at (3.94,1.2) {\(+\)};
            \draw[->] (.65,.55)--(.65,1.08);
            \draw[->] (3.25,2.61)--(2.72,2.61);
            \draw[->] (.7,3)--(.18,3);
            \node at (1.65,1.45) {\(P\)};
            \path[use as bounding box] (-.28,-.22) rectangle (4.22,3.72);
        \end{tikzpicture}
        \caption{}
        \label{fig:lifting-assembly-c}
    \end{subfigure}

    \caption{Building a smooth function with the required zero-set germs. (a) Replace the dashed reference arc near \(g,t\) by the actual germs of \(\mathrm{C}_{F_0}(\R)\), shown in purple. (b) Extend past the endpoints and add a nearby second arc to form \(Z_0\); another real conjugate \(t'\) lies on a separate loop \(Z_1\). (c) Choose negative signs inside the disks and positive signs outside, so the differentials point outward. Positive local multipliers make \(f\) and \(F_0\) agree near the points of \(S \cap \R^2\). The connecting arc lies in \(P\), although the added loops need not. The drawing is schematic.}
    \label{fig:lifting-assembly}
\end{figure}
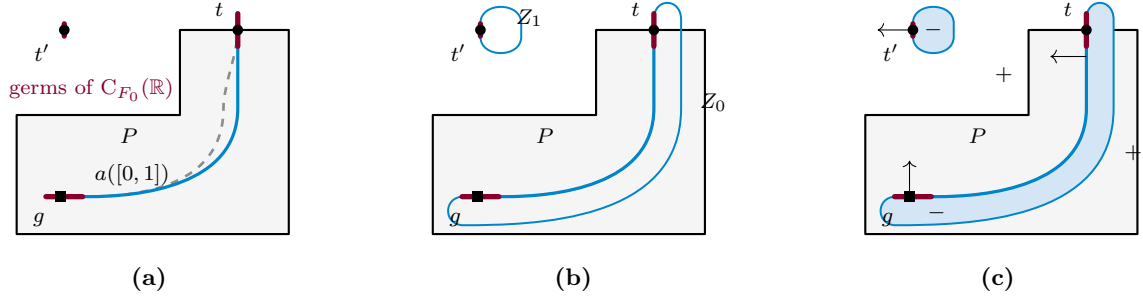

\begin{lemma}[Computable bounds for the reference function]\label{lem:computable-f-bounds}
    For the function \(f\) and arc \(a\) constructed in \cref{lem:computable-f}, we can compute a unit normal field \(\hat{n}\) along \(a\) and the normal strip map
    \[
        \Psi(s,z) \coloneqq a(s)+z\hat{n}(s).
    \]

    \emph{Strip and approximation bounds.} We can compute a positive rational \(\tau\) such that \(\Psi\) is injective with nonzero Jacobian on \([0,1] \times [-\tau,\tau]\), and an integer \(R \ge 1\) such that \([-R,R]^2\) contains the zero set of \(f\) and the strip image. We can also compute a positive rational \(r_0\) such that \(f-F_0\) vanishes on the \(r_0\)-balls about the points of \(S \cap \R^2\), and positive rational constants \(\lambda,L_a,W_a,K_a,V\) satisfying
    \[
        \begin{alignedat}{2}
            \norm{\hat{n}(s)}_2 &= 1, &\qquad \norm{\hat{n}'(s)}_2 &\le K_a, \\
            \partial_z(f \circ \Psi) &\ge \lambda, &\qquad \norm{\partial_s \Psi}_2 &\le L_a, \\
            \abs{\partial_{sz}(f \circ \Psi)} &\le W_a, &\qquad \norm{f-F_0}_{\mathscr{C}^2} &\le V.
        \end{alignedat}
    \]
    The bounds involving \(\hat{n}\) hold on \([0,1]\), those involving \(\Psi\) hold throughout the strip, and the \(\mathscr{C}^2\)-norm is taken on the bounding square.

    \emph{Bounds for keeping the arc inside \(P\).} If \(t \in \operatorname{int}P\), one can compute a positive rational \(\rho\) such that
    \[
        \operatorname{dist}(a([0,1]),\partial P) \ge \rho.
    \]
    If \(t \in \partial P\), one can instead compute positive rational numbers \(r_1,s_0,\ell_{\mathrm{ov}},\kappa,\rho\) and a closed convex cone \(C\), generated by rational vectors, such that
    \[
        \begin{alignedat}{2}
            r_1,\kappa,\rho &> 0, &\qquad 0 < s_0 < s_0+\ell_{\mathrm{ov}} &< 1, \\
            (t+C) \cap B(t,2r_1) &\subset P, &\qquad a([0,s_0+\ell_{\mathrm{ov}}]) &\subset B(t,r_1), \\
            \operatorname{dist}(a'(s),\R^2 \setminus C) &\ge \kappa, 
            &\qquad \operatorname{dist}(a([s_0,1]),\partial P) &\ge \rho.
        \end{alignedat}
    \]
    The derivative margin holds for \(0 \le s \le s_0+\ell_{\mathrm{ov}}\).
\end{lemma}

\begin{proof}
    The preceding construction gives semialgebraic descriptions over the real algebraic numbers. Differentiating its polynomial formulas and implicitly defined graph pieces gives such descriptions for the derivatives used below and the matching derivatives at the joins determine their values there. Consequently, each condition below is decidable, and the constants are found by the search described in the proof of \cref{lem:computable-f}.

    \emph{Choose the normal field and strip.} Since \(f(a(s))=0\) and the zero set of \(f\) is regular, define
    \[
        \hat{n}(s) = \frac{ \nabla f(a(s))}{\norm{\nabla f(a(s))}_2}, \quad 0 \le s \le 1.
    \]
    This is a semialgebraic \(\mathscr{C}^1\) field with \(\norm{\hat{n}(s)}_2=1\). Differentiating \(f(a(s))=0\) shows that \(\hat{n}(s)\) is perpendicular to \(a'(s)\). For \(\Psi(s,z)=a(s)+z\hat{n}(s)\), we therefore have
    \[
        \partial_z(f \circ \Psi)(s,0) = \nabla f(a(s)) \cdot \hat{n}(s) =\norm{\nabla f(a(s))}_2 > 0.
    \]
    The last expression is continuous on \([0,1]\) and has a positive minimum, so the search finds a positive dyadic \(\lambda\) with
    \[
        \norm{\nabla f(a(s))}_2 > 2\lambda, \quad 0 \le s \le 1.
    \]

    Along the arc, \(\hat{n}\) is the outward unit normal used in the tubular neighborhood construction. Note that multiplying \(f_Z\) by a positive function to obtain \(f\) does not change the gradient direction on the zero set. Restricting those normal coordinates to the arc therefore makes \(\Psi\) injective with nonzero Jacobian for sufficiently small strip width. Continuity of \(\partial_z(f \circ \Psi)\), uniformly along the compact arc, also gives
    \[
        \partial_z(f \circ \Psi)(s,z) \ge \lambda, \quad 0 \le s \le 1,\ \abs{z} \le \tau,
    \]
    for all sufficiently small \(\tau>0\). The search therefore finds a positive dyadic \(\tau\) for which injectivity, the nonzero Jacobian, and this derivative bound all hold on the closed strip.

    \emph{Choose the square and upper bounds.} The zero set of \(f\) and the strip image are compact. The search finds a positive integer \(R\) such that both sets lie in \([-R,R]^2\). To obtain \(r_0\), choose a positive dyadic radius smaller than every radius of the finitely many balls on which the preceding construction made \(f=F_0\). Thus \(f-F_0\) vanishes on the required \(r_0\)-balls.

    On the closed strip, the functions \(\norm{\partial_s \Psi}_2\) and \(\abs{\partial_{sz}(f \circ \Psi)}\) are continuous and bounded. Likewise, \(\norm{\hat{n}'}_2\) is bounded on \([0,1]\), and \(f-F_0\) and all its coordinate partial derivatives of order at most two are bounded on the square. The search therefore finds positive integers \(L_a,W_a,K_a,V\) such that
    \[
        \norm{\partial_s\Psi}_2 \le L_a,\qquad
        \abs{\partial_{sz}(f \circ \Psi)} \le W_a,\qquad
        \norm{\hat{n}'}_2 \le K_a,\qquad
        \norm{f-F_0}_{\mathscr{C}^2} \le V
    \]
    on the respective domains.

    \emph{Choose bounds to contain the arc in the gallery.} If \(t \in \operatorname{int}P\), the entire arc is a compact subset of \(\operatorname{int}P\). Its distance from \(\partial P\) is therefore positive, and the search finds a positive dyadic \(\rho\) with \(\operatorname{dist}(a([0,1]), \partial P) \ge \rho\).

    So suppose instead that \(t \in \partial P\). The vector \(a'(0)\) points strictly into the polygon. Choose two rational vectors whose rays lie on either side of this direction and span a narrow closed convex cone within the local interior sector of \(P\) at \(t\). For the cone \(C\) they generate, \(a'(0) \in \operatorname{int}C\), and there is a sufficiently small \(r_1>0\) such that
    \[
        (t+C) \cap B(t,2r_1) \subset P.
    \]
    These choices can be computed by enumerating pairs of linearly independent rational vectors together with positive dyadic radii, testing both \(a'(0) \in \operatorname{int}C\) and the displayed containment. Furthermore, the strictly inward direction and the local polygonal geometry guarantee that a suitable pair and radius occur.

    Since \(a'(0) \in \operatorname{int}C\), we can compute a positive dyadic \(\kappa\) such that
    \[
        \operatorname{dist}(a'(0),\R^2 \setminus C) > 2\kappa.
    \]
    Continuity of \(a'\), together with \(a(0)=t\), implies that for all sufficiently small \(s_0>0\),
    \[
        a([0,2s_0]) \subset B(t,r_1),\qquad
        \operatorname{dist}(a'(s),\R^2 \setminus C) \ge \kappa, \quad 0 \le s \le 2s_0.
    \]
    The search finds such a positive dyadic \(s_0<1/2\). Then set \(\ell_{\mathrm{ov}}=s_0\), whence
    \[
        0 < s_0 < s_0+\ell_{\mathrm{ov}} = 2s_0 < 1,
    \]
    and the two required estimates on \([0,s_0+\ell_{\mathrm{ov}}]\) follow.

    Finally, \(a([s_0,1])\) is a compact subset of \(\operatorname{int}P\), so its distance from \(\partial P\) is positive. The search finds a positive dyadic \(\rho\) such that
    \[
        \operatorname{dist}(a([s_0,1]),\partial P) \ge \rho.
    \]
    The initial interval \([0,2s_0]\) and the remaining interval \([s_0,1]\) overlap on \([s_0,2s_0]\). These choices leave the strip and approximation bounds unchanged and establish all the boundary-target conditions.
\end{proof}

\subsection{Polynomial approximation preserving the sight arc}
For a function \(v\) of class \(\mathscr{C}^r\) on the square \([-R,R]^2\), let \(\norm{v}_{\mathscr{C}^r}\) denote the maximum of the suprema of the absolute values of \(v\) and its coordinate partial derivatives of order at most \(r\).

\begin{lemma}[Approximation with fixed first jets]\label{lem:relative-arc}
    For the reference function just constructed, rational polynomials \(G\) in \(F_0+I(S)^2\) approximate \(f\) arbitrarily closely in \(\mathscr{C}^1\). All sufficiently close polynomials vanishing at \(S\) retain a regular arc from \(g\) to \(t\) inside \(P\).
\end{lemma}

\begin{proof}
    The following estimates prove the lemma and give a quantitative version.

    \emph{Bound the smooth factor.} Compute \(M_S \ge 1\) bounding the absolute complex coordinates of \(S\), and \(\sigma>0\) below the minimum Euclidean distance between distinct complex points of \(S\). Put
    \[
        \begin{alignedat}{2}
            a_0 &= \min(r_0,\sigma/2), &\qquad \mu &= \frac{a_0^{2N}}{2^Nk_N^2(1+M_S)^{2N}}, \\
            A_C &= (2N+1)R^{2N}A_*, &\qquad A_2 &= (2N+1)^2R^{2N}A_*, 
        \end{alignedat}
    \]
    \[
        M = \left\lceil \max\{1,V(\mu^{-1}+3A_2\mu^{-2}+2A_2^2\mu^{-3})\} \right\rceil.
    \]
    Fix a real point \(w\) outside the \(r_0\)-balls about \(S \cap \R^2\). Its distance from every point of \(S\) is at least \(a_0\). For a nonreal point, use its distinct complex conjugate to get distance at least \(\sigma/2\) from \(\R^2\). Now, for each \(p\) choose \(X-p_x\) or \(Y-p_y\) whose value at \(w\) has modulus at least \(a_0/\sqrt{2}\). Their product has value at least \((a_0/\sqrt{2})^N\) and coefficient sum at most \((1+M_S)^N\). If we express the product in the degree-\(N\) basis \(B_i^{(N)}\), then the coefficients are its free coefficients divided by the nonzero integer pivot determinant. Hence
    \[
        (a_0/\sqrt{2})^N \le k_N(1+M_S)^N\max_i \abs{B_i^{(N)}(w)},
    \]
    from which \(A_0 \ge \mu\). Coefficient sums bound the first two derivatives of \(A_0\) by \(A_C,A_2\). Twice differentiating \(u=(f-F_0)/A_0\), and using its zero extension on the balls, gives \(\norm{u}_{\mathscr{C}^2}\le M\).

    \emph{Keep the approximating arc inside \(P\).} Set \(J_a=\frac{2}{\lambda} \left( 2L_a+\frac{2W_a}{\lambda}+K_a \right)\), and choose a positive rational \(\epsilon_0\) satisfying
    \[
        0 < \epsilon_0 \le \min\{\lambda/8,\lambda\tau/4,\lambda\rho/4\}.
    \]
    For a boundary target we must also require \(\epsilon_0 \le \lambda r_1/4\) and \(\epsilon_0 \le \kappa/(2J_a)\). If \(F\) vanishes at \(S\) and \(\norm{F-f}_{\mathscr{C}^1} \le \epsilon_0\), then \(\partial_z(F \circ \Psi) \ge \lambda/2\), and its values at the two sides of the strip have opposite signs. Each normal fiber has a unique zero \(z=\phi(s)\), and by implicit differentiation
    \[
        \abs{\phi} \le 2\epsilon_0/\lambda,\qquad
        \abs{\phi'} \le \frac{2\epsilon_0}{\lambda} \left( \sqrt{2}L_a+\frac{2W_a}{\lambda} \right),\qquad
        \norm{a_F'-a'}_2 \le \epsilon_0 J_a,
    \]
    where \(a_F(s)=a(s)+\phi(s)\hat{n}(s)\). For the bound on \(\phi'\), use \(\partial_s(f \circ \Psi)(s,0)=0\). The error gradient contributes at most \(\sqrt{2}\epsilon_0 L_a\), while evaluating \(\partial_s(f \circ \Psi)\) at \(z=\phi(s)\) instead of \(z=0\) changes it by at most \(W_a \abs{\phi(s)}\). The last inequality follows by differentiating \(a_F\). The graph fixes both endpoints because \(F\) vanishes there. The nonzero Jacobian of \(\Psi\) ensures that \(a_F'\) never vanishes, and injectivity of \(\Psi\) makes \(a_F\) an embedding. The bound \(\partial_z(F \circ \Psi) \ge \lambda/2>0\) implies \(\nabla F(a_F(s)) \ne 0\) for every \(s \in [0,1]\), so the arc lies in the regular locus of \(F=0\).

    For an interior target, the distance bound keeps the whole arc in \(P\). For a boundary target the derivative estimate puts \(a_F'\) in \(C\) on \([0,s_0+\ell_{\mathrm{ov}}]\), while the value estimate puts that portion in \(B(t,2r_1)\). Integrating from the fixed endpoint puts it in \(t+C\), hence in \(P\). On \([s_0,1]\) the distance bound keeps it in \(P\). These intervals overlap, so the same approximation threshold keeps the entire arc inside \(P\).

    \begin{figure}[tbp]
        \centering
        \begin{subfigure}[t]{.49\linewidth}
            \centering
            \begin{tikzpicture}[scale=1.06,every node/.style={font=\small}]
                \LiftingPolygon
                \LiftingArc[draw=liftblue!12,line width=8pt,line cap=butt]
                \LiftingArc[draw=liftblue,densely dotted]
                \foreach \x in {.7,1.05} {
                \draw[black] (\x,.37)--(\x,.74);
                }
                \foreach \y in {1.9,2.25,2.6} {
                \draw[black] (3.08,\y)--(3.42,\y);
                }
                \draw[liftnew] (.65,.55)--(1.05,.55)
                .. controls (1.8,.55) and (3.39,.76) .. (3.31,1.85)
                .. controls (3.25,2.32) and (3.25,2.65) .. (3.25,3);
                \node[liftguard,label={[liftlabel]below left:\(g\)}]
                at (.65,.55) {};
                \node[liftpoint,label={[liftlabel]above left:\(t\)}]
                at (3.25,3) {};
                \node at (1.65,1.45) {\(P\)};
                \draw[draw=liftblue,densely dotted,line width=1.1pt]
                (.2,-.35)--(.65,-.35);
                \node[anchor=west] at (.72,-.35) {\(a([0,1])\)};
                \draw[liftnew] (.2,-.7)--(.65,-.7);
                \node[anchor=west] at (.72,-.7) {\(a_G([0,1])\)};
                \path[use as bounding box] (-.3,-.95) rectangle (4.5,4.18);
            \end{tikzpicture}
            \caption{}
            \label{fig:lifting-persistence-a}
        \end{subfigure}\hfill
        \begin{subfigure}[t]{.48\linewidth}
            \centering
            \begin{tikzpicture}[scale=1.02,every node/.style={font=\small}]
                \fill[gray!8] (0,.25) rectangle (4,3);
                \draw[black,thick] (0,3)--(4,3);
                \fill[liftblue!9] (2,3)--(.55,.25)--(3.45,.25)--cycle;
                \fill[liftblue!17] (2,3)--(1.5,.25)--(2.5,.25)--cycle;
                \draw[black!45,dashed] (2,3) circle (.85);
                \draw[black!45,dashed] (2,3) circle (1.7);
                \draw[draw=liftgreen!35,line width=10pt] (2,2.15)--(2,.4);
                \draw[draw=liftblue,densely dotted,line width=1.1pt] (2,3)--(2,.35);
                \draw[liftnew] (2,3)
                .. controls (2,2.58) and (2.13,2.43) .. (2.1,2.1)
                .. controls (2.08,1.3) and (1.92,.8) .. (2.1,.35);
                \draw[->,draw=liftblue,line width=1.1pt] (1.91,2.9)--(1.91,2.46);
                \draw[->,draw=liftorange,line width=1.1pt] (2.09,2.9)--(2.09,2.46);
                \node[liftpoint,label={above left:\(t\)}] at (2,3) {};
                \node[anchor=west] at (3.55,2.65) {\(P\)};
                \node[anchor=west] at (2.66,3.7) {\(r_1\)};
                \node[anchor=west] at (3.12,4.23) {\(2r_1\)};
                \node at (3.15,-.05) {\(C\)};
                \draw[black!65] (3.15,.1)--(3.05,.38);
                \node[anchor=east] at (1.25,1.98) {overlap};
                \draw[->,black!65] (1.27,2.02)--(1.91,2.08);
                \path[use as bounding box] (-.12,-.65) rectangle (4.35,4.78);
            \end{tikzpicture}
            \caption{}
            \label{fig:lifting-persistence-b}
        \end{subfigure}

        \caption{The two estimates that retain a sight arc under \(\mathscr{C}^1\) approximation. (a) The shaded strip follows the reference arc \(a([0,1])\). Each black normal segment contains one nearby zero of \(G\), producing the arc \(a_G([0,1])\) with the same endpoints. (b) Near \(t\), the inward tangents lie strictly inside \(C\). Away from \(t\), the green neighborhood lies inside \(P\). These regions overlap, so one choice of \(r_1\) and one approximation threshold suffice. The drawing is schematic.}
        \label{fig:lifting-persistence}
    \end{figure}
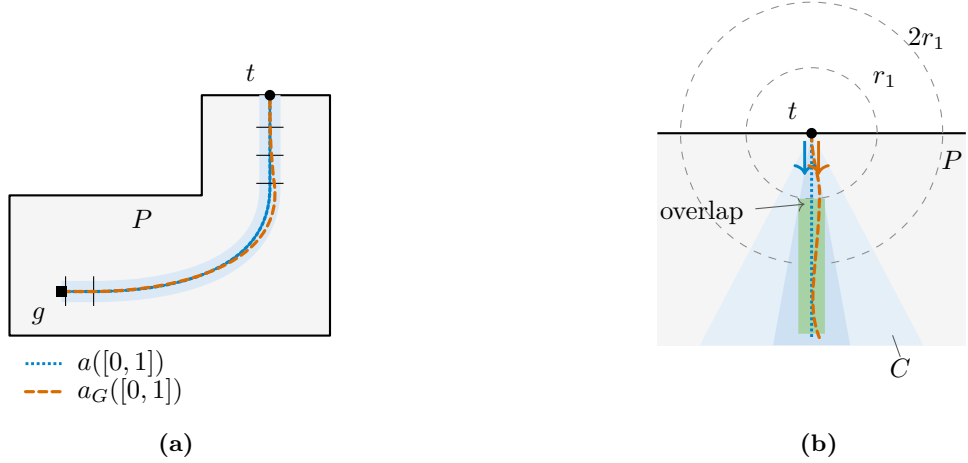

    \emph{Approximate \(u\) by a rational polynomial.} Use \(n,Q_u,d_u,H_u\) from \cref{thm:effective-lifting}. Sample \(u\) on the uniform square grid
    \[
        \left( -R+\frac{2Ri}{n},-R+\frac{2Rj}{n}\right),\quad i,j=0, \ldots, n.
    \] 
    Round each sample to the nearest multiple of \(1/Q_u\), with error at most \(1/(2Q_u)\). Since the grid points are rational and the graph of \(u\) is semialgebraic over the real algebraic numbers, the sample values are real algebraic and this rounding is computable. Rescale each coordinate interval from \([-R,R]\) to \([0,1]\). On the rescaled square, form the \emph{tensor-product Bernstein polynomial}: use products of the one-variable Bernstein basis polynomials in the two coordinates, with the rounded samples as coefficients. Express the resulting polynomial in the original coordinates \(X,Y\), and denote it by \(\beta\).

    We now prove the derivative estimates, including their scaling. For the unrounded polynomial, view the basis weights as independent binomial probabilities. In the original coordinates, the expected absolute difference between each sampled coordinate and the corresponding evaluation coordinate \(X\) or \(Y\) is at most \(R/\sqrt{n}\), whence the value error is at most \(2RM/\sqrt{n}\). Write \(x_i=-R+2Ri/n\) and \(y_j=-R+2Rj/n\). Differentiating in \(X\) gives Bernstein weights of degrees \(n-1,n\), multiplying the scaled sample differences below. For \(0 \le i<n\) and \(0 \le j \le n\), the fundamental theorem of calculus gives
    \[
        \frac{n}{2R}\bigl(u(x_{i+1},y_j)-u(x_i,y_j)\bigr)
        = \int_0^1 u_X\left(x_i+\frac{2R\xi}{n},y_j\right)\,d\xi.
    \]
    For fixed \(X\), let \(I\) be a binomial random variable with parameters \(n-1,(X+R)/(2R)\). The degree-\(n-1\) Bernstein weights in \(X\) are precisely the probabilities of the values of \(I\), so the interval in the integral is \([x_I,x_{I+1}]\). For each \(\xi \in [0,1]\), the binomial mean and variance give
    \[
        \mathbb{E}\abs{x_I+\frac{2R\xi}{n}-X}
        \le \frac{2R}{n}\left(\sqrt{(n-1)/4}+1\right)
        \le \frac{R}{\sqrt{n}}+\frac{2R}{n}.
    \]
    The independently sampled \(Y\)-coordinate has expected absolute difference at most \(R/\sqrt{n}\) from \(Y\). Applying the bounds \(\abs{u_{XX}},\abs{u_{XY}} \le M\) to the integral therefore gives derivative error
    \[
        2RM(1/\sqrt{n}+1/n) \le 4RM/\sqrt{n}.
    \]
    The estimate for the \(Y\)-derivative follows \emph{mutatis mutandis} from the argument above. Rounding adds value error at most \(1/(2Q_u)\) and error at most \(n/(2RQ_u)\) to either first derivative. The choices of \(n,Q_u\) therefore give \(\norm{\beta-u}_{\mathscr{C}^1} \le \varepsilon\).

    Expanding the factors \((R+X)^i(R-X)^{n-i}\) and their \(Y\) counterparts shows
    \[
        \deg \beta \le 2n,\qquad d_u\beta \in \Z[X,Y],\qquad \norm{d_u\beta}_1 \le H_u.
    \]
    The sum of their coefficient bounds, including both binomial sums, is \([4(1+R)^2]^n\), and the numerator of each rounded sample, written with denominator \(Q_u\), has absolute value at most \(Q_uM+1\). Set \(G = F_0+A_0\beta\). Then \(G-F_0 \in I(S)^2\), \(\deg G \le m\), and
    \[
        \norm{G-f}_{\mathscr{C}^1} \le 2A_C \varepsilon = \epsilon_0/2,\qquad D_GG \in \Z[X,Y],\qquad \norm{D_GG}_\infty \le C_G.
    \]
    This proves the approximation assertion of \cref{lem:relative-arc}. Since \(\epsilon_0\) may be chosen arbitrarily small, the same construction gives arbitrarily close \(\mathscr{C}^1\) approximations in \(F_0+I(S)^2\).
\end{proof}

\subsection{Smooth projective completion by rational perturbation}
\begin{lemma}\label{lem:bertini-completion}
    Let \(S \subset \Qbar^2\) be finite and Galois stable, with \(N=\card{S} \ge 1\). In every degree \(m \ge \max\{2,2N-1,N+1\}\), the rational linear system through \(S\) contains a nonempty Zariski-open set of smooth projective plane curves. The discriminant restricts to a nonzero polynomial of degree \(3(m-1)^2\) on this linear system.
\end{lemma}

\begin{proof}
    We consider the space of all homogeneous degree-\(m\) polynomials vanishing at the projective points \([p_x:p_y:1]\) for \(p \in S\). Working over \(\Qbar\), we determine how many independent linear equations on the coefficients are required for the curve to be singular at a specified point.

    To that end, fix a \(z \in \PP^2 \setminus S\). For each point of \(S\), choose a homogeneous linear polynomial vanishing there but not at \(z\). The product of these polynomials has degree \(N\) and vanishes throughout \(S\) by construction. Furthermore, it is nonzero at \(z\). Multiply this product by the \((m-N-1)\)-st power of a homogeneous linear polynomial nonzero at \(z\), and \emph{then} by an arbitrary homogeneous linear polynomial. This produces degree-\(m\) members of the system. In local affine coordinates, the product of the \(N\) chosen linear polynomials and the additional \((m-N-1)\)-st power is nonzero at \(z\). Consequently, varying the final linear polynomial allows us to prescribe the value and both first partial derivatives at \(z\) independently. Requiring all three quantities to vanish therefore imposes three independent linear equations on the coefficients.

    At \(p \in S\), use the first-jet interpolation construction from \cref{lem:equivariant-jets} over \(\Qbar\) which allows us to prescribe both first partial derivatives at \(p\) while keeping the values at every point of \(S\) equal to zero. The interpolants have degree at most \(2N-1 \le m\), whence homogenizing them to degree \(m\) gives members of our system. Since the value at \(p\) is already zero, singularity there imposes just two independent linear conditions.

    Let \(\mathcal{I}\) be the set of pairs \((F,z)\) for which the first jet of \(F\) vanishes at \(z\). This is a closed subset of the product of the coefficient space with \(\PP^2\). Outside \(S\), the three conditions on \(F\) and the two dimensions available to \(z\) give a dimension bound of the dimension of the coefficient space minus one. Over the finite set \(S\), the two conditions give a bound of the dimension of the coefficient space minus two. Thus \(\mathcal{I}\) has dimension strictly smaller than the coefficient space. As projection by \((F,z)\mapsto F\) cannot increase dimension, its image cannot fill that space. The image is also closed because \(\PP^2\) is projective. Its complement is therefore a nonempty Zariski-open subset of the coefficient space consisting of polynomials whose projective zero sets are smooth.

    The discriminant consequently restricts to a nonzero polynomial on this coefficient space. It is homogeneous of degree \(3(m-1)^2\), so its restriction has the same degree; see \cite[Chapter~13, \S1D, pp.~433--434]{GKZ1994}.
\end{proof}

It remains to perturb the approximation to a smooth projective curve without losing the sight arc.

\begin{proof}[Completion of the proof of \cref{thm:effective-lifting}]
    If \(g=t\), the vertical rational line through \(g\) suffices. Assume \(g \ne t\), and let \(G=F_0+A_0\beta\) be the rational polynomial constructed in the proof of \cref{lem:relative-arc}. It has degree at most \(m=2N+2n\), vanishes at every point of \(S\), and satisfies \(\norm{G-f}_{\mathscr{C}^1} \le \epsilon_0/2\). We now perturb \(G\).

    The integral polynomials \(B_1^{(m)},\ldots,B_k^{(m)}\) constructed in \cref{sec:effective-preparation} form a basis over \(\Q\) of the polynomials of degree at most \(m\) vanishing at \(S\). Here \(k=J_m-N\), and each basis polynomial has absolute coefficient sum at most \(J_mb_m\). With \(\delta=3(m-1)^2\) and \(Q_{\Gamma}\) as in \cref{thm:effective-lifting}, write \(\Gamma=\{0,\ldots,\delta\}^k\) for the parameter grid and consider the polynomials
    \[
        G_{\mathbf{a}} = G+\frac{1}{Q_{\Gamma}} \sum_{i=1}^k a_i B_i^{(m)},
        \quad \mathbf{a} = (a_1,\ldots,a_k),\quad a_i \in \{0,1,\ldots,\delta\}.
    \]
    Each of these polynomials still vanishes at \(S\), since every term added to \(G\) does.

    We first check that all these perturbations are small enough to preserve the connecting arc. Since \(R \ge 1\), the degree and coefficient bounds give \(\norm{B_i^{(m)}}_{\mathscr{C}^1} \le mJ_mb_mR^m\) on the square \([-R,R]^2\). The choice of \(Q_{\Gamma}\) therefore gives
    \[
        \norm{G_{\mathbf{a}}-G}_{\mathscr{C}^1} \le \frac{\delta k m J_mb_mR^m}{Q_{\Gamma}} < \frac{\epsilon_0}{4}.
    \]
    Together with the approximation bound for \(G\), this implies
    \[
        \norm{G_{\mathbf{a}}-f}_{\mathscr{C}^1} \le \norm{G_{\mathbf{a}}-G}_{\mathscr{C}^1}
           +\norm{G-f}_{\mathscr{C}^1} < \frac{3\epsilon_0}{4} < \epsilon_0.
    \]
    Thus \cref{lem:relative-arc} applies to every candidate: its zero set contains a regular arc from \(g\) to \(t\) inside \(P\).

    It remains to choose a candidate whose degree-\(m\) homogenization defines a smooth projective curve. If the parameters \(a_1,\ldots,a_k\) are allowed arbitrary values in \(\Qbar\), the polynomials \(G_{\mathbf{a}}\) run through the entire space of degree-at-most-\(m\) polynomials vanishing at \(S\). Indeed, the \(B_i^{(m)}\) form a basis of that space, and \(G\) belongs to it. Their degree-\(m\) homogenizations therefore run through the linear system in \cref{lem:bertini-completion}. By that lemma, the discriminant, expressed as a polynomial in \(a_1, \ldots, a_k\), is not identically zero and has degree at most \(\delta\).

    This discriminant cannot vanish at every point of \(\Gamma\); this is a special case of \cite[Lemma~2.1]{Alon1999}, whose short proof we include. We argue by induction on \(k\), the case \(k=1\) being the root bound for a nonzero univariate polynomial of degree at most \(\delta\). Write the discriminant as \(\sum_{j \le \delta} c_j(a_2,\ldots,a_k)\,a_1^j\). If it vanished on \(\Gamma\), then for each \((a_2,\ldots,a_k) \in \{0,\ldots,\delta\}^{k-1}\) the polynomial in \(a_1\) would have \(\delta+1\) roots, so every \(c_j\) would vanish on \(\{0,\ldots,\delta\}^{k-1}\). Each \(c_j\) has degree at most \(\delta\) in every variable, so by induction each \(c_j\) is identically zero, and so is the discriminant. Hence at least one candidate has nonzero discriminant. We can identify such a candidate by computing the resultant of the three homogeneous partial derivatives of each candidate's homogenization until a nonzero value is obtained.

    So choose such a \(G_{\mathbf{a}}\). Its projective zero set is smooth and therefore absolutely irreducible: two distinct positive-degree components would meet by B\'ezout's theorem and give a singular point, while a repeated component would also make the curve singular. Since \(m \ge 2\), the irreducible homogeneous equation is not divisible by \(Z\). Consequently, \(G_{\mathbf{a}}\) is absolutely irreducible and its projective closure is the smooth curve just selected.

    Finally, multiply \(G_{\mathbf{a}}\) by \(D_GQ_{\Gamma}\) to obtain an integral polynomial, and divide its coefficients by their greatest common divisor to obtain a primitive polynomial \(F\). These operations do not change the zero set. Using \(\norm{D_GG}_\infty \le C_G\) and \(\norm{B_i^{(m)}}_\infty \le J_mb_m\), we obtain
    \[
        \norm{F}_\infty \le \norm{D_GQ_{\Gamma}G_{\mathbf{a}}}_\infty \le Q_{\Gamma}C_G+D_G\delta kJ_mb_m = H.
    \]
    The polynomial \(F\) has degree at most \(m\) and retains the regular connecting arc inside \(P\). This proves \cref{thm:effective-lifting}, and hence \cref{thm:lifting}.
\end{proof}

\begingroup
\small
\apptocmd{\thebibliography}{\setlength{\itemsep}{2pt}}{}{}

\endgroup
\end{document}